\documentclass[twoside]{article}
\usepackage[accepted]{aistats2025}
\usepackage[round]{natbib}

\usepackage{mathtools}
\usepackage{amsmath,amsfonts,amssymb,dsfont,bm,amsthm}
\usepackage{csquotes}
\usepackage{xcolor}
\usepackage{float}
\usepackage{subcaption}
\usepackage{graphicx}
\usepackage{txfonts}
\usepackage{hyperref}
\usepackage{thm-restate} 
\usepackage{soul}

\theoremstyle{definition}

\newtheorem{lemma}{Lemma}

\newtheorem{definition}{Definition}

\newtheorem{problem}{Problem}
\newtheorem{assumption}{Assumption}

\usepackage{slashed}
\newcommand{\ind}{\perp\!\!\!\perp}

\newcommand{\NC}{\mathcal{NC}}
\newcommand{\CD}{\mathcal{CD}}
\newcommand{\Gf}{\G^f}
\newcommand{\Gs}{\G^s}
\newcommand{\Anc}{\text{Anc}}
\newcommand{\Desc}{\text{Desc}}
\newcommand{\G}{\mathcal{G}}
\newcommand{\C}{\mathcal{C}}
\newcommand{\Do}{\text{do}}

\usepackage{tikz}
\usetikzlibrary{babel}
\usetikzlibrary{arrows}
\usetikzlibrary{decorations}
\usetikzlibrary{decorations.pathmorphing, decorations.pathreplacing, decorations.shapes}

\newcommand{\leftselfloop}{
    \begin{tikzpicture}[baseline=(X.base)]
        \node (X) at (0,0) {$\phantom{X}$};
        \path (X) edge [->, loop left, looseness=4, in=155, out=205] node {} (X);
    \end{tikzpicture}
    \!\!\!\!\!\!\!\!
}

\usepackage[vlined]{algorithm2e}
\RestyleAlgo{ruled}
\SetKwComment{Comment}{/* }{ */}

\definecolor{CentraleRed}{rgb}{0.558,0.09, 0.18}
\definecolor{CentraleGray}{rgb}{0.67,0.67, 0.67}
\definecolor{CentraleBlue}{RGB}{0,67,89}
\definecolor{CentraleGrayEq}{RGB}{66,66,66}

\newcommand{\hiddensection}[1]{
  \addtocontents{toc}{\protect\setcounter{tocdepth}{-1}} 
  \section{#1}
  \addtocontents{toc}{\protect\setcounter{tocdepth}{2}} 
}

\newcommand{\hiddensubsection}[1]{
  \addtocontents{toc}{\protect\setcounter{tocdepth}{-1}} 
  \subsection{#1}
  \addtocontents{toc}{\protect\setcounter{tocdepth}{2}} 
}

\newcommand{\hiddensubsubsection}[1]{
  \addtocontents{toc}{\protect\setcounter{tocdepth}{-1}} 
  \subsubsection{#1}
  \addtocontents{toc}{\protect\setcounter{tocdepth}{2}} 
}

\begin{document}

%

%

\twocolumn[

\aistatstitle{Identifiability by common backdoor in summary causal graphs of time series}

\aistatsauthor{Clément Yvernes$^1$ \And Charles K. Assaad$^2$  \And  Emilie Devijver$^1$ \And Eric Gaussier$^1$ }

\aistatsaddress{$^1$ Univ Grenoble Alpes, CNRS, Grenoble INP, LIG, F38000, Grenoble, France  \\  
$^2$ Sorbonne Université, INSERM, Institut Pierre Louis d’Epidémiologie et de Santé Publique, F75012, Paris, France 
 } ]

\begin{abstract}
  The identifiability problem for interventions aims at assessing whether the total effect of some given interventions can be written with a do-free formula, and thus be computed from observational data only. 
  We study this problem, considering multiple interventions and multiple effects, in the context of time series when only abstractions of the true causal graph in the form of summary causal graphs are available. We focus in this study on identifiability by a common backdoor set, and establish, for time series with and without \textit{consistency throughout time, conditions} under which such a set exists. We also provide algorithms of limited complexity 
  to decide whether the problem is identifiable or not.
\end{abstract}

\hiddensection{Introduction}

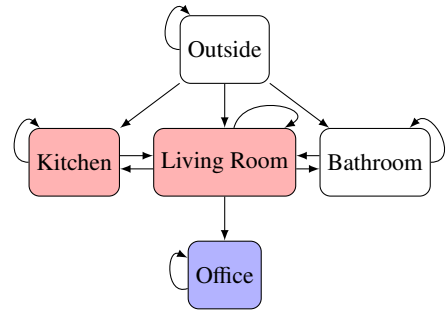
\begin{figure}
\centering
 \begin{tikzpicture}[{black, rectangle, draw, inner sep=0.1cm}]
\tikzset{nodes={draw,rounded corners},minimum height=0.9cm,minimum width=0.9cm, font=\footnotesize}
	
	\node (Outside) at (0,0) {Outside};
	\node[fill=red!30] (LivingRoom) at (0, -1.5) {Living Room};
	\node[fill=red!30] (Kitchen) at (-2,-1.5) {Kitchen};
	\node (Bathroom) at (2,-1.5) {Bathroom};
	\node[fill=blue!30] (Office) at (0,-3) {Office};

	\draw[->,>=latex] (Outside) -- (LivingRoom);
	\draw[->,>=latex] (Outside) -- (Kitchen);
	\draw[->,>=latex] (Outside) -- (Bathroom);
	\draw[->,>=latex] (LivingRoom) -- (Office);
	
	 \begin{scope}[transform canvas={yshift=-.25em}]
         \draw [->,>=latex] (LivingRoom) -- (Kitchen);
         \end{scope}
         \begin{scope}[transform canvas={yshift=.25em}]
         \draw [<-,>=latex] (LivingRoom) -- (Kitchen);
     \end{scope}

     \begin{scope}[transform canvas={yshift=-.25em}]
         \draw [->,>=latex] (LivingRoom) -- (Bathroom);
         \end{scope}
         \begin{scope}[transform canvas={yshift=.25em}]
         \draw [<-,>=latex] (LivingRoom) -- (Bathroom);
     \end{scope}

     \draw[->,>=latex] (Outside) to [out=180,in=135, looseness=2] (Outside);
	\draw[->,>=latex] (LivingRoom) to [out=75,in=30, looseness=2] (LivingRoom);
	\draw[->,>=latex] (Kitchen) to [out=180,in=135, looseness=2] (Kitchen);
	\draw[->,>=latex] (Bathroom) to [out=0,in=45, looseness=2] (Bathroom);
 	\draw[->,>=latex] (Office) to [out=200,in=155, looseness=2] (Office);
	\end{tikzpicture}
    \caption{Summary causal graphs from Thermoregulation \citep{assaad_identifiability_2024,Peters_2013}. Red and blue vertices represent the total effects of interest where the red vertex represents the intervention and the blue vertex represents the response.}
 \label{fig:real_thermoregulation}
\end{figure}

Time series are generated when observations, such as sensor data, are collected over time; as such they are present in various
forms in many different domains. Central to the study of time series in many applications is the estimation of the total effect of different interventions. For example, one may wonder what is the effect on the temperature in a room in a house after having set the temperature on neighbouring rooms, as exemplified in Figure~\ref{fig:real_thermoregulation} which represents five time series of the temperatures, regularly observed through a discrete time grid, of five different variables corresponding to ''outside'', ''kitchen'', ''living room'', ''bathroom'', and ''office''. In this particular example, one is interested in knowing the effect on the temperature in the office (effect, in blue) {at 7pm } after having changed the temperature in the kitchen and the living room  {at 5pm } (interventions, in red).

Knowing the effect of interventions is key to understanding the effect of a treatment in medicine or the effect of a maintenance operation in IT monitoring systems for example. When one cannot perform interventions in practice, for example when these interventions may endanger people's life or when they may disrupt a critical process or be too costly, estimating the effect of interventions usually requires (a) rewriting, if possible, the probability of observing a particular effect at a given time knowing past interventions in terms of probabilities which only involve observed variables, the so-called \textit{do-free formula}, and then (b) estimating these probabilities from observational data only. We focus in this paper on step (a) which is usually referred to as the \textit{identifiability problem}.

The identifiability problem has been studied in various causal graphs in time series. Among such graphs, the summary causal graph represents causal relations within and between time series without any time information, as in Figure~\ref{fig:real_thermoregulation}. As such, they are simpler to define and are a popular tool for experts willing to define causal relations between time series. They are however more challenging to deal with as a lot of crucial information, related to temporal relations, is lost, and, to our knowledge, no general result on identifiability in summary causal graphs has been derived so far. In this study, we provide such a general result for \textit{the identifiability problem by common backdoor, in the context of multiple interventions on multiple effects, with and without assuming \textit{consistency throughout time}}.

The remainder of the paper is structured as follows: related work is discussed in Section \ref{sec:sota}; Section \ref{sec:notions} introduces the main notions while Section \ref{section:IBC} presents our main result regarding identifiability by common backdoor without assuming \textit{consistency throughout time}; Section \ref{sec:Consistency_Time} presents a similar result when \textit{consistency throughout time} holds; lastly, Section~\ref{sec:conclusion} concludes the paper. All proofs are provided in the Supplementary Material.

\hiddensection{State of the art}
\label{sec:sota}
Pearl's framework allows the determination of causal effects from observational data without requiring interventions. This is achieved through graphical criteria such as the backdoor and front-door criteria, which enable the unique derivation of causal effects from a given causal model \citep{Pearl_book2000}.
If the causal directed acyclic graph (DAG) is known, and all variables are observed (causal sufficiency),  the backdoor criterion is complete for total effects in a univariate setting since all parents are observed \citep{Pearl_1995}. \cite{Shpitser_2010} have shown that, however, in case of identifiability, it does not identify all possible adjustment sets. When the backdoor is not complete, \textit{e.g.}, with hidden confounders or multiple interventions, one may relate to the do-calculus \citep{Pearl_1995} and the associated ID algorithm, which is complete \cite{Shpitser_2010}. 

When considering a class of potential causal graphs, the problem is more complex. 
\citet{Maathuis_2013,Perkovic_2016} provide conditions for identifiability based on completed partially directed acyclic graph (CPDAG, \cite{meek1995}), which are necessary and sufficient  for the monovariate case but only sufficient for the multivariate case.
\cite{perkovic_identifying_2020} has developed a necessary and sufficient criterion for causal identification within maximally oriented partially directed acyclic graphs (MPDAGs, \cite{meek1995}) under causal sufficiency, which encompass DAGs, CPDAGs and CPDAGs with   background knowledge. 
When considering latent confounding, 
\cite{NEURIPS2022_17a9ab41,Wang_2023} provide sufficient conditions using partial ancestral graph (PAGs, \cite{Spirtes_2000}).
Another generalization consists of the abstraction of causal graphs, where variables are defined at a higher level of granularity.  Cluster DAGs \citep{Anand_2023} encode  partially understood causal relationships between variables in different  known clusters, representing a group of variables among which causal relationships are not understood or specified. \citet{Anand_2023} provide necessary and sufficient conditions  (considering an extension of the do-calculus) in Cluster DAGs for identifying total effects.

When considering time series, several abstractions are available, such as extended summary causal graphs and summary causal graphs (ECGs and SCGs, respectively). 
\citet{Eichler_2007} provide sufficient conditions for identifiability of the total effect on graphs based on time series that can be applied to SCGs assuming  no instantaneous relations. 
In a more general case, with possible instantaneous relations,
\cite{Assaad_2023} demonstrated that the total effect is always identifiable using summary causal graphs (SCGs), under the assumptions of causal sufficiency and the absence of cycles larger than one in the SCG. 
\cite{assaad2024identifiabilitytotaleffectssummary} provides sufficient conditions based on the front-door criterion when causal sufficiency is not satisfied.
In parallel, under causal sufficiency, \cite{Ferreira_2024} present conditions for identifying direct effect using SCGs assuming a linear model.
More recently, \cite{assaad_identifiability_2024} shown under causal sufficiency that the total effect is always identifiable in ECGs and sufficient conditions for identifiability with common backdoor in SCGs, assuming causal sufficiency, \textit{consistency throughout time} and a single intervention. This contrasts with our setting which considers both multiple interventions and multiple effects, without and with \textit{consistency throughout time}.

\hiddensection{Context}
\label{sec:notions}


For a graph $\G$, if $X\rightarrow Y$, then $X$ is a \emph{parent} of $Y$. A \emph{path} is a sequence of distinct vertices such that there exists an arrow in $\G$ connecting each element to its successor. A \emph{directed path} is a path in which all arrows are pointing towards the last vertex. If there is a directed path from $X$ to $Y$, which we denote by $X \rightsquigarrow Y$ or $Y \leftsquigarrow X$ when there is at least one arrow in the path, then $X$ is an \emph{ancestor} of $Y$, and $Y$ is a \emph{descendant} of $X$. The sets of parents, ancestors and descendants of $X$ in $\G$ are denoted by $\text{Pa}(X,\G)$, $\text{Anc}(X,\G)$ and $\Desc(X,\G)$ respectively. We denote by $\Desc^>(X,\G)$ the set of descendants of $X$ for which there is a non-empty directed path linking $X$ to them. We denote by $\Anc^>(X,\G)$ the set of ancestors of $X$ for which there is a non-empty directed path linking them to $X$. By a slight abuse of notation, we denote $\G \backslash \{Y\}$ the subgraph of $\G$ obtained by removing the vertex $Y$ and its corresponding edges. Furthermore, the mutilated graphs $\G_{\overline{Y}}$, $\G_{\underline{Y}}$ and $\G_{\overline{Y} \underline{Y}}$ represent the graphs obtained by respectively removing from $\G$ all incoming, all outgoing and both incoming and outgoing edges of $Y$. The \emph{skeleton} of $\G$ is the undirected graph given by forgetting all arrow orientations in $\G$. 
The subgraph $\G_{\mid S}$ of the graph $\G$  induced by the set of vertices $S$ consists of the set of vertices in $S$ and all the edges of $\G$ whose endpoints are both in $S$.

A \emph{backdoor path} between $X$ and $Y$ is a path between $X$ and $Y$ in which the first arrow is pointing to $X$. A \emph{directed cycle} is a circular list of distinct vertices in which each vertex is a parent of its successor. We denote by $Cycles(X,\G)$ the set of all directed cycles containing $X$ in $\G$, and by $Cycles^>(X,\G)$ the subset of $Cycles(X,\G)$ with at least 2 vertices (i.e., excluding self-loops). If a path $\pi$ contains $X_i \rightarrow X_j \leftarrow X_k$ as a subpath, then $X_j$ is a \emph{collider} on $\pi$. A path $\pi$ is \emph{blocked} by a subset of vertices $\mathcal{Z}$ if a non-collider in $\pi$ belongs to $\mathcal{Z}$ or if $\pi$ contains a collider of which no descendant belongs to $\mathcal{Z}$. Otherwise, $\mathcal{Z}$ \emph{d-connects} $\pi$.  Given an ordered pair of variables $(X, Y)$ in a DAG, a set of variables $\mathcal{Z}$ satisfies the \emph{standard backdoor criterion} relative to $(X, Y)$ if no vertex in $\mathcal{Z}$ is a descendant of $X$, and $\mathcal{Z}$ blocks every backdoor path between $X$ and $Y$. The following definition presents this criterion in the multivariate case.
 

\begin{definition}[Backdoor Criterion (Multivariate)]
\label{def:critere_backdoor_multivarie}
    Let $\G$ be a causal DAG, $\{X^i\}_i$ and $\{Y^j\}_j$ two distinct sets of variables. Let $\mathcal{Z}$ be a set of variables distinct from $\{X^i\}_i$ and $\{Y^j\}_j$. $\mathcal{Z}$ is said to satisfy the backdoor criterion relative to $\{X^i\}_i$ and $\{Y^j\}_j$ if: 
\begin{enumerate}
    \item No element of $\mathcal{Z}$ is a descendant of any $X^i$ in $\G$; and  \label{def:critere_backdoor:1}
    \item For all pair $(X^i, Y^j)$, $\mathcal{Z}$ blocks all backdoor paths between $X^i$ and $Y^j$ in $\G_{ \underline{(X^{i'})_{i' \neq i}}}$.\footnote{Alternatively, one can equivalently write: $\mathcal{Z} \cup \{X^{j}\}_{j \neq i}$ blocks all backdoor paths between $X^i$ and $Y^j$ in $\G_{ \overline{(X^{i'})_{j' \neq i}}}$, or $\mathcal{Z}$ blocks all backdoor paths between $X^i$ and $Y^j$ in $\G_{ \overline{(X^{i'})_{j' \neq i}} \underline{(X^{i'})_{i' \neq i}}}$.} \label{def:critere_backdoor:2}
\end{enumerate}
\end{definition}
This definition is the multivariate counterpart of the univariate backdoor criterion stated in \cite[Theorem 5.3.2] {Pearl_book2000}. One furthermore has:
\begin{restatable}{lemma}{mylemmaobservationsmultiples}{}
\label{lemma:observations_multiples}
Let $\mathcal{Z}$ be a set of temporal vertices. Then the following statements are equivalent:
\begin{enumerate}
    \item $\mathcal{Z}$ satisfies the backdoor criterion relative to $\{X^i\}_i$ and $\{Y^j\}_j$.
    \item For all $j$, $\mathcal{Z}$ satisfies the backdoor criterion relative to $\{X^i\}_i$ and $\{Y^j\}$.
\end{enumerate}
\end{restatable}



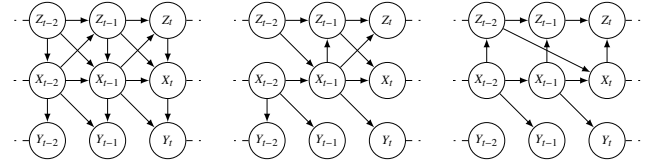
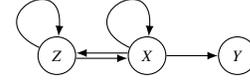
\begin{figure}[t]
\centering

 \begin{subfigure}{0.5\textwidth}
 \centering
 \scalebox{0.80}{
\begin{tikzpicture}[scale=1, transform shape, black, circle, draw, inner sep=0]
 \tikzset{nodes={draw,rounded corners},minimum height=0.6cm,minimum width=0.6cm, font = \tiny}
 \tikzset{latent/.append style={fill=gray!60}}
 
 \node (X) at (1,1) {$X_t$};
 \node (Y) at (1,0) {$Y_t$};
 \node (X-1) at (0,1) {$X_{t-1}$};
 \node (Y-1) at (0,0) {$Y_{t-1}$};
 \node (X-2) at (-1,1) {$X_{t-2}$};
 \node (Y-2) at (-1,0) {$Y_{t-2}$};
 \node (V) at (1,2) {$Z_t$};
 \node (V-1) at (0,2) {$Z_{t-1}$};
 \node (V-2) at (-1,2) {$Z_{t-2}$};
 \draw[->,>=latex] (X-1) to (Y);
 \draw[->,>=latex] (X-2) -- (Y-1);
 \draw[->,>=latex] (X) to  (Y);
 \draw[->,>=latex] (X-1) to (Y-1);
 \draw[->,>=latex] (X-2) to  (Y-2);
 \draw[->,>=latex] (X-2) -- (X-1);
 \draw[->,>=latex] (X-1) -- (X);
 \draw[->,>=latex] (V-2) -- (V-1);
 \draw[->,>=latex] (V-1) -- (V);
 \draw[->,>=latex] (V-2) -- (X-1);
 \draw[->,>=latex] (V-1) -- (X);
 \draw[->,>=latex] (X-2) -- (V-1);
 \draw[->,>=latex] (X-1) -- (V);
 \draw[->,>=latex] (V-2) to  (X-2);
 \draw[->,>=latex] (V-1) to  (X-1);
 \draw[->,>=latex] (V) to  (X);
`
\draw [dashed,>=latex] (V-2) to[left] (-1.55,2);
\draw [dashed,>=latex] (X-2) to[left] (-1.56,1);
 \draw [dashed,>=latex] (Y-2) to[left] (-1.55,0);
 \draw [dashed,>=latex] (V) to[right] (1.55,2);
\draw [dashed,>=latex] (X) to[right] (1.55,1);
\draw [dashed,>=latex] (Y) to[right] (1.55, 0);

\node (a) at (4.65,1) {$X_t$};
\node (b) at (4.65,0) {$Y_t$};
\node (a-1) at (3.65,1) {$X_{t-1}$};
\node (b-1) at (3.65,0) {$Y_{t-1}$};
\node (a-2) at (2.65,1) {$X_{t-2}$};
\node (b-2) at (2.65,0) {$Y_{t-2}$};
\node (c) at (4.65,2) {$Z_t$};
\node (c-1) at (3.65,2) {$Z_{t-1}$};
\node (c-2) at (2.65,2) {$Z_{t-2}$};
 \draw[->,>=latex] (a-1) to (b);
 \draw[->,>=latex] (a-2) -- (b-1);
 \draw[->,>=latex] (a-2) to  (b-2);
 \draw[->,>=latex] (c-2) -- (a-1);
 \draw[->,>=latex] (c-1) -- (a);
 \draw[->,>=latex] (a-1) -- (c);
 \draw[<-,>=latex] (c-1) to  (a-1);
 \draw[->,>=latex] (a-2) -- (a-1);
 \draw[->,>=latex] (a-1) -- (a);
 \draw[->,>=latex] (c-2) -- (c-1);
 \draw[->,>=latex] (c-1) -- (c);
 
\draw [dashed,>=latex] (c-2) to[left] (2.1,2);
\draw [dashed,>=latex] (a-2) to[left] (2.09, 1);
 \draw [dashed,>=latex] (b-2) to[left] (2.1, 0);
\draw [dashed,>=latex] (c) to[right] (5.2,2);
\draw [dashed,>=latex] (a) to[right] (5.2, 1);
\draw [dashed,>=latex] (b) to[right] (5.2, 0);

\node (d) at (8.3,1) {$X_t$};
\node (e) at (8.3,0) {$Y_t$};
\node (d-1) at (7.3,1) {$X_{t-1}$};
\node (e-1) at (7.3,0) {$Y_{t-1}$};
\node (d-2) at (6.3,1) {$X_{t-2}$};
\node (e-2) at (6.3,0) {$Y_{t-2}$};
\node (f) at (8.3,2) {$Z_t$};
\node (f-1) at (7.3,2) {$Z_{t-1}$};
\node (f-2) at (6.3,2) {$Z_{t-2}$};
 \draw[->,>=latex] (d-1) to (e);
 \draw[->,>=latex] (d-2) -- (e-1);
 \draw[->,>=latex] (f-2) -- (d);
 \draw[<-,>=latex] (f-2) to  (d-2);
 \draw[<-,>=latex] (f-1) to  (d-1);
 \draw[<-,>=latex] (f) to  (d);
 \draw[->,>=latex] (d-2) -- (d-1);
 \draw[->,>=latex] (d-1) -- (d);
 \draw[->,>=latex] (f-2) -- (f-1);
 \draw[->,>=latex] (f-1) -- (f);
\draw [dashed,>=latex] (f-2) to[left] (5.75,2);
\draw [dashed,>=latex] (d-2) to[left] (5.74, 1);
 \draw [dashed,>=latex] (e-2) to[left] (5.75, 0);
\draw [dashed,>=latex] (f) to[right] (8.85,2);
\draw [dashed,>=latex] (d) to[right] (8.85, 1);
\draw [dashed,>=latex] (e) to[right] (8.85, 0);

\end{tikzpicture}}
 	\caption{\centering Three FTCGs, $\mathcal{G}_{1}^f$, $\mathcal{G}_{2}^f$ and $\mathcal{G}_{3}^f$.}
 \label{fig:example_FTCG}
 
 \end{subfigure}

  \begin{subfigure}{0.5\textwidth}
 \centering
\begin{tikzpicture}[scale=1, transform shape, black, circle, draw, inner sep=0]
	\tikzset{nodes={draw,rounded corners},minimum height=0.5cm,minimum width=0.5cm, font = \tiny}	
	\tikzset{anomalous/.append style={fill=easyorange}}
	\tikzset{rc/.append style={fill=easyorange}}
 	\node (Z) at (-1.2,0) {$Z$};
	\node (X) at (0,0) {$X$} ;
	\node (Y) at (1.2,0) {$Y$};
 \draw[->,>=latex] (X) -- (Y);
 \begin{scope}[transform canvas={xshift=0, yshift=.1em}]
 \draw [->,>=latex,] (X) -- (Z);
 \end{scope}
 \begin{scope}[transform canvas={xshift=.0, yshift=-.1em}]
 \draw [<-,>=latex,] (X) -- (Z);
 \end{scope}

	\draw[->,>=latex] (X) to [out=155,in=85, looseness=8] (X);
	\draw[->,>=latex] (Z) to [out=155,in=85, looseness=8] (Z);

 \end{tikzpicture}
 \caption{\centering The SCG $\mathcal{G}^s$, reduced from any FTCG in (a).}
 \label{fig:example_SCG}
 \end{subfigure}
 \hfill
 
 \caption{Illustration: (a) three FTCGs; (b) the SCG which can be derived from any FTCG in (a).}
 \label{fig:example_CG}
\end{figure}

%
\hiddensubsection{Causal graphs in time series}
Consider $\mathcal{V}$ a set of $p$ observational time series and $\mathcal{V}^f=\{\mathcal{V}_{t} | t \in \mathbb{Z}\}$ the set of temporal instances of $\mathcal{V}$, where $\mathcal{V}_{t}$ corresponds to the variables of the time series at time $t$. 
We suppose that the time series are generated from an \emph{unknown} dynamic structural causal model (DSCM, \citet{DSCM}), an extension of structural causal models (SCM, \citet{Pearl_book2000}) to time series. This DSCM defines a full-time causal graph (FTCG, see below) which we call the \emph{true} FTCG and a joint distribution $P$ over its vertices which we call the \emph{true} probability distribution, which is compatible with, or Markov relative to, the true FTCG by construction.
The graph that is used to qualitatively represent causal relations described in a DSCM is known as the full-time causal graph (FTCG). 
\begin{definition}[Full-time causal graph (FTCG), Figure~\ref{fig:example_FTCG}]
Let $\mathcal{V}$ be a set of $p$ observational time series and $\mathcal{V}^f=\{\mathcal{V}_{t-\ell} | \ell \in \mathbb{Z}\}$. The \emph{full-time causal graph} (FTCG) $\G^f=(\mathcal{V}^f, \mathcal{E}^f)$ representing a given DSCM is defined by: 
$X_{t-\gamma} \rightarrow Y_t \in \mathcal{E}^f$ if and only if $X$ causes $Y$ at time $t$ with a time lag of $\gamma>0$ if $X=Y$ and with a time lag of $\gamma \geq 0$ for $X \neq Y$. 
\end{definition} 

As common in causality studies on time series, we consider in the remainder acyclic FTCGs with potential self-causes, \textit{i.e.}, the fact that, for any time series $X$, $X_{t-\ell} \, (\ell \in \mathbb{N}^{*})$ may cause $X_t$. Note that acyclicity is guaranteed for relations between variables at different time stamps and that self-causes are present in most time series. We furthermore assume causal sufficiency:
\begin{assumption}[Causal sufficiency] \label{ass:cs}
    There is no hidden common cause between any two observed variables.
\end{assumption}

Experts are used to working with abstractions of causal graphs which summarize the information into a smaller graph that is interpretable, often with the omission of precise temporal information. We consider in this study a known causal abstraction for time series, namely \textit{summary causal graphs}\footnote{Other abstract graphs exist, as extended summary causal graphs; the identifiability problem in such graphs is relatively straightforward and has been studied in \cite{assaad_identifiability_2024}.}. A summary causal graph \citep{Peters_2013, Meng_2020} represents causal relationships among time series, regardless of the time delay between the cause and its effect. 
\begin{definition}[Summary causal graph (SCG), Figure~\ref{fig:example_SCG}]
	\label{Summary_G}
	Let $\G^f=(\mathcal{V}^f,\mathcal{E}^f)$ be an FTCG built from the set of time series $\mathcal{V}$. The  \emph{summary causal graph} (SCG) $\G^s=(\mathcal{V}^s, \mathcal{E}^s)$ associated to $\G^f$ is such that:
     \begin{itemize}
         \item $\mathcal{V}^s$ corresponds to the set of time series $\mathcal{V}$,
         \item $X \rightarrow Y \in \mathcal{E}^s$ if and only if there exists at least one timepoint $t$, and one temporal lag $0\leq \gamma$ such that $X_{t-\gamma} \rightarrow Y_t \in \mathcal{E}^f$.
     \end{itemize} 
In that case, we say that $\G^s$ is \emph{reduced from} $\G^f$. 
\end{definition}
SCGs may include directed cycles and even self-loops. For example, the three FTCGs in Figure \ref{fig:example_FTCG} are acyclic, while the SCG in Figure \ref{fig:example_SCG} has a cycle. 
We use the notation $X \rightleftarrows Y$ to indicate situations where there exist time instants in which $X$ causes $Y$ and $Y$ causes $X$. It is furthermore worth noting that if there is a single SCG reduced from a given FTCG, different FTCGs, with possibly different orientations and skeletons, can yield the same SCG. 
For example, the SCG in Figure \ref{fig:example_SCG}  can be reduced from any FTCG in Figure \ref{fig:example_FTCG}, even though they may have different skeletons (for example, $\G_{1}^f$ and $\G_{3}^f$) and different orientations (for example, $\G_{1}^f$ and $\G_{2}^f$).
In the remainder, for a given SCG $\Gs$, we call any FTCG from which $\Gf$ can be reduced as a \textit{candidate FTCG} for $\Gs$. For example, in Figure \ref{fig:example_CG}, $\G_{1}^f$, $\G_{2}^f$ and $\G_{3}^f$ are all candidate FTCGs for $\G^s$. The class of all candidate FTCGs for $\Gs$ is denoted by $\C(\Gs)$.


\hiddensubsection{Problem Setup}
We focus in this paper on identifying total effects \citep{Pearl_book2000} of multiple interventions on multiple effects, written $P \left(\bigl(Y^j_{t_j} = y^j_{t_j}\bigr)_j \middle\vert \bigl(\Do\bigl(X^i_{t_i} = x^i_{t_i}\bigr) \bigr)_i \right)$ (as well as $P\left((y^j_{t_j})_j \mid (\Do(x^i_{t_i}))_i\right)$ by a slight abuse of notation) \textit{when only the SCG reduced from the true FTCG is known}.  In general, an effect is said to be identifiable from a graph if it can be uniquely computed with a do-free formula from the observed distribution. In our context, this means that the same do-free formula should hold in all candidate FTCG so as to guarantee that it holds for the true one.
\begin{definition}[Identifiability of total effect in an SCG]
\label{identifiability}
Consider an SCG $\Gs$. The total effect $P\left((y^j_{t_j})_j \mid (\Do(x^i_{t_i}))_i\right)$ is \emph{identifiable} if and only if it can be rewritten with a do-free formula that is valid for any FTCG in $\C(\Gs)$.
\end{definition}
In this paper, we focus on a specific notion of identification in which there exists a set of variables $\mathcal{Z}$ that satisfies the backdoor criterion in every FTCG of $\C(\Gs)$. This corresponds to a multivariate backdoor adjustment which states that total effects are identifiable with variables satisfying the multivariate backdoor criterion in Definition~\ref{def:critere_backdoor_multivarie}, as formulated below:
\begin{restatable}{proposition}{mypropbackdooradjustmentmultivarie}[(Multivariate) Backdoor Adjustment]
\label{prop:backdoor_adjustment_multivarie}
    Consider an SCG $\Gs = (\mathcal{V}^s, \mathcal{E}^s)$ and the total effect $P\left((y^j_{t_j})_j \mid (\Do(x^i_{t_i}))_i\right)$. If a  set $\mathcal{Z} \subseteq \mathcal{V}^s$ satisfies the multivariate backdoor criterion relative to $\{X^i_{t_i}\}_i$ and $\{Y^j_{t_j}\}_j$ in the true FTCG, then the total effect is identifiable by:
    \begin{equation}
    \label{eq:backdoor_adjustment}
        P\left((y^j_{t_j})_j \mid (\Do(x^i_{t_i}))_i\right)
        = \sum_z P\left((y^j_{t_j})_j  \mid (x^i_{t_i})_i, z\right) P(z),
    \end{equation}
\noindent where $P$ is the true probability distribution. The sum over $z$ is taken over all possible valuations of $\mathcal{Z}$.
\end{restatable}
This leads to the following definition of \textit{identifiability by common backdoor}, which implies that if there exists a common backdoor set, then it is a valid backdoor set for the true FTCG and the true probability distribution, so that the total effect is identifiable with this set.
\begin{definition}[Identifiability by Common Backdoor]
\label{def:IBC}
    Let $\Gs = (\mathcal{V}^s, \mathcal{E}^s)$ be an SCG. An effect $P\left((y^j_{t_j})_j \mid (\Do(x^i_{t_i}))_i\right)$ is said to be \emph{identifiable by common backdoor in $\Gs$} if there exists a set $\mathcal{Z} \subseteq \mathcal{V}^s$ such that for all $\Gf \in \mathcal{C}(\Gs), 
     ~ \mathcal{Z}$ satisfies the backdoor criterion relative to $\{X^i_{t_i}\}_i$ and $\{Y^j_{t_j}\}_j$.
\end{definition}

As the common backdoor set for multiple effects is also a backdoor set for each effect, as stated in Lemma \ref{lemma:observations_multiples}, our problem can be reformulated {with a single effect. Note that the problem with the multiple effects is identifiable if and only if all problems with single effects are identifiable. In that case, the union of all common backdoors for single effects is a valid common backdoor for all effects. Our problem thus takes the form}:

\begin{problem}
\label{pb:mains_objective}
    Consider an SCG $\Gs$. 
    We aim to find out {operational}\footnote{{That is conditions one can rely on in practice. In particular the number of candidate FTCGs in $\C(\Gs)$ is usually too costly to enumerate (it may even be infinite); operational conditions should thus not rely on the enumeration of all FTCGs.}} necessary and sufficient conditions to identify the total effect $P \left(y_t \mid \Do( x^1_{t-\gamma_1}), \dots , \Do( x^n_{t-\gamma_n}) \right)$ by common backdoor when having access solely to the SCG $\Gs$.
\end{problem}
Note that if $Y$ is not a descendant of one of the intervening variables $X^i$ in $\Gs$ {or if $\gamma_i < 0$}, then $x^i_{t-\gamma_i}$ can be removed from the conditioning set through, \textit{e.g.}, the adjustment for direct causes \citep{Pearl_book2000}.  In the extreme case where $Y$ is not a descendant of any element of $\{X^i\}_i$, then $P \left(y_t \mid \Do( x^1_{t-\gamma_1}), \dots , \Do( x^n_{t-\gamma_n}) \right) = P(y_t)$. In the remainder, we thus assume that $Y$ is a descendant of each element in $\{X^i\}_i$ in $\Gs$ {and that $\gamma_i \geq 0$ for all $i$}. 

\hiddensection{Identifiability by Common Backdoor \label{section:IBC}}

 We provide in this section the main results of this paper, which is a graphical necessary and sufficient condition  for identifiability of the causal effect by common backdoor, and a solution to compute it in practice. The classical \textit{consistency throughout time}, assuming that causal relations are the same at different time instants, is not assumed here and its discussion is postponed to Section \ref{sec:Consistency_Time}. All the proofs are deferred to Section \ref{sec:proof:4} in the Supplementary Material.
 
\hiddensubsection{Necessary and sufficient condition based on the cone of descendants \label{subsection:equiv_based_CD}}

We first introduce the \emph{cone of descendants}, the set of vertices that are descendants of at least one interventional variable $X^i_{t-\gamma_i}$ in at least one candidate FTCG. The cone of descendants, and the related notion of non-conditionable set defined below, define a set of variables which cannot be elements of a common backdoor set as they violate the first condition of the multivariate backdoor adjustment given in  Definition~\ref{def:critere_backdoor_multivarie}. As such, they cannot be used as conditioning variables in the do-free formula rewriting the interventions (Equation~\ref{eq:backdoor_adjustment} in Proposition~\ref{prop:backdoor_adjustment_multivarie}).

\begin{algorithm}[t]
\setlength{\rightskip}{-.5cm} 
\caption{Computation of $(t_{\NC}(S))_{S \in \mathcal{V}^s}$}\label{algo:calcul_t_NC}
\SetKwInput{KwData}{Input}
\SetKwInput{KwResult}{Output}
\KwData{ $\Gs = (\mathcal{V}^s, \mathcal{E}^s)$ an SCG and the effect $P(y_t \mid \text{do}(x^1_{t-\gamma_1}), \dots, \text{do}(x^n_{t-\gamma_n}))$.}
\KwResult{$(t_{\NC}(S))_{S \in \mathcal{V}^s}$}
$L \gets [X^i_{t-\gamma_i}]_{i \in \{ 1, \dots,n\}, \gamma_i \text{ sorted in decreasing order}}$ \;
$(t_{\NC}(S)) \gets +\infty \quad \forall S \in \mathcal{V}^s$\;
$S.seen \gets False \quad \forall S \in \mathcal{V}^s$\;
\For{$X^i_{t - \gamma_i} \in L$}{
    \ForEach{unseen $D \in \Desc^>(X^i,\Gs) $}{
        $t_{\NC}(D) \gets \min \{t_1 \mid t_1 \geq t-\gamma_i \text{ and } D_{t_1} \in \NC\}$\;
        $D.seen \gets true$ \;
    }
}
\end{algorithm}

\begin{definition}[Cone of Descendants]
\label{def:CD}
    Let $\Gs$ be an SCG and $P(y_t \mid \text{do}(x^1_{t-\gamma_1}), \dots, \text{do}(x^n_{t-\gamma_n}))$ be the considered effect. We define the \emph{cone of descendants} as follows:
    $$\CD \coloneqq \bigcup_{\Gf \in \mathcal{C}(\Gs)} \bigcup_{i \in \{1,\ldots,n\}} \Desc\left( X^i_{t - \gamma_i}, \Gf_i \right),$$
where $\Gf_i \coloneqq \G^f_{ \overline{(X^{j})_{j \neq i}} \underline{(X^{j})_{j \neq i}}}$ is the mutilated graph obtained by removing all incoming and outgoing edges of $(X^{j})_{j \neq i}$.
\end{definition}
We then define the set of non-conditionable variables. 
\begin{definition}[Set of non-conditionable variables]
    Let $\Gs = (\mathcal{V}^s, \mathcal{E}^s)$ be an SCG and $P(y_t \mid \text{do}(x^1_{t-\gamma_1}), \dots, \text{do}(x^n_{t-\gamma_n}))$ be the considered effect. The \emph{set of non-conditionable variables}  is defined by 
    $$\NC \coloneqq \CD \backslash \{ X^i_{t-\gamma_i}\}_{i=1}^n.$$
    We also denote,  for a time series $F \in \mathcal{V}^s$, 
    $$t_{\NC}(F) \coloneqq \min \{ t_1 \mid F_{t_1} \in \NC \}$$
    the first time the time series $F$ is in the non-conditionable set $\NC$, with the convention $\min \{\emptyset\} = +\infty$. 
\end{definition}

In the next lemma, we show that $\{t_{\NC}(F)\}_{F \in \mathcal{V}^S}$ gives a simple characterization of these sets.

\begin{restatable}{lemma}{mylemmadefequivCD}{(Characterization of the cone of descendants)}
\label{lemma:def_equiv_CD}
Let $\Gs = (\mathcal{V}^s, \mathcal{E}^s)$ be an SCG and let $P(y_t \mid \text{do}(x^1_{t-\gamma_1}), \dots, \text{do}(x^n_{t-\gamma_n}))$ be the considered effect. With the convention $\{F_{t_1}\}_{t_1 \geq +\infty} = \emptyset$, we have:
    \begin{align*}
    \CD &=  \bigcup_{Z \in \mathcal{V}^S} \{Z_{t_1}\}_{t_1 \geq t_{\NC}(Z)} \cup \left\{ X^i_{t-\gamma_i}\right\}_i;\\
 \NC &=  \bigcup_{Z \in \mathcal{V}^S} \{Z_{t_1}\}_{t_1 \geq t_{\NC}(Z)} \backslash \left\{ X^i_{t-\gamma_i}\right\}_i.
\end{align*}
Moreover, $(t_{\NC}(F))_{F \in \mathcal{V}^s}$ {can be} computed through Algorithm \ref{algo:calcul_t_NC}, {detailed} in Appendix \ref{sec:proof:4}, {which} complexity is 
linear.
\end{restatable}

Theorem~\ref{th:equiv_IBC_multivarie} below shows that identifiability by common backdoor is directly related to the existence of collider-free backdoor path remaining in the cone of descendants. 

\begin{restatable}{theorem}{mytheoremequivIBC}{}
    \label{th:equiv_IBC_multivarie}
    Let $\Gs = (\mathcal{V}^s, \mathcal{E}^s)$ be an SCG and $P(y_t \mid \Do(x^1_{t-\gamma_1}), \dots, \Do(x^n_{t-\gamma_n}))$ be the considered effect. Then the two statements are equivalent:
    \begin{enumerate}
        \item The effect is identifiable by common backdoor in $\Gs$. \label{th:equiv_IBC_multivarie:1}
        \item For all  $X^i_{t-\gamma_i}$ and all FTCGs $\Gf$ belonging to $\C(\Gs)$, there does not exist a collider-free backdoor path going from $X^i_{t-\gamma_i}$ to $Y_t$ that remains in $\CD$ in $\Gf_i$. \label{th:equiv_IBC_multivarie:2}
    \end{enumerate}
    
    In that case, a backdoor set is given by $\C \coloneqq \mathcal{V}^f \backslash \CD$, and we have 
        \begin{equation*}
        P((y^j_{t_j})_j \mid (\Do(x^i_{t-\gamma_i}))_i))
        = \sum_{z \in \mathcal{C}} P\left((y^j_{t_j})_j  \mid (x^i_{t-\gamma_i})_i, z\right) P(z).
    \end{equation*}
\end{restatable}

\begin{figure}[t]
    \begin{tikzpicture}[scale = 1]
        \coordinate (Y) at (0,0);
        \coordinate (Xi) at (-1.5,3);
        \coordinate (Xj) at (1.5,1.25);
        \coordinate (D) at (.05,1.9);
        \coordinate (C) at (.45,1.9);
        
        \draw[dashed] (Xi) -- (-2,3);
        \draw[dashed] (-2,3) -- (-2,0);
        \draw[dashed] (Xi) -- (.25,3);
        \draw[dashed] (.25,3) -- (.25,0);
        
        \draw[dashed] (Xj) -- (-.5,1.25);
        \draw[dashed] (-.5,1.25) -- (-.5,0);
        \draw[dashed] (Xj) -- (2,1.25);
        \draw[dashed] (2,1.25) -- (2,0);

        \draw[->] (-2.5,3.25) -- (-2.5,-.25);
        \fill (-2.5,3)  circle (.5pt) node[left] {\small $t-\gamma_i$};
        \fill (-2.5,1.25)  circle (.5pt) node[left] {\small $t-\gamma_j$};
        \fill (-2.5,0)  circle (.5pt) node[left] {\small $t$};
        
        \draw[green!50!black, thick, decorate, decoration={snake, pre length = 2 mm, amplitude=.4mm, segment length=3mm, post=lineto, post length=1mm}, -] (Xi) .. controls (-1.2,2.2) and (-0.5,1.9) .. (D);
        \draw[green!50!black, thick, {latex[scale=.9]}-] (D) -- (C);
        \draw[green!50!black, thick, decorate, decoration={snake, amplitude=.4mm, segment length=3mm, post length = 1 mm}] (C) .. controls (1.5,1.5) and (-0.15,.8) .. (Y);

        \fill (Y)  circle (1pt) node[below] {$Y_t$};
        \fill (Xi) circle (1pt) node[above] {$X^i_{t-\gamma_i}$};
        \fill (Xj) circle (1pt) node[above right] {$X^j_{t-\gamma_j}$};
        \fill (D)  circle (1pt) node[above left] {\small $D_{t_d}$};
        \fill (C)  circle (1pt) node[above right] {\small $C_{t_c}$};
    \end{tikzpicture}
    \begin{tikzpicture}[overlay, shift={(-1,3.75)}, scale=0.3]
        \coordinate (X) at (1,1);
        \fill (X)  circle (3pt) node[above] {$X_{t_x}$};
        \draw[dashed] (X) -- (0,1);
        \draw[dashed] (0,1) -- (0,0);
        \draw[dashed] (X) -- (2,1);
        \draw[dashed] (2,1) -- (2,0);
        
        \node[anchor=west] at (2.1,0.5) {\smaller$ \bigcup_{\Gf \in \mathcal{C}(\Gs)} \Desc(X_{t_x}, \Gf)$};
    
        
        \draw[green!50!black, thick, decorate, decoration={snake, amplitude=.4mm, segment length=3mm}] (0.5,-.95) -- (1.5,-.95);
        \node[anchor= west] at (2.1,-.95) {\smaller Path \(\pi^f\) leaving \(\CD\)};

        \node at (1, -1.95) {\small $C_{t_c}$};
        \node[anchor= west] at (2.1,-1.95) {\smaller First vertex of \(\pi^f\) in \(\mathcal{V}^f \backslash \CD\)};

        \node at (1, -2.95) {\small $D_{t_d}$};
        \node[anchor= west] at (2.1,-2.95) {\smaller Predecessor of \(C\) in \(\pi^f\)};
        
    \end{tikzpicture}
    
    \caption{Illustration of the cone of descendants. The path in \textcolor{green!50!black}{green} demonstrates that all paths leaving the Cone of Descendants ($\CD$) exit it through an arrow pointing in the direction of the start of the path.
 }\label{fig:CD}
\end{figure}
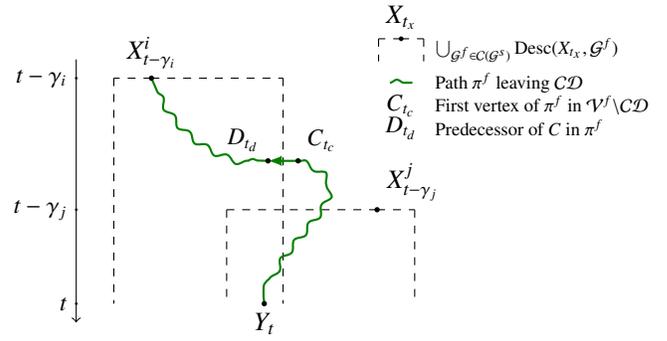

The intuition behind Theorem \ref{th:equiv_IBC_multivarie} is that the effect is identifiable by a common backdoor in $\Gs$ if and only if all collider-free backdoor paths from the intervention to $Y_t$ conform to the structure presented in Figure \ref{fig:CD}.

\hiddensubsection{In practice \label{ssct:in practice}}

To determine whether the causal effect is identifiable, we propose a polynomial-time algorithm that efficiently tests for the existence of collider-free backdoor paths which remains in $\CD$. However, instead of enumerating all FTCGs and collider-free backdoor paths in $\CD$, which would be computationally prohibitive, we introduce a more refined approach to characterize their existence. Specifically, we differentiate between paths with and without forks. {In particular, collider-free backdoor paths without forks can be easily and efficiently identified. The situation is a bit more complex for collider-free backdoor paths with forks but one can still efficiently identify them through a divide-and-conquer approach on nodes in $\NC$ connected to both an intervention and the effect through a directed path. All these elements are detailed in the following subsections.}

\subsubsection{Collider-free backdoor paths without fork in \texorpdfstring{$\CD$}{CD}}
\label{subsubsection:CfbWF}
First, Lemma \ref{lemma:IBC_enumeration_chemins_diriges}  characterizes efficiently the existence of collider-free backdoor paths that do not contain a fork.

\begin{restatable}{lemma}{mylemmaIBCenumerationcheminsdiriges}{(Characterization of collider-free backdoor paths without fork)}
\label{lemma:IBC_enumeration_chemins_diriges}
    Let $\Gs$ be an SCG and $P(y_t \mid \text{do}(x^1_{t-\gamma_1}), \dots, \text{do}(x^n_{t-\gamma_n}))$ be the considered effect. The following statements are equivalent:
    \begin{enumerate}
        \item There exists an intervention $X^i_{t-\gamma_i}$ and an FTCG $\Gf$ belonging to $\C(\Gs)$  such that the mutilated graph  $\Gf_i$ contains  $X^i_{t- \gamma_i} \leftsquigarrow Y_t$ which remains in $\CD$.
        \label{lemma:IBC_enumeration_chemins_diriges:1}
        \item There exists an intervention $X^i_{t-\gamma_i}$ such that $ \gamma_i = 0$ and $X^i \in \Desc \left(Y, \Gs_{\mid \mathcal{S}} \right)$, where $\mathcal{S} \vcentcolon = \bigcup_{i \in \{1,\ldots,n\}} \Desc \left( X^i, \Gs\right)$.
        \label{lemma:IBC_enumeration_chemins_diriges:2}
    \end{enumerate}
\end{restatable}


Condition \ref{lemma:IBC_enumeration_chemins_diriges:2} in Lemma \ref{lemma:IBC_enumeration_chemins_diriges} can be computed efficiently. Indeed, it is sufficient to traverse $\Gs_{\mid \mathcal{S}}$, and when an intervention is encountered, verify whether it occurs at time t. This can be achieved with the standard time complexity of graph traversal algorithms $\mathcal{O}(\left| \mathcal{V}^s \right| + \left| \mathcal{E}^s \right|)$ \citep[Chapter~22]{Cormen}. 

\hiddensubsubsection{Collider-free backdoor paths containing a fork in \texorpdfstring{$\CD$}{CD}\label{subsubsection:CfbF}}

We first introduce an accessibility concept essential to the enumeration of forked paths. 

\begin{definition}[$\mathcal{NC}$-accessibility]
\label{def:F-E-acc}
Let $\mathcal{V}$ be a  set of $p$ observational time series, 
 $\Gs = (\mathcal{V}^s, \mathcal{E}^s)$ be an SCG, $P(y_t \mid \text{do}(x^1_{t-\gamma_1}), \dots, \text{do}(x^n_{t-\gamma_n}))$ be the considered effect and $V_{t_v}\in \mathcal{V}^f$.  We say that \emph{$F_{t_1}\in \mathcal{V}^f \backslash \{V_{t_v}\}$ is $V_{t_v}$-$\NC$-accessible} if there exists an FTCG belonging to $\mathcal{C} (\Gs)$ in which there exists a directed path from $F_{t_1}$ to $V_{t_v}$ which remains\footnote{$V_{t_v}$ is not in $\NC$ when it is an intervention.} in $\NC$ except perhaps for $V_{t_v}$.
    We denote 
    $$t^{\NC}_{V_{t_v}}(F) \coloneqq \max \{ t_1 \mid F_{t_1} \text{is $V_{t_v}$-$\NC$-accessible} \},$$
    with the convention $\max\{\emptyset\} = -\infty$.
\end{definition}

This leads to characterise  efficiently the existence of a collider-free backdoor path with a fork that remains in $\CD$, as proposed in  Lemma \ref{lemma:equiv_existence_chemin_fork_CD_without_consistency_through_time}.

\begin{restatable}{lemma}{mylemmaequivexistencecheminforkCDwithoutconsistencythroughtime}{(Characterization of collider-free backdoor paths with fork)}
\label{lemma:equiv_existence_chemin_fork_CD_without_consistency_through_time}
    Let $\mathcal{V}$ be a set of $p$ observational time series. Let $\Gs$ be an SCG and $P(y_t \mid \text{do}(x^1_{t-\gamma_1}), \dots, \Do(x^n_{t-\gamma_n}))$ be the considered effect such that for all $\Gf$ belonging to $\C(\Gs)$, the mutilated graph $\Gf_i$ does not contain {a directed path from $Y_t$ to an intervention} $X^i_{t- \gamma_i} \leftsquigarrow Y_t$ which remains in $\CD$. The following statements are equivalent:
    \begin{enumerate}
        \item  There exists an intervention $X^i_{t-\gamma_i}$, $F_{t_f} \in \mathcal{V}^f$ and an FTCG $\Gf$ belonging to $\C(\Gs)$  such that the mutilated graph $\Gf_i$ contains the path $X^i_{t- \gamma_i} \leftsquigarrow F_{t_f} \rightsquigarrow Y_t$ which remains in $\CD$.
        \label{lemma:equiv_existence_chemin_fork_CD_without_consistency_through_time:1}

        \item There exists an intervention $X^i_{t- \gamma_i}$ and $F_{t_f} \in \mathcal{V}^f$ such that $F_{t_f}$ is $X^i_{t- \gamma_i}$-$\NC$-accessible and $Y_t$-$\NC$-accessible.
        \label{lemma:equiv_existence_chemin_fork_CD_without_consistency_through_time:2}
    \end{enumerate}
\end{restatable}

The intuition behind Lemma \ref{lemma:equiv_existence_chemin_fork_CD_without_consistency_through_time} is the divide-and-conquer strategy mentioned above. 
{It exploits the fact that in an FTCG, a backdoor path with a fork is the concatenation of two directed paths. As directed paths can be efficiently identified, one can also efficiently prove the existence of an FTCG with a collider-free backdoor path by constructing it with two directed paths (see the proof in the Supplementary Material).}
%
Lemma \ref{lemma:calcul_V_E_acc} furthermore shows that there is an efficient algorithm to compute $t^{\NC}_{V_{t_v}}(F)$ for all time series $F \in \mathcal{V}^s$.


\begin{restatable}{lemma}{mylemmacalculVEacc}{}
    \label{lemma:calcul_V_E_acc}
    Let $\mathcal{V}$ be a set of $p$ observational time series,  $\Gs = (\mathcal{V}^s, \mathcal{E}^s)$ be an SCG,  $P(y_t \mid \text{do}(x^1_{t-\gamma_1}), \dots, \Do(x^n_{t-\gamma_n}))$ be the considered effect and $V_{t_v}\in \mathcal{V}^f$. The set $\{t^{\NC}_{V_{t_v}}(F) \mid F \in \mathcal{V}^S \}$ {can be} computed by Algorithm \ref{algo:calcul_V_E_acc} {which} complexity is $\mathcal{O} \left(\left| \mathcal{E}^s \right| + \left| \mathcal{V}^s \right| \log \left| \mathcal{V}^s \right| \right)$\footnote{In amortized time. With a binary heap, the complexity of the algorithm is $\mathcal{O} ((\left| \mathcal{E}^s \right| + \left| \mathcal{V}^s \right|) \log \left| \mathcal{V}^s \right|)$.}.
\end{restatable}

\begin{algorithm}[t]
\setlength{\rightskip}{-.5cm} 
\caption{Computation of $(t^{\NC}_{V_{t_v}}(S))_{S \in \mathcal{V}^s}$}\label{algo:calcul_V_E_acc}
\SetKwInput{KwData}{Input}
\SetKwInput{KwResult}{Output}
\KwData{ $\Gs = (\mathcal{V}^s, \mathcal{E}^s)$ an SCG and the effect $P(y_t \mid \text{do}(x^1_{t-\gamma_1}), \dots, \text{do}(x^n_{t-\gamma_n}))$ and $V_{t_v} \in \mathcal{V}^f$}
\KwResult{$(t^{\NC}_{V_{t_v}}(S))_{S \in \mathcal{V}^s}$}
$Q \gets \text{PriorityQueue}(V_{t_v})$ \;
$t^{\NC}_{V_{t_v}}(S) \gets -\infty \quad \forall S \in \mathcal{V}^s$\;
$S.seen \gets False \quad \forall S \in \mathcal{V}^s$\;
\While{$Q \neq \emptyset$}{
    $S_{t_s} \gets Q.\text{pop\_element\_with\_max\_time\_index}()$\;
    \ForEach{unseen $P \in \text{Pa}(S, \Gs)$}{
        $t^{\NC}_{V_{t_v}}(P) \gets \max \{t_1 \mid t_1 \leq t_s \text{ and } P_{t_1} \in \NC \backslash \{V_{t_v}\}\}$\;
        \lIf{$t^{\NC}_{V_{t_v}}(P) \neq -\infty$}{
            $Q.insert(P_{t^{\NC}_{V_{t_v}}(P)})$
        }
        $P.seen \gets true$ \;
    }
}
\end{algorithm}

From the viewpoint of complexity, there is little interest in replacing Algorithm \ref{algo:calcul_V_E_acc} with a formula.\footnote{We refer to a formula as a condition involving a combination of descendant, ancestor, and cycle sets, as exemplified in Theorem~\ref{th:CondNecSuf_IBC_monovarie}.} Indeed, calculating $t^{\NC}_{V_{t_v}}(F)$ for all $F$ is done with at least a complexity of $\mathcal{O}\left( \left|  \mathcal{V}^s \right| + \left| \mathcal{E}^s \right| \right)$ because it is necessary to traverse $\Gs$. 

As shown by Lemma \ref{lemma:CharactVNCAcc}, for a time series $F\in \mathcal{V}^s$, knowing $t^{\NC}_{V_{t_v}}(F)$ and $t_{\NC}(F)$ is enough to characterise efficiently the set $\{ t_1 \mid F_{t_1} \text{is $V_{t_v}$-$\NC$-accessible} \}$.

\begin{restatable}{lemma}{mylemmaCharactVNCAcc}{}
\label{lemma:CharactVNCAcc}
    Let $\mathcal{V}$ be a set of $p$ observational time series,  $\Gs = (\mathcal{V}^s, \mathcal{E}^s)$ be an SCG, $P(y_t \mid \text{do}(x^1_{t-\gamma_1}), \dots, \Do(x^n_{t-\gamma_n}))$ be the considered effect, $V_{t_v}\in \mathcal{V}^f$  and $F\in \mathcal{V}^s$. Then the following statements are equivalent:
    \begin{enumerate}
        \item $F_{t_f}$ is $V_{t_v}$-$\NC$-accessible.
        \label{lemma:CharactVNCAcc:1}
        \item $t_{\NC}(F) \leq t_f \leq t^{\NC}_{V_{t_v}}(F)$ and $F_{t_f} \notin \{ X^i_{t-\gamma_i}\}_i$.
        \label{lemma:CharactVNCAcc:2}
    \end{enumerate}
\end{restatable}

{When Lemma \ref{lemma:CharactVNCAcc} is applied with $t_f = t_{\NC}(F)$, it  identifies whether the time series $F$ is $V_{t_v}$-$\NC$-accessible for at least one time index. Corollary \ref{cor:1} uses this idea to reduce the second condition of Lemma \ref{lemma:equiv_existence_chemin_fork_CD_without_consistency_through_time} to straightforward calculations.
}

\begin{restatable}{corollary}{mycorUn}{}
\label{cor:1}
    Let $\mathcal{V}$ be a set of $p$ observational time series,  $\Gs = (\mathcal{V}^s, \mathcal{E}^s)$ be an SCG,  $P(y_t \mid \text{do}(x^1_{t-\gamma_1}), \dots, \Do(x^n_{t-\gamma_n}))$ be the considered effect, $V_{t_v} \in \mathcal{V}^f$ and $F \in \mathcal{V}^s$. Let $X^i_{t - \gamma_i}$ be a fixed intervention. The following statements are equivalent:
    \begin{enumerate}
        \item There exists $t_f$ such that $F_{t_f}$ is $X^i_{t - \gamma_i}$-$\NC$-accessible and $Y_t$-$\NC$-accessible.
        \label{cor:1:1}
        
        \item $t_{\NC}(F) \leq t^{\NC}_{X^i_{t - \gamma_i}}(F)$ and $t_{\NC}(F) \leq t^{\NC}_{Y_t}(F)$.
        \label{cor:1:2}
    \end{enumerate}
\end{restatable}

The second statement of Corollary \ref{cor:1} can be assessed by executing Algorithm \ref{algo:calcul_V_E_acc} twice, knowing that Algorithm \ref{algo:calcul_t_NC} has already been run. Consequently, the overall complexity is $\mathcal{O}(\left| \mathcal{E}^s \right| + \left| \mathcal{V}^s \right| \log \left| \mathcal{V}^s \right|)$.

\hiddensubsubsection{An efficient algorithm \label{sssct:algo_IBC}}
{By combining the results of Sections \ref{subsubsection:CfbWF} and \ref{subsubsection:CfbF} as in Algorithm \ref{algo:calcul_IBC}, one can directly assess whether the effect is identifiable or not:}


\begin{restatable}{theorem}{mythforalgoIBC}{}
    \label{th:th_for_algo_2}
    Let $\Gs = (\mathcal{V}^s, \mathcal{E}^s)$ be an SCG and $P(y_t \mid \Do(x^1_{t-\gamma_1}), \dots, \Do(x^n_{t-\gamma_n}))$ be the considered effect. Then the two statements are equivalent:
    \begin{itemize}
        \item The effect is identifiable by common backdoor in $\Gs$.
        \item Algorithm \ref{algo:calcul_IBC} outputs True. In that case, a common backdoor set is given by $\C \coloneqq \mathcal{V}^f \backslash \CD$.
    \end{itemize}
    The algorithm has a polynomial complexity of $\mathcal{O}\left(n \cdot (\left| \mathcal{E}^s \right| + \left| \mathcal{V}^s \right| \log \left| \mathcal{V}^s \right|)\right)$. 
\end{restatable}

\begin{algorithm}[t]
\caption{Identifiability by common backdoor.}
\label{algo:calcul_IBC}

\SetKwInput{KwData}{Input}
\SetKwInput{KwResult}{Output}
\SetKwFunction{FMain}{AccessibleSet}
\KwData{ $\Gs = (\mathcal{V}^s, \mathcal{E}^s)$ an SCG and the effect $P(y_t \mid \text{do}(x^1_{t-\gamma_1}), \dots, \text{do}(x^n_{t-\gamma_n}))$.}
\KwResult{A boolean indicating whether the effect is identifiable by common backdoor or not.}

\tcp{Enumeration of directed paths.} 
$\mathcal{S} \gets \bigcup_{i} \Desc \left( X^i, \Gs\right)$ \;
\If{$\exists i \in \{1,\ldots, n\}$ s.t. $ X^i \in \Desc \left(Y, \Gs_{\mid \mathcal{S}} \right)$ and $\gamma_i = 0$}{
    \Return{False}
}

\tcp{Enumeration of fork paths.} 
$(t_{\NC}(S))_{S \in \mathcal{V}^s} \gets $ Algorithm \ref{algo:calcul_t_NC} \;
\ForEach{$V_{t_v} \in \{ Y_t, X^1_{t-\gamma_1}, \cdots,X^n_{t-\gamma_n} \}$}{
    $(t^{\NC}_{V_{t_v}}(S))_{S \in \mathcal{V}^s} \gets$ Algorithm \ref{algo:calcul_V_E_acc} \;
}

\ForEach{$F \in \mathcal{V}^s, X^i_{t- \gamma_i} \in (X^j_{t- \gamma_j})_j$}{
    \If{$t_{\NC}(F) \leq t^{\NC}_{X^i_{t - \gamma_i}}(F)$ and $t_{\NC}(F) \leq t^{\NC}_{Y_t}(F)$}{
        \Return{False}
    }
}

\Return{True}
\end{algorithm}

{Note that} an efficient algorithm that reaches a pseudo-linear complexity is discussed in the Appendix. As before, from the point of view of complexity, there is little interest in replacing the efficient implementation of Algorithm \ref{algo:calcul_IBC} with a formula. Indeed, we can not expect having a complexity better than $\mathcal{O}\left( n + \left| \mathcal{E}^s \right| + \left| \mathcal{V}^s \right| \right)$ because in the worst case, it is necessary to traverse $\Gs$ and consider all interventions. 

\hiddensection{With consistency throughout time}\label{sec:Consistency_Time}
In practice, it is usually impossible to work with general FTCGs in which causal relations may change from one time instant to another, and people have resorted to the \textit{consistency throughout time} assumption (also referred to as Causal Stationarity in \citet{Runge_2018}), to obtain a simpler class of FTCGs. 

\begin{assumption}[\textit{Consistency throughout time}]
	\label{ass:Consistency_Time}
An FTCG $\G^f$ is said to be \emph{consistent throughout time} if all the causal relationships remain constant in direction throughout time. 
\end{assumption}

Under this assumption, the number of candidate FTCGs for a fixed SCG $\mathcal{G}^s$ is smaller, meaning that conditions to be identifiable are weaker and thus that more effects should be identifiable. We detail in  Section  \ref{sec:ctt:id} necessary and sufficient conditions to be identifiable, and in Section \ref{sec:mono} rewrite those conditions in the particular case of a single intervention.
All the proofs are deferred to Section \ref{sec:proof:5} in the Supplementary Material.

\hiddensubsection{Identifiability}
\label{sec:ctt:id}

Theorem \ref{th:equiv_IBC_multivarie} remains valid under Assumption \ref{ass:Consistency_Time}. Lemma \ref{lemma:IBC_enumeration_chemins_diriges} also holds because Assumption \ref{ass:Consistency_Time} only affects paths that traverse different time indices. The enumeration of collider-free backdoor paths with a fork that remains in $\CD$ is however more complex, as detailed below. 

\begin{restatable}{lemma}{mylemmaequivexistencecheminforkNC}{}
\label{lemma:equiv_existence_chemin_fork_NC}
    Let $\Gs = (\mathcal{V}^s, \mathcal{E}^s)$ be an SCG and $P(y_t \mid \text{do}(x^1_{t-\gamma_1}), \dots, \Do(x^n_{t-\gamma_n}))$ be the considered effect such that for all $\Gf$ belonging to $\C(\Gs)$, the mutilated graph $\Gf_i$ does not contain {a directed path from $Y_t$ to an intervention} $X^i_{t- \gamma_i} \leftsquigarrow Y_t$ which remains in $\CD$. The following statements are equivalent:
    \begin{enumerate}
        \item There exist an intervention $X^i_{t- \gamma_i}$,  $F_{t'} \in \mathcal{V}^f$ and an FTCG $\Gf$ belonging to $\C(\Gs)$ such that the mutilated graph $\Gf_i$ contains the path $X^i_{t- \gamma_i} \leftsquigarrow F_{t'} \rightsquigarrow Y_t$ which remains in $\CD$.
        \label{lemma:equiv_existence_chemin_fork_NC:1}

        \item At least one of the following conditions is satisfied:
        \begin{enumerate}
            \item There exist {an intervention} $X^i_{t- \gamma_i}$ and  $F \in \mathcal{V}^s$ such that 
            $F_{t_{\NC}(F)}$ is well defined, $X^i_{t- \gamma_i}$-$\NC$-accessible and $Y_t$-$\NC$-accessible, and
            $\left\{
            \begin{aligned}
                &F \neq Y \text{, or} \\
                &t - \gamma_i \neq t_{\NC(F)}.\\
            \end{aligned}
            \right.$
            \label{lemma:equiv_existence_chemin_fork_NC:2a}
            
            \item There exists {an intervention} $X^i_{t- \gamma_i}$ such that $t - \gamma_i = t_{\NC(Y)}$ and at least one of the following properties is satisfied:
                    \begin{enumerate}
                        \item $Y_{t_{\NC}(Y)}$ is $X^i_{t- \gamma_i}$ - $\NC$-accessible without using $X^i_{t- \gamma_i} \leftarrow Y_{t- \gamma_i}$ and $Y_t$ - $\NC$-accessible.
                        \label{lemma:equiv_existence_chemin_fork_NC:2b:i}
                        
                        \item $Y_{t_{\NC}(Y)}$ is $X^i_{t- \gamma_i}$ - $\NC$-accessible and $Y_t$ - $\NC$-accessible without using $X^i_t \rightarrow Y_t$.
                        \label{lemma:equiv_existence_chemin_fork_NC:2b:ii}
                    \end{enumerate}
            \label{lemma:equiv_existence_chemin_fork_NC:2b}
        \end{enumerate}
        \label{lemma:equiv_existence_chemin_fork_NC:2}
    \end{enumerate}
\end{restatable}

Lemma \ref{lemma:equiv_existence_chemin_fork_NC} characterizes the existence of a collider-free backdoor path containing a fork. While the conditions outlined are more complex than those in Corollary~\ref{cor:1}, they play the same role, and still require only a small number of calls to $\NC$-accessibility. Consequently, one can replace the conditions in the final loop of Algorithm~\ref{algo:calcul_IBC} with conditions 2.(a) and 2.(b) of Lemma~\ref{lemma:equiv_existence_chemin_fork_NC} to derive an algorithm for identifiability by common backdoor in $\Gs$ under \textit{consistency throughout time}, as stated in the following theorem which is the counterpart of Theorem~\ref{th:th_for_algo_2}.

\begin{restatable}{theorem}{mythforalgoIBCwith_consistency}{}
    \label{th:th_for_algo_3}
    Let $\Gs$ be an SCG that satisfies Assumption \ref{ass:Consistency_Time} and $P(y_t \mid \Do(x^1_{t-\gamma_1}), \dots, \Do(x^n_{t-\gamma_n}))$ be the considered effect. Then the two statements are equivalent:
    \begin{itemize}
        \item The effect is identifiable by common backdoor in $\Gs$.

        \item An adaptation of Algorithm \ref{algo:calcul_IBC} outputs True.
    \end{itemize}

    In that case, a backdoor set is given by $\C \coloneqq \mathcal{V}^f \backslash \CD$.

    In its simpler form, the adaptation of Algorithm~\ref{algo:calcul_IBC} still has a polynomial complexity of $\mathcal{O}\left(n \cdot (\left| \mathcal{E}^s \right| + \left| \mathcal{V}^s \right| \log \left| \mathcal{V}^s \right|)\right)$. A more efficient algorithm that reaches a pseudo-linear complexity is discussed in the Supplementary Material, Section \ref{sct:IBC_pseudo_lineaire}.

\end{restatable}

\hiddensubsection{Monovariate case}
\label{sec:mono}

\cite{assaad_identifiability_2024} {recently} established a set of sufficient conditions to characterize the identifibility by common backdoor in the monovariate case. We have reformulated these conditions and proved that they are necessary and sufficient in the following theorem. 

\begin{restatable}{theorem}{mythIBCmonovarie}{}
\label{th:CondNecSuf_IBC_monovarie}

    Let $\Gs$ be an SCG such that $X\in \Anc(Y,\Gs)$ and we consider the total effect $P(y_t\mid \Do(x_{t-\gamma}))$. The effect is identifiable by common backdoor in $\Gs$ if and only if:
\begin{itemize}
     \item For $\gamma = 0$, there is no collider-free backdoor path from $X$ to $Y$ in $\Gs_{\mid \Desc(X, \Gs)}$.

    \item For $\gamma = 1$, 
    \begin{align*}
                \text{either } & \text{Pa}(X, \Gs) \cap \Desc\left( X, \Gs \right) \subseteq \{X\},  \\
                \text{or } &
                    \Gs_{\mid \Desc\left( X, \Gs \right) \cap \left( \text{Par} \left( X, \Gs \right) \cup \text{Par} \left( Y, \Gs \right) \right)} \in \left\{ X \leftrightarrows Y, \leftselfloop X \leftrightarrows Y \right\};
    \end{align*}
    \item For $\gamma \geq 2$,
   $\text{Cycles}^>\left(X, \Gs  \right) =  \emptyset$.
\end{itemize}
    \noindent In these cases, a common backdoor set is given by $\mathcal{A}_{\gamma}$, where:
   \begin{align*}
        \mathcal{A}_0 &\coloneqq \bigcup_{\pi^s \text{ backdoor}} \Bigl\{ Z_t  \mid Z \in  \pi^s \backslash \Desc\left( X,\Gs \right) \Bigr\} \cup \Bigl\{ \left( Z_{t'} \right)_{ t' < t} \mid Z \in \Gs \Bigr\}\\
        \mathcal{A}_1 &\coloneqq  \mathcal{V}^f \backslash \CD;\\
        \mathcal{A}_{\gamma} &\coloneqq  \left\{ (Z_{t'})_{t' \leq t-\gamma} \mid Z \in \Anc\left(X, \Gs \right) \right\} \backslash \left\{ X_{t-\gamma} \right\}, \text{for }\gamma \geq 2.
   \end{align*}

\end{restatable}

Note that a similar theorem can be derived in the monovariate case without \textit{consistency throughout time}; we detail this in the Supplementary Material, Section \ref{sec:app:mono}. 

\hiddensection{Conclusion}
\label{sec:conclusion}

This work has established complete conditions for identifying the effect of multiple interventions on multiple effects by common backdoor in summary causal graphs of time series. Specifically, Theorem \ref{th:equiv_IBC_multivarie} shows that, both with and without \textit{consistency throughout time}, the problem reduces to testing the existence of collider-free backdoor paths from interventions to effects that remain within a specified set, the cone of descendants. This complete characterization allowed us to derive {efficient} algorithms to determine whether an effect is identifiable by common backdoor, again both with and without \textit{consistency throughout time}. All the provided proofs are constructive, meaning that whenever an effect is not identifiable by common backdoor, it is possible to explicitly exhibit a collider-free backdoor path that remains within $\CD$. Furthermore, assuming \textit{consistency throughout time} and considering only a single intervention and a single effect, we have shown that the sufficient conditions for identifiability with common backdoor recently established in~\cite{assaad_identifiability_2024} are also necessary.

Coming back to  the toy example in Figure \ref{fig:real_thermoregulation}, let us denote by {$B_t$}, $O_t$, $K_t$ and $L_t$ respectively the temperature in {the bathroom}, in the office, in the kitchen and in the living room at time $t$. Then, the effect $P(O_t = o_t \mid \Do(K_{t-1} = k_{t-1}), \Do(L_{t-1} = \ell_{t-1}))$ is not identifiable by common backdoor because the path $L_{t-1} \leftarrow B_{t-1} \rightarrow L_t \rightarrow O_t$ remains within the set $\CD = \{K_{t-1}, K_t, L_{t-1}, L_t, B_{t-1}, B_t, O_{t-1}, O_{t} \}$. On the other hand, the effect $P(O_t= o_t \mid \Do(K_{t-1} = k_{t-1}), \Do(L_{t-1} = \ell_{t-1}), \Do(L_{t} = \ell_t))$ is identifiable by common backdoor.

\newpage
\bibliography{References}

\appendix

\onecolumn
\aistatstitle{Supplementary Material}

\tableofcontents

\newpage

\section{How to build an FTCG belonging to \texorpdfstring{$\C(\Gs)$}{C(Gs)} ?}

Many proofs in this paper use constructive arguments to demonstrate the existence of an FTCG belonging to $\C(\Gs)$ that contains a given path $\pi^f$. To facilitate the understanding of these arguments, Lemma \ref{lemma:helper:construct_FTCG_with_path} shows how these FTCGs are constructed by adding missing arrows.

\subsection{Without Assumption \ref{ass:Consistency_Time}}

\begin{lemma}
    \label{lemma:helper:construct_FTCG_with_path}
    Let $\Gs$ be an SCG and $\pi^f$ be a path over $\mathcal{V}^f$. If $\pi^f$ is a DAG, all its arrows respect time orientation and its reduction is a subgraph of $\Gs$, then there exists an FTCG $\Gf$ in $\C(\Gs)$ which contains $\pi^f$. 
\end{lemma}

\begin{proof}
    $\Gs$ is a SCG, by definition it is the reduction of an FTCG $\Gf_{\star} =(\mathcal{V}^f_{\star},\mathcal{E}^f_{\star})$. Let $t_{\min}$ be the minimum time index seen by $\pi^f$. Let us consider the graph $\Gf = (\mathcal{V}^f \coloneqq \mathcal{V}^f_{\star}, \mathcal{E}^f)$ whose edges are constructed as follows:
    \begin{enumerate}
        \item All edges from $\pi^f$ are set in $\mathcal{E}^f$.
        \label{lemma:helper:construct_FTCG_with_path:step1}
        
        \item For all the edges $A \rightarrow B$ in $\Gs$ that are not reductions of arrows in $\pi^f$, add the edge $A_{t_{\min} -1} \rightarrow B_{t}$ into $\mathcal{E}^f$. 
        \label{lemma:helper:construct_FTCG_with_path:step2}

    \end{enumerate}
    
    We have the following properties:
    \begin{itemize}

        \item $\Gf$ contains $\pi^f$ because all arrows of $\pi^f$ are set in $\mathcal{E}^f$.
                
        \item $\Gf$ is a DAG. Indeed, there is no cycle due to instantaneous arrows because all instantaneous arrows come  from $\pi^f$ and $\pi^f$ is a DAG. Moreover, there is no cycle due to delayed arrows because all of them follow the flow of time.
                
        \item $\Gf$ belongs to $\C(\Gs)$: Let us consider $\G^r = (\mathcal{V}^r, \mathcal{E}^r)$ the reduction of $\Gf$. By definition, $\mathcal{V}^r = \mathcal{V} = \mathcal{V}^s$ . Les us prove that $\mathcal{E}^r = \mathcal{E}^s$ by showing the two inclusions:
        \begin{itemize}
            \item $\mathcal{E}^r \subseteq \mathcal{E}^s$: Let us consider an arrow $a^r$ of $G^r$ and $a^f$ a corresponding arrow in $\Gf$. We distinguish two cases:
            \begin{itemize}
                \item If $a^f$ is in $\pi^f$ then $a^r$ is in $\pi^r$ and $\pi^r$ is a subgraph of $\Gs$. Therefore  $a^r \in \mathcal{E}^s$.

                \item Otherwise, $a^r$ have been added to $\mathcal{E}^f$ during step \ref{lemma:helper:construct_FTCG_with_path:step2}. Therefore $a^r \in \mathcal{E}^s$.
            \end{itemize}
            In both cases $a^r \in \mathcal{E}^s$, therefore $\mathcal{E}^r \subseteq \mathcal{E}^s$.
            
            \item $\mathcal{E}^s \subseteq \mathcal{E}^r$:  Let us consider an arrow $A \rightarrow B$ of $G^s$. we distinguish two cases:
            \begin{itemize}
                \item If $A \rightarrow B$ is the reduction of an arrow $a^f$ of $\pi^f$, then $a^f \in \mathcal{E}^f$. Therefore $A \rightarrow B$ is an arrow of $\G^r$

                \item Otherwise, $A_{t_{\min} -1} \rightarrow B_{t}$ is in $\mathcal{E}^f$. Therefore $A \rightarrow B$ is an arrow of $\G^r$
            \end{itemize} 
            In both cases $A \rightarrow B$ is an arrow of $\G^r$, therefore $\mathcal{E}^s \subseteq \mathcal{E}^r$.
        \end{itemize}  
        Therefore, $\G^r = \G^s$ and thus $\Gf$ belongs to $\C(\Gs)$.
    \end{itemize}
    Therefore, $\Gf$ is an FTCG belonging to $\C(\Gs)$ which contains $\pi^f$.
\end{proof}

\subsection{With Assumption \ref{ass:Consistency_Time}}

The construction given in Lemma \ref{lemma:helper:construct_FTCG_with_path} does not work to build an FTCG that satisfies Assumption \ref{ass:Consistency_Time}. Indeed, some arrows are missing. By copying the construction of lemma \ref{lemma:helper:construct_FTCG_with_path} at each time step, if $\pi^f$ verify the correct properties, therefore it is possible to construct an FTCG belonging to $\C(\Gs)$ which contains $\pi^f$ and satisfy assumption \ref{ass:Consistency_Time}. This reasoning is encapsulated by Lemmma \ref{lemma:helper:construct_FTCG_with_path_with_cst}.

\begin{lemma}
    \label{lemma:helper:construct_FTCG_with_path_with_cst}
    Let $\Gs$ be an SCG and $\pi^f$ be a path over $\mathcal{V}^f$. If $\pi^f$ is a DAG, all its arrows respect time orientation, its reduction is a subgraph of $\Gs$ and all its instantaneous arrows respect Assumption \ref{ass:Consistency_Time}, then there exists an FTCG $\Gf$ in $\C(\Gs)$ which contains $\pi^f$ and satisfies Assumption \ref{ass:Consistency_Time}. 
\end{lemma}

\begin{proof}
    $\Gs$ is a SCG, by definition it is the reduction of an FTCG $\Gf_{\star} =(\mathcal{V}^f_{\star},\mathcal{E}^f_{\star})$. Let $t_{\min}$ be the minimum time index seen by $\pi^f$. Let us consider the graph $\Gf = (\mathcal{V}^f \coloneqq \mathcal{V}^f_{\star}, \mathcal{E}^f)$ whose edges are constructed as follows:
    \begin{enumerate}
        \item All edges from $\pi^f$ are set in $\mathcal{E}^f$.
        \label{lemma:helper:construct_FTCG_with_path_with_cst:step1}
        
        \item For all the edges $A \rightarrow B$ in $\Gs$ that are not reductions of arrows in $\pi^f$, add the edge $A_{t_{\min} -1} \rightarrow B_{t}$ into $\mathcal{E}^f$. 
        \label{lemma:helper:construct_FTCG_with_path_with_cst:step2}

        \item Copy these arrows at each time step.
        \label{lemma:helper:construct_FTCG_with_path_with_cst:step3}
    \end{enumerate}

    We have the following properties:
    \begin{itemize}

        \item $\Gf$ contains $\pi^f$ because all arrows of $\pi^f$ are set in $\mathcal{E}^f$.
                
        \item $\Gf$ is a DAG: All instantaneous arrows come from $\pi^f$ and its copies. Since $\pi^f$ is a DAG and all its instantaneous arrows respect Assumption \ref{ass:Consistency_Time}, we know that there is no cycle made of instantaneous arrows. Moreover, delayed arrows follow the flow of time, therefore, they cannot form a cycle. Therefore, $\Gf$ does not contain any cycle.
        
        Moreover, there is no cycle due to delayed arrows because all of them follow the flow of time.
                
        \item $\Gf$ belongs to $\C(\Gs)$: By the same reasoning of Lemma \ref{lemma:helper:construct_FTCG_with_path}, at Step \ref{lemma:helper:construct_FTCG_with_path_with_cst:step2}, the reduction $\G^r_{\text{at Step \ref{lemma:helper:construct_FTCG_with_path_with_cst:step2}}}$ of $\G^f_{\text{at Step \ref{lemma:helper:construct_FTCG_with_path_with_cst:step2}}}$ is equal to $\G^s$. Adding copies of arrows does not change the reduction. Therefore, the reduction of $\Gf$ is equal to $\Gs$ i.e. $\Gf$ belongs to $\C(\Gs)$.

        \item $\Gf$ satisfies Assumption \ref{ass:Consistency_Time} because copying all arrows ensure that all causal relationships remain constant throughout time
    \end{itemize}
    Therefore, $\Gf$ is an FTCG belonging to $\C(\Gs)$ which contains $\pi^f$ and satisfies Assumption \ref{ass:Consistency_Time}.
\end{proof}

\newpage

\section{Proofs of Section \ref{section:IBC}}\label{sec:proof:4}

\subsection{Proofs of Section \ref{subsection:equiv_based_CD}}

\mylemmadefequivCD*

\begin{proof}
    \textbf{We start by proving the characterizations of $\CD$ and $\NC$.} \newline
    \noindent By the transitivity of descendant relationships, we have: \newline
    $$
    \CD = \bigcup_{\Gf \in \mathcal{C}(\Gs)} \bigcup_{i \in \{1,\ldots,n\}} \Desc\left( X^i_{t - \gamma_i}, \Gf \right)
    $$

    Le us denote  $\mathcal{S} \coloneqq  \bigcup_{Z \in \mathcal{V}^S} \{Z_{t_1}\}_{t_1 \geq t_{\NC}(Z)} \cup \left\{ X^i_{t-\gamma_i}\right\}_i$. We prove $\CD = \mathcal{S}$ by double inclusion:
    \begin{itemize}
        \item $\CD \subseteq \mathcal{S}$: Let $Z_{t_1} \in \CD$. We distinguish two cases:
        \begin{itemize}
            \item If $Z_{t_1} \in \left\{ X^i_{t-\gamma_i}\right\}_i$, then $Z_{t_1} \in \mathcal{S}$ by definition of $\mathcal{S}$.
            \item Otherwise, $Z_{t_1} \in \NC$, hence $t_1 \geq t_{\NC}(Z)$ thus $Z_{t_1} \in \mathcal{S}$. 
        \end{itemize}
        In both cases, $Z_{t_1} \in \mathcal{S}$. Therefore $\CD \subseteq \mathcal{S}$.
        
        \item $\mathcal{S} \subseteq \CD$: Let $Z_{t_1} \in \mathcal{S}$. We distinguish two cases:
        \begin{itemize}
            \item If $Z_{t_1} \in \left\{ X^i_{t-\gamma_i}\right\}_i$, then $Z_{t_1} \in \CD$ because in any FTCG $\Gf \in \C(\Gs)$, $ X^i_{t-\gamma_i} \in \Desc(X^i_{t-\gamma_i}, \Gf)$.
            \item Otherwise, $t_1 \geq t_{\NC}(Z)$, hence $t_{\NC}(Z) < + \infty$ thus there exists $i \in \{1,\cdots,n\}$ and an FTCG $\Gf \in \C(\Gs)$ in which there exists $\pi^f$ a directed path from $X^i_{t-\gamma_i}$ to $Z_{t_{\NC}(Z)}$. We can construct another FTCG belonging to $\C(\Gs)$ containing $\pi^f$ except that the last arrow points to $Z_{t_1}$ instead of $Z_{t_{\NC}(Z)}$. Thus, $Z_{t_1} \in \CD$.
        \end{itemize}
        In both cases, $Z_{t_1} \in \CD$. Therefore $\mathcal{S} \subseteq \CD $.
    \end{itemize}

    \noindent Therefore, $\CD =  \bigcup_{Z \in \mathcal{V}^S} \{Z_{t_1}\}_{t_1 \geq t_{\NC}(Z)} \cup \left\{ X^i_{t-\gamma_i}\right\}_i$ and $\NC = \CD \backslash \left\{ X^i_{t-\gamma_i}\right\}_i  = \bigcup_{Z \in \mathcal{V}^S} \{Z_{t_1}\}_{t_1 \geq t_{\NC}(Z)} \backslash \left\{ X^i_{t-\gamma_i}\right\}_i$

    \textbf{We now prove that Algorithm \ref{algo:calcul_t_NC} is correct:} \newline

    Let $F$ be a time series. We distinguish two cases:
    \begin{itemize}
        \item If $\{ t_1 \mid F_{t_1} \in \NC\} = \emptyset$: By definition, $t_{\NC}(F) = + \infty$. In that case, for every $i$, $F \notin \Desc^>(X^i,\Gs)$. Therefore the algorithm never updates its value for $t_{\NC}(F)$ and outputs the correct result.
        \item Otherwise $\{ t_1 \mid F_{t_1} \in \NC\} \neq \emptyset$ and $F_{t_{\NC}(F)} \in \NC$. In that case, let us consider $X^{i_0}_{t-\gamma_{i_0}}$, the first intervention encountered by the algorithm for which there exists an FTCG $\Gf$ which contains the path $X^{i_0} \rightsquigarrow Y_t$. Thus, $F \in \Desc^>(X^{i_0},\Gs)$ and by definition of $i_0$, $F$ is unseen. Therefore, the algorithm computes $t_{\NC}(F) \gets \min \{t_1 \mid t_1 \geq t-\gamma_{i_0} \text{ and } F_{t_1} \in \NC\}$. It is the correct answer. Indeed, by definition,  $t_{\NC}(F) \leq \min \{t_1 \mid t_1 \geq t-\gamma_{i_0} \text{ and } F_{t_1} \in \NC\}$. By contradiction, we show that $t_{\NC}(F) \geq \min \{t_1 \mid t_1 \geq t-\gamma_{i_0} \text{ and } F_{t_1} \in \NC\}$. If $t_{\NC}(F) < \min \{t_1 \mid t_1 \geq t-\gamma_{i_0} \text{ and } F_{t_1} \in \NC\}$, then $t_{\NC}(F) < t- \gamma_{i_0}$. Hence, there exists $i_1$ such that $\exists~ \Gf \mid X^{i_1}_{t-\gamma_{i_1}} \rightsquigarrow F_{t_{\NC}(F)}$. Causality does not move backwards in time, thus $t-\gamma_{i_1} \leq t_{\NC}(F) < t-\gamma_{i_0}$. Since interventions are sorted by increasing time indices, $X^{i_1}_{t-\gamma_{i_1}}$ would have already been encountered by the algorithm, which contradicts the definition of $X^{i_0}_{t-\gamma_{i_0}}$.
    \end{itemize}

    In both cases, the algorithm computes the correct time, therefore, it is correct.

    \textbf{We now prove that Algorithm \ref{algo:calcul_t_NC} works in linear complexity:} \newline

The line $t_{\NC}(D) \gets \min \{t_1 \mid t_1 \geq t-\gamma_i \text{ and } D_{t_1} \in \NC\}$ of Algorithm \ref{algo:calcul_t_NC} hides computations that could be expensive. In fact, it can be computed efficiently in $\mathcal{O}(1)$ during the run of Algorithm \ref{algo:calcul_t_NC}. To do so, it is necessary to pre-run Algorithm \ref{algo:calcul_t_NC:1} in order to compute the set $\mathcal{H} \coloneqq \{ \min \{t_1 \mid t_1 \geq t-\gamma_i \text{ and } X^i_{t_1} \notin \{ X^i_{t-\gamma_i}\}_i\}\}_{i \in \{1,\cdots,n\}}$ in $\mathcal{O}(n)$. For each intervention $X^i_{t-\gamma_i}$, $\mathcal{H}$ gives the smallest time $t_1$ after $t-\gamma_i$ such that $X^i_{t_1}$ is not an intervention. When $\mathcal{H}$ is known, it is possible to compute $\min \{t_1 \mid t_1 \geq t-\gamma_i \text{ and } D_{t_1} \in \NC\}$ in $\mathcal{O}(1)$. To show this, we distinguish two cases for unseen vertex $D$ in the foreach loop of Algorithm ~\ref{algo:calcul_t_NC}:
\begin{itemize}
    \item $D_{t - \gamma_i}$ is not an intervention: since $D_{t - \gamma_i}$ is a descendant of an intervention in some FTCG,  $D_{t - \gamma_i} \in \NC$. Thus,  $\min \{t_1 \mid t_1 \geq t-\gamma_i \text{ and } D_{t_1} \in \NC\} = t-\gamma_i$. Therefore, Algorithm \ref{algo:calcul_t_NC} need to compute $t_{\NC}(D) \gets t - \gamma_i $ and that is done in $\mathcal{O}(1)$.
    
    \item $D_{t - \gamma_i}$ is an intervention: let $t_D \coloneqq \min \{t_1 \mid t_1 \geq t-\gamma_i \text{ and } D_{t_1} \notin \{ X^i_{t-\gamma_i}\}_i \} \in \mathcal{H}$. Since $D$ is a strict descendant of $X^i$ in $\Gs$, there exists an FTCG in $\C(\Gs)$ with a directed path $\pi^f$ from $X_{t - \gamma_i}$ to $D_{t - \gamma_i}$. Then by changing the last arrow of $\pi^f$, we can construct another FTCG belonging to $\C(\Gs)$ which contains a path from an intervention to $D_{t_D}$. Hence, $D_{t_D}$ is a descendant of an intervention in some FTCG and it is not an intervention; therefore, $D_{t_D} \in NC$. Since for all $t_1 \in \{t-\gamma_i, \dots, t_D -1 \}$, $D_{t_1}$ is an intervention by definition of $t_D$, we know that for all $t_1 \in \{t-\gamma_i, \dots, t_D -1 \}$, $D_{t_1} \notin \NC$. Therefore $t_D = \min \{t_1 \mid t_1 \geq t-\gamma_i \text{ and } D_{t_1} \in \NC\}$. Algorithm \ref{algo:calcul_t_NC} thus needs to compute $t_{\NC}(D) \gets t_D $ and it can be done in $\mathcal{O}(1)$ by picking the value in $\mathcal{H}$.
\end{itemize}

\begin{algorithm}[H]
    \caption{Computation of $\{t_1 \mid t_1 \geq t-\gamma_i \text{ and } X^i_{t_1} \notin \{ X^i_{t-\gamma_i}\}_i\}_{i \in \{1,\cdots,n\}}$}
    \label{algo:calcul_t_NC:1}
    \SetKwInput{KwData}{Input}
    \SetKwInput{KwResult}{Output}
    
    \KwData{List of interventions, sorted by decreasing time indices.}
    \KwResult{$\{ \min \{t_1 \mid t_1 \geq t-\gamma_i \text{ and } X^i_{t_1} \notin \{ X^i_{t-\gamma_i}\}_i\}\}_{i \in \{1,\cdots,n\}}$, denoted as $\{t_{X^i_{t-\gamma_i}}\}_i$}
    
    $L \gets$ List of lists of interventions, grouped by time series, preserving the time index ordering\;
    
    \For{$l \in L$}{
        $X^i_{t-\gamma_i} \gets l[0]$\;
        $t_{X^i_{t-\gamma_i}} \gets t - \gamma_i + 1$\;
        
        \ForEach{$X^j_{t-\gamma_j}$ in $l[1:]$}{
            Let $X^i_{t-\gamma_i}$ be the predecessor of $X^j_{t-\gamma_j}$ in $l$\;
            \eIf{$t-\gamma_i +1 = t-\gamma_j$}{
                $t_{X^j_{t-\gamma_j}} \gets t_{X^i_{t-\gamma_i}}$\;
            }{
                $t_{X^j_{t-\gamma_j}} \gets t-\gamma_j + 1$\;
            }
        }
    }
\end{algorithm}

During the run of Algorithm \ref{algo:calcul_t_NC}, all computations in the loops are done in $\mathcal{O}(1)$. By sharing the unseen set among all calculations of $\Desc^>(X^i,\Gs)$, all time series that are not interventions are seen at most one time and all interventions are seen at most two times\footnote{At most one time during the computation of $\Desc^>(X^i,\Gs)$ of another intervention and at most one time during the computation of $\Desc^>(X^i,\Gs)$ for itself.} Similarly, arrows from a non interventional time series are seen at most one time and arrows from an intervention time series are seen at most two time. By taking into account the complexity of Algorithm \ref{algo:calcul_t_NC:1}, the complexity of Algorithm \ref{algo:calcul_t_NC} is $\mathcal{O}(n+\left| \mathcal{E}^s \right|+\left| \mathcal{V}^s \right|)$.

\end{proof}

\mytheoremequivIBC*

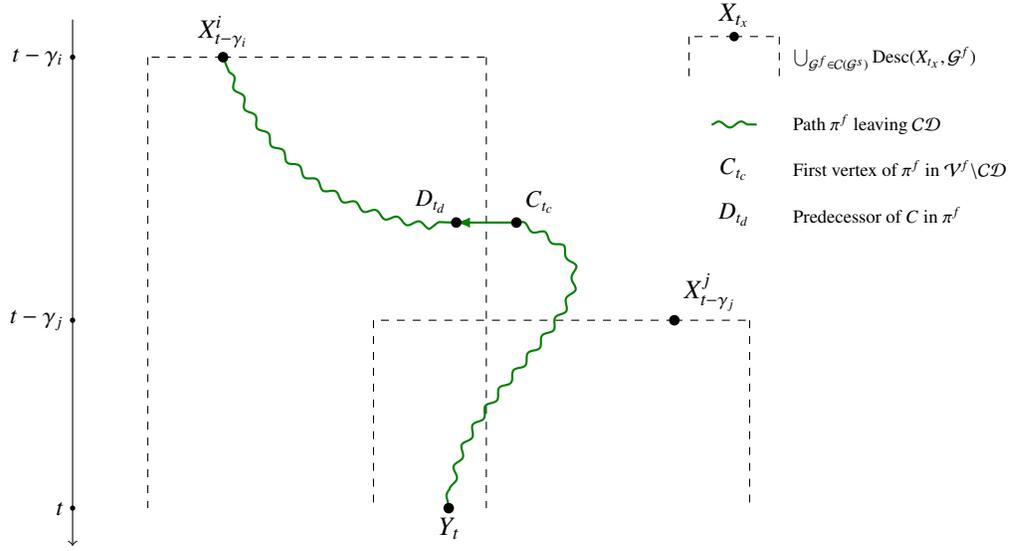
\begin{figure}[t]
    \centering
    \begin{tikzpicture}[scale = 2]
        \coordinate (Y) at (0,0);
        \coordinate (Xi) at (-1.5,3);
        \coordinate (Xj) at (1.5,1.25);
        \coordinate (D) at (.05,1.9);
        \coordinate (C) at (.45,1.9);
        
        \draw[dashed] (Xi) -- (-2,3);
        \draw[dashed] (-2,3) -- (-2,0);
        \draw[dashed] (Xi) -- (.25,3);
        \draw[dashed] (.25,3) -- (.25,0);
        
        \draw[dashed] (Xj) -- (-.5,1.25);
        \draw[dashed] (-.5,1.25) -- (-.5,0);
        \draw[dashed] (Xj) -- (2,1.25);
        \draw[dashed] (2,1.25) -- (2,0);

        \draw[->] (-2.5,3.25) -- (-2.5,-.25);
        \fill (-2.5,3)  circle (.5pt) node[left] {\small $t-\gamma_i$};
        \fill (-2.5,1.25)  circle (.5pt) node[left] {\small $t-\gamma_j$};
        \fill (-2.5,0)  circle (.5pt) node[left] {\small $t$};
        
        \draw[green!50!black, thick, decorate, decoration={snake, pre length = 2 mm, amplitude=.4mm, segment length=3mm, post=lineto, post length=1mm}, -] (Xi) .. controls (-1.2,2.2) and (-0.5,1.9) .. (D);
        \draw[green!50!black, thick, {latex[scale=.9]}-] (D) -- (C);
        \draw[green!50!black, thick, decorate, decoration={snake, amplitude=.4mm, segment length=3mm, post length = 1 mm}] (C) .. controls (1.5,1.5) and (-0.15,.8) .. (Y);

        \fill (Y)  circle (1pt) node[below] {$Y_t$};
        \fill (Xi) circle (1pt) node[above] {$X^i_{t-\gamma_i}$};
        \fill (Xj) circle (1pt) node[above right] {$X^j_{t-\gamma_j}$};
        \fill (D)  circle (1pt) node[above left] {\small $D_{t_d}$};
        \fill (C)  circle (1pt) node[above right] {\small $C_{t_c}$};
    \end{tikzpicture}
    \begin{tikzpicture}[overlay, shift={(-.9,6.2)}, scale=0.6]
        \coordinate (X) at (1,1);
        \fill (X)  circle (3pt) node[above] {$X_{t_x}$};
        \draw[dashed] (X) -- (0,1);
        \draw[dashed] (0,1) -- (0,0);
        \draw[dashed] (X) -- (2,1);
        \draw[dashed] (2,1) -- (2,0);
        
        \node[anchor=west] at (2.1,0.5) {\smaller$ \bigcup_{\Gf \in \mathcal{C}(\Gs)} \Desc(X_{t_x}, \Gf)$};
    
        
        \draw[green!50!black, thick, decorate, decoration={snake, amplitude=.4mm, segment length=3mm}] (0.5,-.95) -- (1.5,-.95);
        \node[anchor= west] at (2.1,-.95) {\smaller Path \(\pi^f\) leaving \(\CD\)};

        \node at (1, -1.95) {\small $C_{t_c}$};
        \node[anchor= west] at (2.1,-1.95) {\smaller First vertex of \(\pi^f\) in \(\mathcal{V}^f \backslash\CD\)};

        \node at (1, -2.95) {\small $D_{t_d}$};
        \node[anchor= west] at (2.1,-2.95) {\smaller Predecessor of \(C\) in \(\pi^f\)};
        
    \end{tikzpicture}
    
    \caption{Proof idea of Theorem \ref{th:equiv_IBC_multivarie}}\label{fig:CD_proof}
\end{figure}

\begin{proof}
    \textbf{We start by proving that all backdoor paths from an intervention to $Y_t$ that leave $\CD$ are blocked by $\C \coloneqq \mathcal{V}^f \backslash \CD$.} Let us consider an FTCG $\Gf$ and a backdoor path from $X^i_{t-\gamma_i}$ to $Y_t$  in $\Gf_i$ that leaves $\CD$.  Let $C_{t_c}$ be the first vertex of $\pi^f$ in $\mathcal{V}^f \backslash \CD$ and $D_{t_d}$ be the predecessor of $C_{t_c}$ in $\pi^f$.\footnote{$t_d$ may be different from $t_c$.} This construction is presented in Figure \ref{fig:CD_proof}. Necessarily $\pi^f_{\mid \{C_{t_c}, D_{t_d}\}}  = D_{t_d} \leftarrow C_{t_c}$. Indeed, otherwise $C_{t_c}$ would be in $\CD$. Hence, $C_{t_c}$ is not a collider on $\pi^f$, thus $\pi^f$ is blocked by $\{ C_{t_c}\}$. Therefore $\C = \mathcal{V}^f \backslash \CD$ blocks $\pi^f$.

    \noindent \textbf{Let us prove the two implications of the theorem:}
    \begin{itemize}
        \item $\ref{th:equiv_IBC_multivarie:1} \Rightarrow \ref{th:equiv_IBC_multivarie:2}$: Let us prove the contrapositive: Let $\Gf$ be an FTCG such that $\Gf_i$ contains a collider-free backdoor path from $X^i_{t-\gamma_i}$ to $Y_t$ that remains within $\CD$. This path does not have a collider, so the only way to block it is by conditioning on one of its vertices. However, all its vertices are within $\CD$, meaning they are descendants of an intervention in at least one FTCG. Thus, the path cannot be blocked. \textbf{Therefore, the effect is not identifiable by common backdoor.}
        
        \item $\ref{th:equiv_IBC_multivarie:2} \Rightarrow \ref{th:equiv_IBC_multivarie:1}$: Let us assume that condition \ref{th:equiv_IBC_multivarie:2} holds. We will show that $\C = \mathcal{V}^f \backslash \CD$ is a valid common backdoor set:
        \begin{itemize}
            \item By definition, \textbf{$\C$ does not contain any descendant of any intervention in any FTCG.}
            \item Let $\pi^f$ be a backdoor path from an $X^i_{t-\gamma_i}$ to $Y_t$ in some $\Gf_i$. By condition \ref{th:equiv_IBC_multivarie:2}, $\pi^f$ leaves $\CD$. Therefore \textbf{$\C$ blocks $\pi^f$.}
        \end{itemize}
        Therefore, \textbf{$\C = \mathcal{V}^f \backslash \CD$ is a valid common backdoor set}.
    \end{itemize}
\end{proof}

\subsection{Proofs of Section \ref{ssct:in practice}}
\subsubsection{Proofs of Section \ref{subsubsection:CfbWF}}

\mylemmaIBCenumerationcheminsdiriges*

\begin{proof}
    \noindent \textbf{Let us prove the two implications:}
    \begin{itemize}
        \item $\ref{lemma:IBC_enumeration_chemins_diriges:1} \Rightarrow \ref{lemma:IBC_enumeration_chemins_diriges:2}$: Let us consider $X^i_{t-\gamma_i}$ and an FTCG $\Gf$ belonging to $\C(\Gs)$  such that the mutilated graph  $\Gf_i$ contains  $ \pi^f \coloneqq X^i_{t- \gamma_i} \leftsquigarrow Y_t$ which remains in $\CD$. By definition, we already know that $t- \gamma_i \leq t$. Causality does not move backwards in time thus $t- \gamma_i \geq t$. Therefore, $t-\gamma_i = t$ and all vertices of $\pi^f$ are a time $t$. Therefore the reduction $\pi^s$ of $\pi^f$ is still a path in $\Gs$. Moreover, $\pi^f$ remains in $\CD$ thus $\pi^s \subseteq \mathcal{S}$. Therefore, \textbf{$X^i_{t-\gamma_i} \in \Desc(Y, \Gs_{\mid \mathcal{S}})$ and the intervention is at time t.}
        
        \item $\ref{lemma:IBC_enumeration_chemins_diriges:2} \Rightarrow \ref{lemma:IBC_enumeration_chemins_diriges:1}$: Let $X^i_t$ be the intervention at time $t$ such that $X^i \in \Desc(Y, \Gs_{\mid \mathcal{S}})$. Let us consider $\pi^s = X^i \leftarrow V^2 \leftarrow \cdots \leftarrow V^{n-1} \leftarrow Y$, the smallest path from $Y$ to an interventional time series in $\Gs_{\mid \mathcal{S}}$. For all $i \in \{2, \cdots, n-1 \}$, $V^i_t$ is not an intervention because otherwise we could find a smaller path that $\pi^f$. We can construct an FTCG containing the path $\pi^f = X^i_t \leftarrow V^2_t \leftarrow \cdots \leftarrow V^{n-1}_t \leftarrow Y_t$. $\pi^f$ remains in $\CD$ because for all $i$, $V^i$ belongs to $\mathcal{S}$. Therefore, \textbf{There exists a $X^i_{t-\gamma_i}$ and an FTCG $\Gf$ belonging to $\C(\Gs)$  such that the mutilated graph  $\Gf_i$ contains  $X^i_{t- \gamma_i} \leftsquigarrow Y_t$ which remains in $\CD$.}
    \end{itemize}
\end{proof}

\subsubsection{Proofs of Section \ref{subsubsection:CfbF}}

\mylemmaequivexistencecheminforkCDwithoutconsistencythroughtime*

\begin{proof}
    \noindent \textbf{Let us prove the two implications of the theorem:}
    \begin{itemize}
        \item $\ref{lemma:equiv_existence_chemin_fork_CD_without_consistency_through_time:1}
        \Rightarrow \ref{lemma:equiv_existence_chemin_fork_CD_without_consistency_through_time:2}$:
        Let $X^i_{t-\gamma_i}$  be the intervention and $\Gf$ the FTCG belonging to $\C(\Gs)$  such that the mutilated graph  $\Gf_i$ contains the path $X^i_{t- \gamma_i} \leftsquigarrow F_{t_f} \rightsquigarrow Y_t$ which remains in $\CD$. $\Gf$ proves that \textbf{$F_{t_f}$ is $X^i_{t- \gamma_i}$-$\NC$-accessible and $Y_t$-$\NC$-accessible.}
        \item $\ref{lemma:equiv_existence_chemin_fork_CD_without_consistency_through_time:2}
        \Rightarrow \ref{lemma:equiv_existence_chemin_fork_CD_without_consistency_through_time:1}$:
        Let $X^i_{t- \gamma_i}$ and $F_{t_f}$ be such that $F_{t_f}$ is $X^i_{t- \gamma_i}$-$\NC$-accessible and $Y_t$-$\NC$-accessible. Therefore there exists $\Gf_1$ in which there is $\pi^f_1 \vcentcolon= X^i_{t- \gamma_i} \leftsquigarrow F_{t_f}$ which remains in $\NC$ except for $X^i_{t- \gamma_i}$ and there exists $\Gf_2$ in which there is $\pi^f_2 \vcentcolon= F_{t_f} \rightsquigarrow Y_t$ which remains in $\NC$. We distinguish two cases:
        \begin{itemize}
            \item If $\pi^f_1 \cap \pi^f_2 = \{ F_{t_f}\}$, then we can build an FTCG $\Gf_3 \in \C(\Gs)$ which contains $\pi^f_3 \coloneqq X^i_{t- \gamma_i} \leftsquigarrow F_{t_f} \rightsquigarrow Y_t$ the concatenation of $\pi^f_1$ and $\pi^f_2$. $\pi^f_3$ is a path without cycle because $\pi^f_1 \cap \pi^f_2 = \{ F_{t_f}\}$. It remains in $\NC$ because its vertices come from $\pi^f_1$ and $\pi^f_2$ which remain in $\NC$. Therefore $\Gf_{3i}$ contains the path $X^i_{t- \gamma_i} \leftsquigarrow F_{t_f} \rightsquigarrow Y_t$ which remains in $\CD$.

            \item Otherwise, $\pi^f_1 \cap \pi^f_2$ contains at least two element. We know that $X^i_{t-\gamma_i} \notin \pi^f_1 \cap \pi^f_2$ because $\pi^f_2$ remains in $\NC$. Similarly, we know that $Y_t \notin \pi^f_1 \cap \pi^f_2$ because there would be a directed path from $Y_t$ to $X^i_{t-\gamma_i}$. Let us consider $V_{t_v}$ the latest element of $\pi^f_1$ in $\pi^f_1 \cap \pi^f_2$. In $\Gf_1$ there is ${\pi^f_1}' \coloneqq X^i_{t- \gamma_i} \leftsquigarrow V_{t_v}$, in $\Gf_2$ there is ${\pi^f_2}' \coloneqq  V_{t_v} \rightsquigarrow  Y_t$ and ${\pi^f_1}' \cap {\pi^f_2}' = \{ V_{t_v}\}$.Therefore, with the same reasoning as in the first case, we can construct $\Gf_3$ such that $\Gf_{3i}$ contains the path $X^i_{t - \gamma_i} \leftsquigarrow V_{t_v} \rightsquigarrow Y_t$ which remains in $\CD$.
        \end{itemize}
        In all cases, \textbf{there exists a $X^i_{t-\gamma_i}$ and an FTCG $\Gf$ belonging to $\C(\Gs)$  such that the mutilated graph $\Gf_i$ contains the path $X^i_{t- \gamma_i} \leftsquigarrow F_{t_f} \rightsquigarrow Y_t$ which remains in $\CD$.}
    \end{itemize}
\end{proof}

\mylemmacalculVEacc*

\begin{proof}
    \noindent Let $V_{t_v}$ be a temporal variable.\newline
    
    \textbf{Let us prove that Algorithm \ref{algo:calcul_V_E_acc} terminates:}

    \noindent Algorithm \ref{algo:calcul_V_E_acc} \textbf{terminates} because at each step of the \textbf{while} loop, (number of unseen times series, length of $Q$) is strictly decreasing with respect to the {lexicographic ordering}.\newline

    \textbf{Let us prove that Algorithm \ref{algo:calcul_V_E_acc} is correct:}

             Firstly, by induction, \textbf{we show that the algorithm computes the correct value for each seen time series:}
        \begin{itemize}
            \item At the first step of the while loop, for each parent $P$ of $V$ in $\Gs$, the algorithm will compute $t^{\NC}_{V_{t_v}}(P) \gets \max \{t_1 \mid t_1 \leq t_v \text{ and } P_{t_1} \in \NC \backslash \{V_{t_v}\}\}$.\newline
            
            Let $P$ be a parent of $V$ in $\Gs$. We will show that $\{ t_1 \mid P_{t_1} \text{is $V_{t_v}$-$\NC$-accessible} \}  = \{t_1 \mid t_1 \leq t_v \text{ and } P_{t_1} \in \NC \backslash \{V_{t_v}\}\}$ by showing the two inclusions:
            \begin{itemize}
                \item For all $t_1 \in \{t_1 \mid t_1 \leq t_v \text{ and } P_{t_1} \in \NC \backslash \{V_{t_v}\}\}$, we can construct an FTCG $\Gf_{t_1} \in \C(\Gs)$ which contains $P_{t_1} \rightarrow V_{t_v}$. Therefore, all $P_{t_1} \in \{ P_{t_1} \mid t_1 \leq t_v \text{ and } P_{t_1} \in \NC \backslash \{V_{t_v}\}\}$ are $V_{t_v}$-$\NC$-accessible, \textit{i.e.}, $\{t_1 \mid t_1 \leq t_v \text{ and } P_{t_1} \in \NC \backslash \{V_{t_v}\}\} \subseteq \{ t_1 \mid P_{t_1} \text{is $V_{t_v}$-$\NC$-accessible} \}$

                \item  Let $t_1$ be such that $P_{t_1}$ is $V_{t_v}$-$\NC$-accessible. Causality does not move backwards in time thus $t_1 \leq t_v$. Moreover, $P_{t_1} \in \NC \backslash \{V_{t_v}\}$. Indeed, $P_{t_1}$ is $V_{t_v}$-$\NC$-accessible thus $P_{t_1} \in \NC$ and $P_{t_1} \neq V_{t_v}$ because no FTCG can contain a self loop. Therefore $ \{ t_1 \mid P_{t_1} \text{is $V_{t_v}$-$\NC$-accessible} \} \subseteq \{t_1 \mid t_1 \leq t_v \text{ and } P_{t_1} \in \NC \backslash \{V_{t_v}\}\}$ 
            \end{itemize}
            Therefore $\{ t_1 \mid F_{t_1} \text{is $V_{t_v}$-$\NC$-accessible} \}  = \{t_1 \mid t_1 \leq t_v \text{ and } P_{t_1} \in \NC \backslash \{V_{t_v}\}\}$, thus $t^{\NC}_{V_{t_v}}(P) =\max \{t_1 \mid t_1 \leq t_v \text{ and } P_{t_1} \in \NC \backslash \{V_{t_v}\}\}$.\newline
            
            \textbf{Therefore the algorithm computes correct values at the first step of the loop.}

            \item Let us suppose that the algorithm is correct until the $(n-1)$-th loop step. It pops $S_{t_s} = S_{t^{\NC}_{V_{t_v}}(S)}$ with $t^{\NC}_{V_{t_v}}(S)$ being the element of $Q$ with maximum time index. Let $P$ be an unseen parent of $S$ in $\Gs$, the algorithm will compute $t^{\NC}_{V_{t_v}}(P) \gets \max \{t_1 \mid t_1 \leq t_s \text{ and } P_{t_1} \in \NC \backslash \{V_{t_v}\}\}$.  We will show that $\{ t_1 \mid P_{t_1} \text{is $V_{t_v}$-$\NC$-accessible} \}  = \{t_1 \mid t_1 \leq t_v \text{ and } P_{t_1} \in \NC \backslash \{V_{t_v}\}\}$ by showing the two inclusions:
            \begin{itemize}
                \item $S_{t^{\NC}_{V_{t_v}}(S)}$ is $\NC$-accessible thus there exists an FTCG $\Gf \in \C(\Gs)$ which contains $S_{t^{\NC}_{V_{t_v}}(S)} \rightsquigarrow V_{t_v}$. For all $t_1 \in \{t_1 \mid t_1 \leq t_s \text{ and } P_{t_1} \in \NC \backslash \{V_{t_v}\}\}$, we can construct an FTCG $\Gf_{t_1} \in \C(\Gs)$ which contains $P_{t_1} \rightarrow S_{t^{\NC}_{V_{t_v}}(S)} \rightsquigarrow V_{t_v}$. Therefore, all $P_{t_1} \in \{ P_{t_1} \mid t_1 \leq t_s \text{ and } P_{t_1} \in \NC \backslash \{V_{t_v}\}\}$ are $V_{t_v}$-$\NC$-accessible, \textit{i.e.}, $\{t_1 \mid t_1 \leq t_v \text{ and } P_{t_1} \in \NC \backslash \{V_{t_v}\}\} \subseteq \{ t_1 \mid P_{t_1} \text{is $V_{t_v}$-$\NC$-accessible} \}$

                \item  Let $t_1$ be such that $P_{t_1}$ is $V_{t_v}$-$\NC$-accessible. Thus, $P_{t_1} \in \NC \backslash \{V_{t_v}\}$. Let us show that $t_1 \leq t_s$: there exists an FTCG $\Gf$ which contains $\pi^f \coloneqq P_{t_1} \rightsquigarrow V_{t_v}$ which remains in $\NC$. Let $U_{t_u}$ be the successor of $P_{t_1}$ in $\pi^f$. Causality does not move backwards in time thus $t_1 \leq t_u$. We distinguish two cases:
                \begin{itemize}
                    \item If $U = S$, then we already have $t_1 \leq t_s$.

                    \item Otherwise, $t_u \leq t_s$. Indeed, otherwise, thanks to the priority queue, U would have been seen before $S$ and $P$ would not be an unseen vertex. Therefore $t_1 \leq t_u \leq t_s$.
                \end{itemize}
                In all cases $t_1 \leq t_s$. Therefore, $ \{ t_1 \mid P_{t_1} \text{is $V_{t_v}$-$\NC$-accessible} \} \subseteq \{t_1 \mid t_1 \leq t_s \text{ and } P_{t_1} \in \NC \backslash \{V_{t_v}\}\}$.
            \end{itemize}
            Therefore $\{ t_1 \mid P_{t_1} \text{is $V_{t_v}$-$\NC$-accessible} \}  = \{t_1 \mid t_1 \leq t_s \text{ and } P_{t_1} \in \NC \backslash \{V_{t_v}\}\}$, thus $t^{\NC}_{V_{t_v}}(P) =\max \{t_1 \mid t_1 \leq t_v \text{ and } P_{t_1} \in \NC \backslash \{V_{t_v}\}\}$.\newline
            
            \textbf{Therefore the algorithm computes correct values at the $n$-th step of the loop.}
        \end{itemize}
        \textbf{Therefore, by induction principle, the algorithm computes the correct value for each seen time series.}

        Secondly, \textbf{we show that the algorithm sees all time series $S$ such that $t^{\NC}_{V_{t_v}}(S) \neq -\infty$:}

        Let $S$ be a time series such that $t^{\NC}_{V_{t_v}}(S) \neq -\infty$. By definition, $S_{t^{\NC}_{V_{t_v}}(S)}$ is $V_{t_v}$-$\NC$-accessible. Therefore there exists an FTCG in which there exists a path $\pi^f \coloneqq S_{t^{\NC}_{V_{t_v}}(S)} \rightsquigarrow V_{t_v}$ which remains in $\NC$. By cutting unnecessary parts and using Lemma \ref{lemma:helper:construct_FTCG_with_path} or \ref{lemma:helper:construct_FTCG_with_path_with_cst}, we can assume that $\pi^f$ does not pass twice by the same time series except perhaps for $V$ if $S = V$. Thus there is a directed path $\pi^S$ from $S$ to $V$ in $\Gs$ in which all vertices $U$ have a finite $t^{\NC}_{V_{t_v}}(U)$. Since, if the algorithm sees $U$ and $t^{\NC}_{V_{t_v}}(U) \neq -\infty$ then all parents of $U$ are seen by the algorithm and since all the parents of $V$ are seen by the algorithm, by induction on $\pi^s$, we can conclude that $S$ is seen.

        Therefore, \textbf{the algorithm sees all time series $S$ such that $t^{\NC}_{V_{t_v}}(S) \neq -\infty$}.

        Finally, \textbf{we show that Algorithm \ref{algo:calcul_V_E_acc} is correct:}
        Let $S$ be a time series. We distinguish two cases:
        \begin{itemize}
            \item If $S$ is seen by the algorithm, then the correct value is computed by the algorithm.
            \item Otherwise, the algorithm computes $t^{\NC}_{V_{t_v}}(S) \gets - \infty$ thanks to the initialisation step, which is the correct value.
        \end{itemize}

        \textbf{Therefore, Algorithm \ref{algo:calcul_V_E_acc} is correct.}\newline
    
    \noindent \textbf{Let us prove that Algorithm \ref{algo:calcul_V_E_acc} has a complexity of  $\mathcal{O} \left(\left| \mathcal{E}^s \right| + \left| \mathcal{V}^s \right| \cdot \log  \left| \mathcal{V}^s \right|\right)$:}

    $t^{\NC}_{V_{t_v}}(P) \gets \max \{t_1 \mid t_1 \leq t_s \text{ and } P_{t_1} \in \NC \backslash \{V_{t_v}\}\}$ can be computed in $\mathcal{O}(1)$ during the algorithm. Indeed, to do so we need to run first an algorithm similar to Algorithm \ref{algo:calcul_t_NC:1}.

    Let us assume that the priority queue is implemented using a Fibonacci heap. In this case, inserting an element takes $\mathcal{O}(1)$ amortized time, and extracting the max element takes $\mathcal{O}(\log \mid Q \mid)$ amortized time. During the execution of the algorithm, the priority queue $Q$ contains at most $\left| \mathcal{V}^s \right|$ elements, where $\left| \mathcal{V}^s \right|$ is the number of time series (or nodes) in $\Gs$. Therefore, each extraction costs $\mathcal{O}(\log \left| \mathcal{V}^s \right|)$, and since there are at most $\left| \mathcal{V}^s \right|$ extractions, the total cost of all extractions is $\mathcal{O}(\left| \mathcal{V}^s \right| \log \left| \mathcal{V}^s \right|)$.
    Additionally, each edge in $\Gs$ is processed at most once during the computation of the unseen parents of $S$ in $\Gs$ and the update of $t^{\NC}{V{t_v}}(P)$. These computations are performed in $\mathcal{O}(1)$ time. Thus, the total cost of processing all edges is $\mathcal{O}(\left| \mathcal{E}^s \right|)$.

    \textbf{Therefore, the overall complexity of the algorithm is $\mathcal{O}(\left| \mathcal{E}^s \right| + \left| \mathcal{V}^s \right| \log \left| \mathcal{V}^s \right|)$.}
    
\end{proof}

\mylemmaCharactVNCAcc*

\begin{proof}
    \noindent \textbf{Let us prove the two implications:}
    \begin{itemize}
        \item $\ref{lemma:CharactVNCAcc:1} \Rightarrow \ref{lemma:CharactVNCAcc:2}$: Let $F_{t_f}$ and $V_{t_v}$ be such that $F_{t_f}$ is $V_{t_v}$-$\NC$-accessible. Thus $F_{t_f} \in \NC$. Hence, $F_{t_f} \notin \{ X^i_{t-\gamma_i}\}_i$ and $t_{\NC}(F) \leq t_f$ by definition of $t_{\NC}(F)$. By definition of $t^{\NC}_{V_{t_v}}(F)$, $ t_f \leq t^{\NC}_{V_{t_v}}(F)$. Therefore $t_{\NC}(F) \leq t_f \leq t^{\NC}_{V_{t_v}}(F)$ and $F_{t_f} \notin \{ X^i_{t-\gamma_i}\}_i$.

        \item $\ref{lemma:CharactVNCAcc:1} \Rightarrow \ref{lemma:CharactVNCAcc:2}$: Let $F_{t_f}$ and $V_{t_v}$ be such that $t_{\NC}(F) \leq t_f \leq t^{\NC}_{V_{t_v}}(F)$ and $F_{t_f} \notin \{ X^i_{t-\gamma_i} \}_i$. Since $t_{\NC}(F) \leq t_f$, it follows that $t_{\NC}(F) \neq +\infty$. By Lemma \ref{lemma:def_equiv_CD}, we have $\NC = \bigcup_{Z \in \mathcal{V}^S} \{Z_{t_1}\}_{t_1 \geq t_{\NC}(Z)} \backslash \{ X^i_{t-\gamma_i} \}_i$, and given that $F_{t_f} \notin \{ X^i_{t-\gamma_i} \}_i$, it follows that $F_{t_f} \in \NC$. Since $t_f \leq t^{\NC}_{V_{t_v}}(F)$, it follows that $t^{\NC}_{V_{t_v}}(F) \neq -\infty$. Thus, $F_{t^{\NC}_{V_{t_v}}(F)}$ is $V_{t_v}$-$\NC$-accessible. Hence, there exists an FTCG in which there is a path $ \pi^f: V_{t_v} \leftsquigarrow F_{t^{\NC}_{V_{t_v}}(F)}$ that remains in $\NC$. By changing the last arrow of this path, we can construct an FTCG which contains the path $V_{t_v} \leftsquigarrow F_{t_f}$. This path remains in $\NC$ because $F_{t_f} \in \NC$ and all other vertices of the paths are in $\pi^f$ which remains in $\NC$. Therefore $F_{t_f}$ is $V_{t_v}$-$\NC$-accessible.

    \end{itemize}
\end{proof}

\mycorUn*

\begin{proof}
    Let us show the two implications of the corollary:
    \begin{itemize}
        \item $\ref{cor:1:1} \Rightarrow \ref{cor:1:2}$: If there exists $t_f$ such that $F_{t_f}$ is $X^i_{t - \gamma_i}$-$\NC$-accessible and $Y_t$-$\NC$-accessible then $F_{t_{\NC}(F)}$ is $X^i_{t - \gamma_i}$-$\NC$-accessible and $Y_t$-$\NC$-accessible. Thus, by Lemma \ref{lemma:CharactVNCAcc}, $t_{\NC}(F) \leq t_{\NC}(F) \leq t^{\NC}_{X^i_{t - \gamma_i}}(F)$ and $F_{t_{\NC}(F)} \notin \{ X^i_{t-\gamma_i}\}_i$ and $t_{\NC}(F) \leq t_{\NC}(F) \leq t^{\NC}_{Y_t}(F)$ and $F_{t_{\NC}(F)} \notin \{ X^i_{t-\gamma_i}\}_i$. Hence $t_{\NC}(F) \leq t^{\NC}_{X^i_{t - \gamma_i}}(F)$ and $t_{\NC}(F) \leq t^{\NC}_{Y_t}(F)$.

        \item $\ref{cor:1:2} \Rightarrow \ref{cor:1:1}$: If $t_{\NC}(F) \leq t^{\NC}_{X^i_{t - \gamma_i}}(F)$ and $t_{\NC}(F) \leq t^{\NC}_{Y_t}(F)$, then $t_{\NC}(F) \neq + \infty$ and $t^{\NC}_{X^i_{t - \gamma_i}}(F) \neq -\infty$. Therefore, $F_{t_{\NC}(F)} \notin \{ X^i_{t-\gamma_i}\}_i$ and  $t_{\NC}(F) \leq t_{\NC}(F) \leq t^{\NC}_{X^i_{t - \gamma_i}}(F)$ and $F_{t_{\NC}(F)} \notin \{ X^i_{t-\gamma_i}\}_i$ and $t_{\NC}(F) \leq t_{\NC}(F) \leq t^{\NC}_{Y_t}(F)$. Therefore, by Lemma \ref{lemma:CharactVNCAcc}, $F_{t_{\NC}(F)}$ is $X^i_{t - \gamma_i}$-$\NC$-accessible and $Y_t$-$\NC$-accessible. Therefore, there exists $t_f$ such that $F_{t_f}$ is $X^i_{t - \gamma_i}$-$\NC$-accessible and $Y_t$-$\NC$-accessible.
    \end{itemize}
    
\end{proof}

\subsubsection{Proofs of Section \ref{sssct:algo_IBC}}

\mythforalgoIBC*

\begin{proof}
    Algorithm \ref{algo:calcul_IBC} uses directly the characterizations of Lemma \ref{lemma:IBC_enumeration_chemins_diriges} and Lemma \ref{lemma:equiv_existence_chemin_fork_CD_without_consistency_through_time}. It outputs $False$ if and only if there exists an FTCG $\Gf$ in which there is a backdoor path from $X^i_{t -\gamma_i}$ to $Y_t$, otherwise it outputs $True$. Therefore, by Theorem \ref{th:equiv_IBC_multivarie}, Algorithm \ref{algo:calcul_IBC} is correct.

    Algorithm \ref{algo:calcul_IBC} tests the existence of directed paths in $\mathcal{O}(\left| \mathcal{V}^s \right| + \left| \mathcal{E}^s \right|)$, runs Algorithm \ref{algo:calcul_t_NC} in $\mathcal{O}(\left| \mathcal{V}^s \right| + \left| \mathcal{E}^s \right|)$, calls Algorithm \ref{algo:calcul_V_E_acc} $n+1$ times in $\mathcal{O}(n(\left| \mathcal{E}^s \right| + \left| \mathcal{V}^s \right| \log \left| \mathcal{V}^s \right|))$ and does $\mathcal{O}(n \cdot \left| \mathcal{V}^s \right|)$. Therefore, it overall complexity is $\mathcal{O}\left(n( \left| \mathcal{E}^s \right| + \left| \mathcal{V}^s \right| \log \left| \mathcal{V}^s \right|)\right)$.
\end{proof}

\newpage
\section{Proofs of Section \ref{sec:Consistency_Time}}\label{sec:proof:5}

\mylemmaequivexistencecheminforkNC*

\begin{proof}
    Let us prove the direct implication (\ref{lemma:equiv_existence_chemin_fork_NC:1} $\Rightarrow$ \ref{lemma:equiv_existence_chemin_fork_NC:2}). Let $F$ and $X^i_{t - \gamma_i}$ be such that there exists an FTCG $\Gf$ belonging to $\C(\Gs)$ such that $\Gf_i$ contains the path $X^i_{t- \gamma_i} \leftsquigarrow F_{t'} \rightsquigarrow Y_t$ which remains in $\CD$. Two cases arise:
    \begin{itemize}
        \item If $F \neq Y$ or $t - \gamma_i \neq t_{\NC}(F)$, the accessibility in $\Gf$ ensures that $F_{t_{\NC}(F)}$ is both $X^i_{t - \gamma_i}$-$\NC$-accessible and $Y_t$-$\NC$-accessible. Thus, we have proven Proposition \ref{lemma:equiv_existence_chemin_fork_NC:2a} in this case.
    
        \item Otherwise, $Y = F$ and $t - \gamma_i = t_{\NC}(F) = t_{\NC}(Y)$. Let us prove by contradiction that Lemma \ref{lemma:equiv_existence_chemin_fork_NC} \ref{lemma:equiv_existence_chemin_fork_NC:2b} holds. If it does not hold, then $\Gf$ must necessarily contain the edges $X^i_{t - \gamma_i} \leftarrow Y_{t_{\NC}(Y)}$ and $X^i_t \rightarrow Y_t$, which contradicts Assumption \ref{ass:Consistency_Time}.
    \end{itemize}

    \noindent Conversely, let us prove the indirect direction (\ref{lemma:equiv_existence_chemin_fork_NC:2} $\Rightarrow$ \ref{lemma:equiv_existence_chemin_fork_NC:1}). To prove this implication, we must show that  Lemma \ref{lemma:equiv_existence_chemin_fork_NC} \ref{lemma:equiv_existence_chemin_fork_NC:2a} implies  Lemma \ref{lemma:equiv_existence_chemin_fork_NC}  \ref{lemma:equiv_existence_chemin_fork_NC:1}, and that  Lemma \ref{lemma:equiv_existence_chemin_fork_NC}  \ref{lemma:equiv_existence_chemin_fork_NC:2b} implies  Lemma \ref{lemma:equiv_existence_chemin_fork_NC}  \ref{lemma:equiv_existence_chemin_fork_NC:1}. For each of these proofs, we are given $\Gf_1$, which contains $\pi^f_1 \vcentcolon= X^i_{t- \gamma_i} \leftsquigarrow F_{t_{\NC}(F)}$, which remains in $\NC$ except for $X^i_{t- \gamma_i}$, and $\Gf_2$, which contains $\pi^f_2 \vcentcolon= F_{t_{\NC}(F)} \rightsquigarrow Y_t$, which remains in $\NC$. To prove Proposition \ref{lemma:equiv_existence_chemin_fork_NC:1}, it suffices to construct $\Gf_3$ such that $\Gf_{3i}$ contains $\pi^f_3 \vcentcolon= X^i_{t- \gamma_i} \leftsquigarrow F_{t'} \rightsquigarrow Y_t$, which remains in $\NC$. We observe that if $\pi^f_3$ passes only through vertices visited by $\pi^f_1$ or $\pi^f_2$, then  its existence in $\Gf_{3i}$ is equivalent to its existence in $\Gf_{3}$.

    \noindent Furthermore, without loss of generality, we can assume that the only intersection between $\pi^f_1$ and $\pi^f_2$ is $F_{t_{\NC}(F)}$. Indeed, consider $V_{t_v}$, the last element of $\pi^f_2$ in $\pi^f_1 \cap \pi^f_2$. By contradiction, we show that $V_{t_v} \neq Y_t$. If $V_{t_v} = Y_t$, then $\Gf_1$ contains $X^i_{t- \gamma_i} \leftsquigarrow Y_t$. However, by assumption, no FTCG contains a backdoor path without a fork from an intervention to $Y_t$. Thus, $V_{t_v} \neq Y_t$, and we can therefore work on ${\pi^f_1}' \vcentcolon= X^i_{t- \gamma_i} \leftsquigarrow V_{t_{\NC}(V)}$ and ${\pi^f_2}' \vcentcolon=V_{t_{\NC}(V)} \rightsquigarrow Y_t$. Therefore, concatenating $\pi^f_1$ and $\pi^f_2$ to create an FTCG does not create a cycle. Only Assumption \ref{ass:Consistency_Time} can be violated.

    If it is possible to concatenate $\pi^f_1$ and $\pi^f_2$ while maintaining Assumption \ref{ass:Consistency_Time}, then we prove Proposition \ref{lemma:equiv_existence_chemin_fork_NC:1}. If this is not the case, we show that, thanks to the assumptions of \ref{lemma:equiv_existence_chemin_fork_NC:2a} or \ref{lemma:equiv_existence_chemin_fork_NC:2b}, we can construct a suitable $\pi^f_3$. We thus merge the proofs of (\ref{lemma:equiv_existence_chemin_fork_NC:2a} $\Rightarrow$ \ref{lemma:equiv_existence_chemin_fork_NC:1}) and (\ref{lemma:equiv_existence_chemin_fork_NC:2b} $\Rightarrow$ \ref{lemma:equiv_existence_chemin_fork_NC:1}) as the reasoning is identical. 
    
    \noindent Therefore, we assume that it is not possible to concatenate $\pi^f_1$ and $\pi^f_2$ because Assumption \ref{ass:Consistency_Time} would be violated.
    We denote by $V^1_{t_1} \rightarrow V^2_{t_1}$ the last arrow of $\pi^f_1$ that contradicts $\pi^f_2$, and $V^1_{t_2} \leftarrow V^2_{t_2}$ the last arrow of $\pi^f_2$ that contradicts $V^1_{t_1} \rightarrow V^2_{t_1}$. We proceed by case distinction:
    \begin{itemize}
        \item If $V^1_{t_1} \neq F_{t_{\NC}(F)}$ and $V^2_{t_1} \neq X^i_{t- \gamma_i}$, we further distinguish three cases:
        \begin{itemize}
            \item If $t_1 < t_2$, then by adding the arrow $V^2_{t_1} \rightarrow V^1_{t_2}$, we can construct $\pi^f_3 = X^i_{t- \gamma_i} \leftsquigarrow V^2_{t_1}  \rightarrow V^1_{t_2} \rightsquigarrow Y_t$ (see Figure \ref{subfig:CheminsForkNC:a}).
            
            \item If $t_1 > t_2$, then by adding the arrow $V^1_{t_2} \rightarrow V^2_{t_1}$, we can construct $\pi^f_3 = X^i_{t- \gamma_i} \leftsquigarrow V^2_{t_1}  \leftarrow V^1_{t_2} \rightsquigarrow Y_t$ (see Figure \ref{subfig:CheminsForkNC:b}).

            \item The case $t_1 = t_2$ is excluded because the only intersection between $\pi^f_1$ and $\pi^f_2$ is $F_{t_{\NC}(F)}$.
        \end{itemize}
        
        \item If $V^1_{t_1} \neq F_{t_{\NC}(F)}$ and $V^2_{t_1} = X^i_{t- \gamma_i}$. $\Gf_2$ contains $V^1_{t_2} \rightsquigarrow Y_t$, which can be transformed into $V^1_{t-\gamma_i} \rightsquigarrow Y_t$ by changing the first arrow, yielding $\pi^f_3 = X^i_{t-\gamma_i} \leftarrow V^1_{t-\gamma_i} \rightsquigarrow Y_t$ (see Figure \ref{subfig:CheminsForkNC:c}).
        
        \item If $V^1_{t_1} = F_{t_{\NC}(F)}$ and $V^2_{t_1} \neq X^i_{t- \gamma_i}$. $\Gf_2$ contains $F_{t_2} \rightsquigarrow Y_t$, which can be transformed into $F_{t_{\NC}(F)}\rightsquigarrow Y_t$ by only changing the first arrow, yielding $\pi^f_3 = X^i_{t-\gamma_i} \leftsquigarrow V^2_{t_{\NC}(F)} \leftarrow F_{t_{\NC}(F)}  \rightsquigarrow Y_t$ (see Figure \ref{subfig:CheminsForkNC:d}).
        
        \item If $V^1_{t_1} = F_{t_{\NC}(F)}$ and $V^2_{t_1} = X^i_{t- \gamma_i}$, we further distinguish two cases:
        \begin{itemize}
            \item If $F \neq Y$, then $t - \gamma_i = t_{\NC(F)}$, because the case $t - \gamma_i \neq t_{\NC(F)}$ does not contradict Assumption \ref{ass:Consistency_Time}. $\Gf_2$ contains $F_{t_2} \rightsquigarrow Y_t$, which can be transformed into $F_{t_{\NC}(F)}\rightsquigarrow Y_t$, yielding $\pi^f_3 = X^i_{t-\gamma_i} \leftarrow V_{t_{\NC}(F)} \leftarrow F_{t_{\NC}(F)}  \rightsquigarrow Y_t$ (see Figure \ref{fig:CheminsForkNC:e}).
            
            \item Otherwise, if $F = Y$, the only case that can contradict Assumption \ref{ass:Consistency_Time} is $t - \gamma_i = t_{\NC(F)}$, i.e., we are strictly under the assumptions of Proposition \ref{lemma:equiv_existence_chemin_fork_NC:2b}. According to Proposition \ref{lemma:equiv_existence_chemin_fork_NC:2b}, $Y_{t_{\NC}(Y)}$ is $Y_t$ and is $\NC$-accessible without using $X^i_t \rightarrow Y_t$. Therefore, we can construct $\pi^f_3$ without contradicting Assumption \ref{ass:Consistency_Time}.
        \end{itemize}
    \end{itemize}

\end{proof}
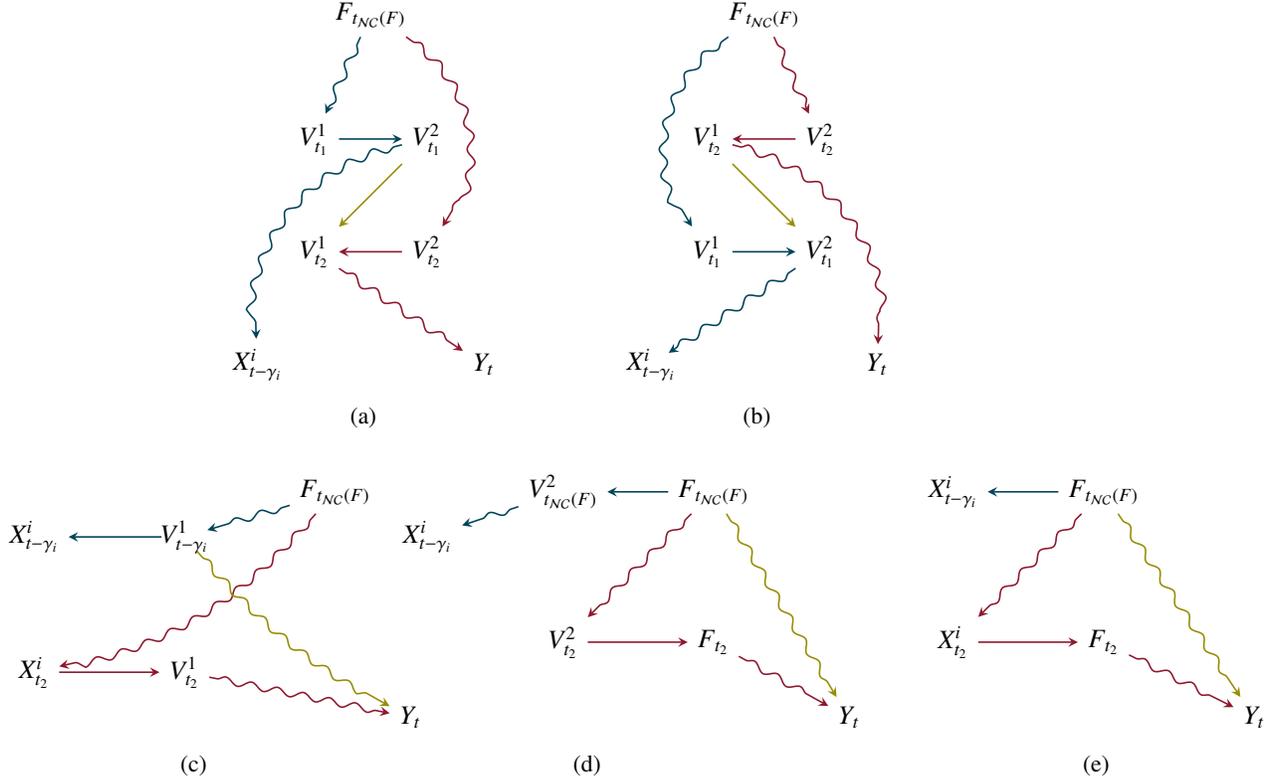
\begin{figure}
    \centering
    \begin{subfigure}{0.3\textwidth}
        \centering
        \begin{tikzpicture}[scale = 1.5, ->,>=stealth,auto,node distance=3cm,semithick]
            \node at (0,0) (V1) {$V^1_{t_2}$};
            \node at (1,0) (V2) {$V^2_{t_2}$};
            \node at (1,1) (V3) {$V^2_{t_1}$};
            \node at (0,1) (V4) {$V^1_{t_1}$};
            \node at (0.5,2.1) (F) {$F_{t_{\NC}(F)}$};
            \node at (-0.5,-1) (X) {$X^i_{t- \gamma_i}$};
            \node at (1.5,-1) (Y) {$Y_t$};

            \path (V4) edge[CentraleBlue] (V3)
                  (F)  edge[CentraleBlue, decorate, decoration={snake, pre length = 3pt, post length=2pt, amplitude=1.5pt}] (V4)
                  (V3) edge[CentraleBlue, bend right = 40, decorate, decoration={snake, pre length = 3pt, post length=2pt, amplitude=1.5pt}] (X)
                  (V2) edge[CentraleRed]  (V1)
                  (F)  edge[CentraleRed, bend left = 45, decorate, decoration={snake, pre length = 3pt, post length=5pt, amplitude=1.5pt}] (V2)
                  (V1) edge[CentraleRed, decorate, decoration={snake, pre length = 3pt, post length=5pt, amplitude=1.5pt}]  (Y)
                  (V3) edge[olive] (V1);
        \end{tikzpicture}
        \caption{\label{subfig:CheminsForkNC:a}}
    \end{subfigure}
    \begin{subfigure}{0.3\textwidth}
        \centering
        \begin{tikzpicture}[scale = 1.5, ->,>=stealth,auto,node distance=3cm,semithick]
            \node at (0,1) (V1) {$V^1_{t_2}$};
            \node at (1,1) (V2) {$V^2_{t_2}$};
            \node at (1,0) (V3) {$V^2_{t_1}$};
            \node at (0,0) (V4) {$V^1_{t_1}$};
            \node at (0.5,2.1) (F) {$F_{t_{\NC}(F)}$};
            \node[outer sep = -5] at (-0.5,-1) (X) {$X^i_{t- \gamma_i}$};
            \node at (1.5,-1) (Y) {$Y_t$};

            \path (V4) edge[CentraleBlue] (V3)
                  (F)  edge[CentraleBlue, bend right = 45, decorate, decoration={snake, pre length = 3pt, post length=2pt, amplitude=1.5pt}] (V4)
                  (V3) edge[CentraleBlue, decorate, decoration={snake, pre length = 3pt, post length=2pt, amplitude=1.5pt}] (X)
                  (V2) edge[CentraleRed]  (V1)
                  (F)  edge[CentraleRed, decorate, decoration={snake, pre length = 3pt, post length=2pt, amplitude=1.5pt}] (V2)
                  (V1) edge[CentraleRed,bend left = 40, decorate, decoration={snake, pre length = 3pt, post length=5pt, amplitude=1.5pt}]  (Y)
                  (V1) edge[olive] (V3);
        \end{tikzpicture}
        \caption{\label{subfig:CheminsForkNC:b}}
    \end{subfigure}
    \newline
    
    \begin{subfigure}{0.3\textwidth}
    \centering
        \begin{tikzpicture}[scale = 2, ->,>=stealth,auto,node distance=3cm,semithick]
            \node[outer sep = -4] at (1,.7) (V1) {$V^1_{t-\gamma_i}$};
            \node at (1,-.2) (V2) {$V^1_{t_2}$};
            \node at (0,.7) (X1) {$X^i_{t-\gamma_i}$};
            \node at (0,-.2) (X2) {$X^i_{t_2}$};
            \node at (2,1) (F) {$F_{t_{\NC}(F)}$};
            \node at (2.5,-.5) (Y) {$Y_t$};

            \path (F)  edge[CentraleBlue, decorate, decoration={snake, pre length = 3pt, post length=5pt, amplitude=1.3pt}] (V1)
                  (V1) edge[CentraleBlue] (X1)
                  (F)  edge[CentraleRed, bend left = 18, decorate, decoration={snake, pre length = 3pt, post length=5pt, amplitude=1.3pt}] (X2)
                  (X2) edge[CentraleRed] (V2)
                  (V2) edge[CentraleRed, decorate, decoration={snake, pre length = 3pt, post length=5pt, amplitude=1.3pt}] (Y)
                  (V1)  edge[olive, bend right = 10, decorate, decoration={snake, pre length = 3pt, post length=5pt, amplitude=1.3pt}] (Y);

        \end{tikzpicture}
        \caption{\label{subfig:CheminsForkNC:c}}
    \end{subfigure}
    \begin{subfigure}{0.3\textwidth}
    \centering
        \begin{tikzpicture}[scale = 2, ->,>=stealth,auto,node distance=3cm,semithick]
            \node at (0,0) (V2) {$V^2_{t_2}$};
            \node at (1,0) (F2) {$F_{t_2}$};
            \node at (1,1) (F1) {$F_{t_{\NC}(F)}$};
            \node at (0,1) (V1) {$V^2_{t_{\NC}(F)}$};
            \node at (-.9,.7) (X) {$X^i_{t-\gamma_i}$};
            \node at (1.9,-.5) (Y) {$Y_t$};

            \path (F1) edge[CentraleBlue] (V1)
                  (V1) edge[CentraleBlue, decorate, decoration={snake, pre length = 3pt, post length=5pt, amplitude=1.3pt}] (X)
                  (F1) edge[CentraleRed, decorate, decoration={snake, pre length = 3pt, post length=5pt, amplitude=1.3pt}] (V2)
                  (V2) edge[CentraleRed] (F2)
                  (F2) edge[CentraleRed, decorate, decoration={snake, pre length = 3pt, post length=5pt, amplitude=1.3pt}] (Y)
                  (F1) edge [olive, decorate, decoration={snake, pre length = 3pt, post length=5pt, amplitude=1.3pt}] (Y);
        \end{tikzpicture}
        \caption{\label{subfig:CheminsForkNC:d}}
    \end{subfigure}
    \hfill
    \begin{subfigure}{0.3\textwidth}
    \centering
    \begin{tikzpicture}[scale = 2, ->,>=stealth,auto,node distance=3cm,semithick]
        \node at (1,1) (F) {$F_{t_{\NC}(F)}$};
        \node at (1,0) (V2) {$F_{t_2}$};
        \node at (0,1) (X1) {$X^i_{t-\gamma_i}$};
        \node at (0,0) (X2) {$X^i_{t_2}$};
        \node at (2,-.5) (Y) {$Y_t$};

        \path (F) edge[CentraleBlue] (X1)
              (F)  edge[CentraleRed,decorate, decoration={snake, pre length = 3pt, post length=5pt, amplitude=1.3pt}] (X2)
              (X2) edge[CentraleRed] (V2)
              (V2) edge[CentraleRed, decorate, decoration={snake, pre length = 3pt, post length=5pt, amplitude=1.3pt}] (Y)
              (F)  edge[olive, decorate, decoration={snake, pre length = 3pt, post length=5pt, amplitude=1.3pt}] (Y);
    \end{tikzpicture}
    \caption{\label{fig:CheminsForkNC:e}}
    \end{subfigure}
    
    \caption{ $\pi^f_1$ is represented in {\color{CentraleBlue} blue}, $\pi^f_2$ is represented in {\color{CentraleRed} red}, and the modification to be made to construct $\pi^f_3$ in $\Gf_3$ is represented in {\color{olive} green}.} \label{fig:CheminsForkNC}
\end{figure}

\mythIBCmonovarie*

\begin{proof}
The proof is the combination of the following theorems.
\end{proof}

\subsection{Proof of Theorem \ref{th:CondNecSuf_IBC_monovarie}: The case \texorpdfstring{$\gamma = 0$}{γ = 0}}

\begin{restatable}{theorem}{mythIBCmonovarieGammaZero}{}
\label{th:CondNecSuf_IBC_monovarie_gamma_zero}

    Let $\Gs$ be an SCG such that $X\in \Anc(Y,\Gs)$ and we consider the total effect $P(y_t\mid \Do(x_{t}))$. The effect is identifiable by common backdoor in $\Gs$ if and only if there is no collider-free backdoor path from $X$ to $Y$ in $\Gs_{\mid \Desc(X, \Gs)}$. In that case, a common backdoor set is given by $\mathcal{A}_{0}$, where:
   $$
       \mathcal{A}_0 \coloneqq \bigcup_{\pi^s \text{ backdoor}} \Bigl\{ Z_t  \mid Z \in  \pi^s \backslash \Desc\left( X,\Gs \right) \Bigr\} \cup \Bigl\{ \left( Z_{t'} \right)_{t' < t} \mid Z \in \Gs \Bigr\}
   $$

\end{restatable}

\begin{proof}
    \noindent By Theorem \ref{th:equiv_IBC_multivarie}, $P(y_t \mid \Do(x_t))$ is identifiable by common backdoor  if and only if for every FTCG $\Gf$ belonging to $\C(\Gs)$, there does not exist a collider-free backdoor path going from $X_t$ to $Y_t$ that remains in $\CD$. To finish the proof, we will show that the existence of an FTCG belonging to $\C(\Gs)$ which contains a collider-free backdoor path from $X_t$ to $Y_t$ is equivalent to the existence of collider-free backdoor path from $X$ to $Y$ in $\Gs$:
    
    \noindent By applying Algorithm \ref{algo:calcul_t_NC} and Lemma \ref{lemma:def_equiv_CD}, we can show that $\CD = \bigcup_{Z \in \Desc(X,\Gs)} \left\{ Z_t\right\}$. Let us show the two implications:
    \begin{itemize}
        \item Let $\Gf$ be an FTCG of $\C(\Gs)$ that contains $\pi^f$, a collider-free backdoor path from $X_t$ to $Y_t$ which remains in $\CD$. Since $ \pi^f \subseteq \CD = \bigcup_{Z \in \Desc(X,\Gs)} \left\{ Z_t\right\}$ all the vertices of $\pi^f$ are at time $t$. Therefore, the reduction $\pi^s$ of $\pi^s$ is still a collider-free backdoor path in $\Gs$.  Since $ \pi^f \subseteq \bigcup_{Z \in \Desc(X,\Gs)} \left\{ Z_t\right\}$, $\pi^s$ is  a collider-free backdoor path in $\Gs_{\mid \Desc(X, \Gs)}$.

        \item Let $\pi^s \coloneqq \left[ V^1, \cdots, V^n \right]$ be a collider-free backdoor path from $X$ to $Y$ in $\Gs_{\mid \Desc(X, \Gs)}$. We can construct an FTCG $\Gf$ belonging to $\C(\Gs)$ which contains $\pi^f \coloneqq \left[ V^1_t, \cdots, V^n_t \right]$. Since $\CD = \bigcup_{Z \in \Desc(X,\Gs)} \left\{ Z_t\right\}$ and $\pi^s \subseteq \Desc(X, \Gs)$, $\pi^f$ is a collider-free backdoor path from $X_t$ to $Y_t$ and it remains in $\CD$.
    \end{itemize}

    \textbf{Therefore, $P(y_t \mid \Do(x_t))$ is identifiable by common backdoor if and only if there is no collider-free backdoor path from $X$ to $Y$ in $\Gs_{\mid \Desc(X, \Gs)}$.}

    \textbf{Let us show that $\mathcal{A}_0$ is a correct adjustment set for the case $\gamma = 0$:}
    \noindent By definition, $\mathcal{A}_0$ is defined as follows:
    \begin{align*}
    \mathcal{A}_0 \vcentcolon&= \mathcal{A}_0^t \cup \mathcal{A}_0^-\\
   \mathcal{A}_0^t \vcentcolon&=  \bigcup_{\pi^s \text{ backdoor}} \Bigl\{ Z_t  \mid Z \in  \pi^s \backslash \Desc\left( X,\Gs \right) \Bigr\}\\
   \mathcal{A}_0^- \vcentcolon&= \Biggl\{ \left( Z_{t'} \right)_{ t' < t} \mid Z \in \Gs \Biggr\}.
    \end{align*}
    
    Let us assume that $P(y_t \mid \Do(x_t))$ is identifiable by common backdoor, then we have:
\begin{itemize}
    \item $\mathcal{A}_0$ contains no descendants of $X_t$ in any FTCG $\Gf \in \mathcal{C} \left( \Gs\right)$: Indeed, $\mathcal{A}_0^-$ cannot contain any descendants of $X_t$ in any FTCG because causality does not move backwards in time. By contradiction, we show that $\mathcal{A}_0^t$ cannot contain any descendants of $X_t$ in any FTCG. Let $\Gf \in \mathcal{C} \left( \Gs\right)$ and suppose $Z_{t} \in \Desc\left( X_t, \Gf \right)$. Thus, $Z \in \Desc\left( X, \Gs \right)$, which contradicts the definition of $\mathcal{A}_0^t$.
    
    \item $\mathcal{A}_0$ blocks all backdoor paths from $X_t$ to $Y_t$ in any FTCG $\Gf \in \C(\Gs)$: Let $\Gf \in \mathcal{C} \left( \Gs\right)$ be an FTCG and $\pi^f$ a backdoor path in $\Gf$. We distinguish two cases:
    \begin{itemize}
        \item If $\pi^f$ passes through $\mathcal{A}_0^-$, then it is blocked by $\mathcal{A}_0^-$.
        \item Otherwise, $\pi^f$ remains at time $t$. We then consider $\pi^s$, its reduction in $\Gs$. Since $\pi^f$ remains at time $t$,we know that $\pi^s$ is a path in $\Gs$. We distinguish two cases:
        \begin{itemize}
            \item If $\pi^s \subseteq \Desc(X, \Gs)$, then $\pi^s$ contains a collider $C$ because $P(y_t \mid \Do(x_t))$ is identifiable by common backdoor and $\Gs_{\mid \Desc(X, \Gs)}$ contains no collider-free backdoor path from $X$ to $Y$. On $\pi^f$, $C_t$ is a blocking collider because $C_t$ is a descendant of $X_t$, and therefore all its descendants are descendants of $X_t$, and $\mathcal{A}_0$ contains no descendants of $X_t$.
            
            \item Otherwise, let $Z$ be the first element of $\pi^s$ that is not a descendant of $X$. $Z$ is not a collider, so $Z_t$ is not a collider on $\pi^f$. However, $Z_t$ belongs to $\mathcal{A}_0^t$, hence $\mathcal{A}_0^t$ blocks $\pi^f$.
        \end{itemize}
    \end{itemize}
\end{itemize}

\textbf{Thus, $\mathcal{A}_0$ is a valid common backdoor set when the effect is identifiable by common backdoor.}
\end{proof}

It is worth noting that the size of $\mathcal{A}_0$ can be reduced. Indeed, as the proof shows, for each $\pi^s$, it is sufficient to take the first element of $\pi^s$ that is not a descendant of $X$, rather than taking $\Bigl\{ Z_t \mid Z \in \pi^s \backslash \Desc\left( X,\Gs \right) \Bigr\}$.

\subsection{Proof of Theorem \ref{th:CondNecSuf_IBC_monovarie}: A necessary condition for the cases \texorpdfstring{$\gamma \geq 1$}{γ >= 1}}

We begin by proving lemmas that provide a necessary condition for identifiability by common backdoor for all cases where $\gamma \geq 1$.

\begin{restatable}{lemma}{mylemmaCycleAvecAncY1}{}
\label{lemma:gammaGeq1:CycleAvecAncY}
    Let $\Gs$ be an SCG such that $X \in \Anc\left( Y, \Gs \right)$ and $P(y_t\mid \Do(x_{t-\gamma})$ is the considered effect for $\gamma \geq 1$. If $\text{Cycles}^>(X, \Gs \backslash\left\{ Y\right\}) \cap \Anc\left( Y, \Gs \right) \neq \emptyset$, then the effect is not identifiable by common backdoor.
\end{restatable}

\begin{proof}
    If $\text{Cycles}^>(X, \Gs \backslash\left\{ Y\right\}) \cap \Anc\left( Y, \Gs \right) \neq \emptyset$, then there exists a node $Z$ such that $\Gs$ contains the structure represented in Figure \ref{subfig:CycleAvecAncY:1}. It is always possible to choose $Z$ such that the paths $X \rightsquigarrow Z$ and $Z \rightsquigarrow Y$ only have $Z$ in common. We construct a backdoor path from $X_{t-\gamma}$ to $Y_t$ without colliders that remains in $\CD$ (cf. \ref{subfig:CycleAvecAncY:2}).

\begin{figure}[H]
    \centering
    \begin{subfigure}{0.4\textwidth}
        \centering
        \begin{tikzpicture}[->,>=stealth',shorten >=1pt,auto,node distance=3cm, semithick]
            \tikzstyle{every node}=[fill=white,draw=none,text=black, minimum height = 21pt]
            \node (X) {$X$};
            \node (Z) [right of=X] {$Z$};
            \node (Y) [right of=Z] {$Y$};
            \node (B) [below of=Z] {}; 
            \path (Z) edge[decorate, decoration={snake, pre length = 3pt, post length=5pt, amplitude=2pt}] (Y)
                  (X) edge[bend left = 35, decorate, decoration={snake, pre length = 3pt, post length=5pt, amplitude=2pt}] (Z)
                  (Z) edge[bend left = 35, decorate, decoration={snake, pre length = 3pt, post length=5pt, amplitude=2pt}] (X);
        \end{tikzpicture}
        \caption{Structure belonging to the SCG.}
        \label{subfig:CycleAvecAncY:1}
    \end{subfigure}
    \hfill
    \begin{subfigure}{0.4\textwidth}
        \centering
        \begin{tikzpicture}[->,>=stealth',shorten >=1pt,auto,node distance=3cm,
                            semithick]
            \tikzstyle{every node}=[fill=white,draw=none,text=black, minimum height = 21pt]

            \node (Xt1) {$X_{t-\gamma}$};
            \node (Zt1) [below of=Xt1] {$Z_{t-\gamma}$};
            \node (Yt1) [below of=Zt1] {$Y_{t-\gamma}$};
            \node (Xt) [right of=Xt1] {$X_t$};
            \node (Zt) [below of=Xt] {$Z_{t}$};
            \node (Yt) [below of=Zt] {$Y_{t}$};

            \path (Zt1) edge[CentraleRed, decorate, decoration={snake, pre length = 3pt, post length=5pt}] (Xt1)
                  (Zt1) edge[CentraleRed, decorate, decoration={snake, pre length = 3pt, post length=5pt}] (Yt);
        \end{tikzpicture}
        \caption{FTCG $\mathcal{G}^f \in \mathcal{C} \left( \mathcal{G}^s \right)$}
        \label{subfig:CycleAvecAncY:2}
    \end{subfigure}
    \caption{Illustration of Lemma \ref{lemma:gammaGeq1:CycleAvecAncY}. Figure \ref{subfig:CycleAvecAncY:1} represents a structure contained in the SCG. The FTCG presented in Figure \ref{subfig:CycleAvecAncY:2} is indeed constructible. In fact, to construct the arrow $Z_{t-\gamma} \rightsquigarrow Y_t$, it is sufficient to place all its arrows at time $t - \gamma$ and use the last one to go to time $t$. Moreover, Assumption \ref{ass:Consistency_Time} is respected since the paths $X \rightsquigarrow Z$ and $Z \rightsquigarrow Y$ only have $Z$ in common. The missing arrows are placed at a greater $\gamma$ and are therefore not represented. The collider-free backdoor path represented in \ref{subfig:CycleAvecAncY:2} remains within $\CD$.
}
    \label{fig:CycleAvecAncY}
\end{figure}
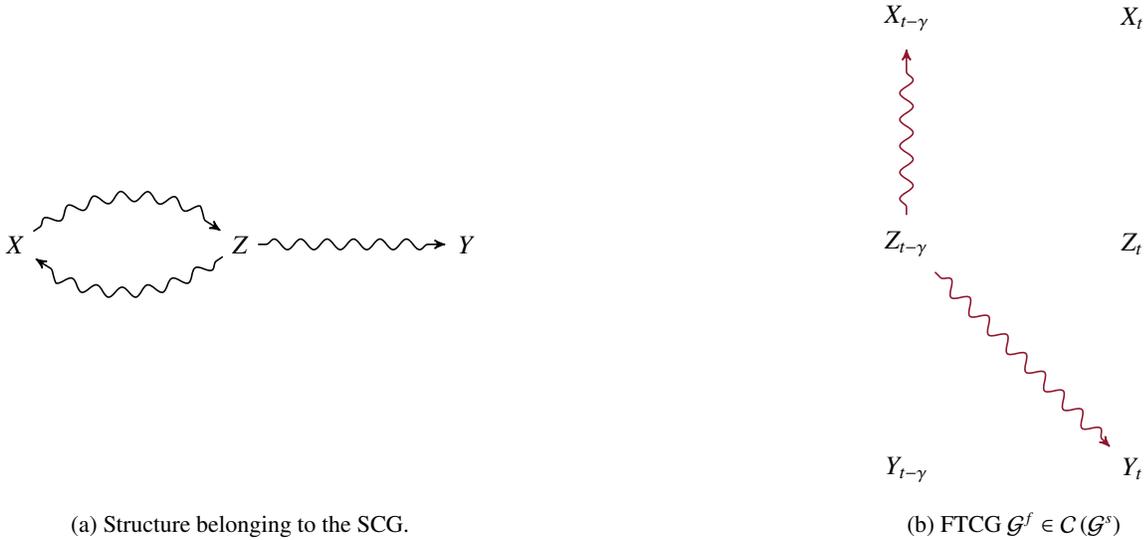
\end{proof}

\begin{restatable}{lemma}{mylemmaCycleSansAncY}{}
\label{lemma:gammaGeq1:CycleSansAncY}
    Let $\Gs$ be an SCG such that $X \in \Anc\left( Y, \Gs \right)$ and $P(y_t\mid \Do(x_{t-\gamma})$ is the considered effect for $\gamma \geq 1$. If $\text{Cycles}^>(X, \Gs \backslash\left\{ Y\right\}) \neq \emptyset$ and $\text{Cycles}^>(X, \Gs \backslash\left\{ Y\right\}) \cap \Anc\left( Y, \Gs \right) = \emptyset$, then the effect is not identifiable by common backdoor.
\end{restatable}

\begin{proof}
If $\text{Cycles}^>(X, \Gs \backslash\left\{ Y\right\}) \neq \emptyset$ and $\text{Cycles}^>(X, \Gs \backslash\left\{ Y\right\}) \cap \Anc\left( Y, \Gs \right) = \emptyset$, then the SCG contains a structure defined in Figure \ref{subfig:CycleSansAncY:1}. We then construct an FTCG $\Gf \in \mathcal{C} \left( \Gs \right)$ that contains a collider-free backdoor path that remains in $\CD$. Figure \ref{fig:CycleSansAncY} presents this construction.

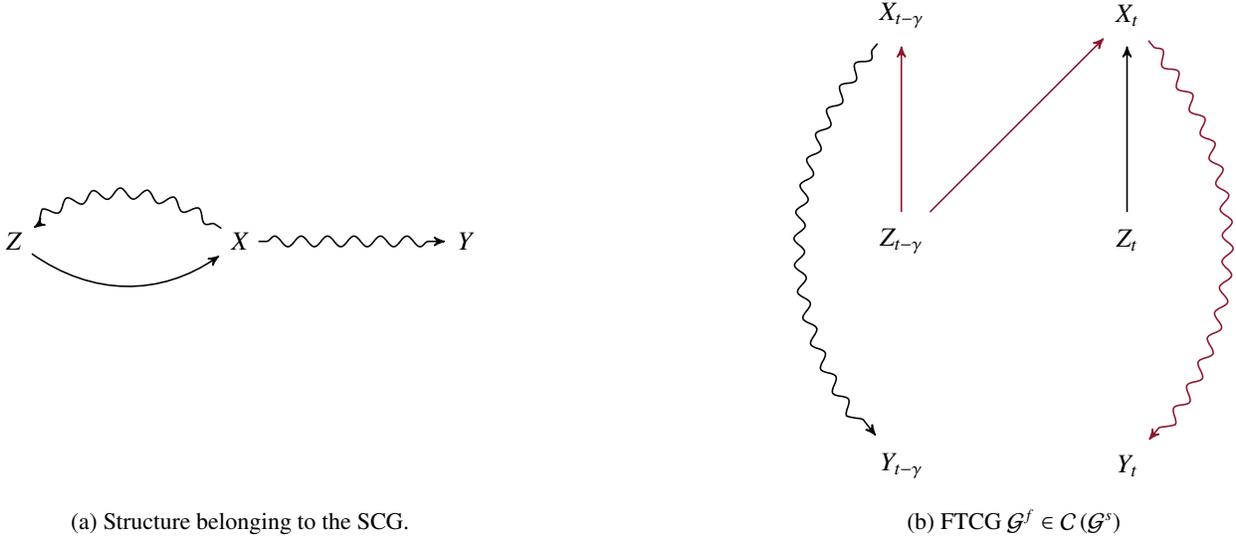
\begin{figure}[H]
    \centering
    \begin{subfigure}{0.4\textwidth}
        \centering
        \begin{tikzpicture}[->,>=stealth',shorten >=1pt,auto,node distance=3cm, semithick]
            \tikzstyle{every node}=[fill=white,draw=none,text=black, minimum height = 21pt]
            \node (X)  {$X$};
            \node (Z)  [left of=X]  {$Z$};
            \node (Y)  [right of=X] {$Y$};
            \node (Blank) [below of = X] {}; 
            
            \path (X)  edge[bend right = 35, decorate, decoration={snake, pre length = 3pt, post length=2.5pt, amplitude=2pt}] (Z)
                  (Z)  edge[bend right = 35] (X)
                  (X)  edge[decorate, decoration={snake, pre length = 3pt, post length=5pt, amplitude=2pt}] (Y);
        \end{tikzpicture}
        \caption{Structure belonging to the SCG.}
        \label{subfig:CycleSansAncY:1}
    \end{subfigure}
    \hfill
    \begin{subfigure}{0.4\textwidth}
        \centering
        \begin{tikzpicture}[->,>=stealth',shorten >=1pt,auto,node distance=3cm,
                            semithick]
            \tikzstyle{every node}=[fill=white,draw=none,text=black, minimum height = 21pt]

            \node (Xtg) {$X_{t-\gamma}$};
            \node (Ztg) [below of=Xtg] {$Z_{t-\gamma}$};
            \node (Ytg) [below of=Ztg] {$Y_{t-\gamma}$};
            \node (Xt)  [right of=Xtg] {$X_{t}$};
            \node (Zt) [below of=Xt] {$Z_{t}$};
            \node (Yt) [below of=Zt] {$Y_{t}$};

            \path (Ztg) edge[CentraleRed] (Xtg)
                  (Zt)  edge (Xt)
                  (Ztg) edge[CentraleRed] (Xt)
                  (Xtg) edge[bend right = 40,decorate, decoration={snake, pre length = 3pt, post length=5pt, amplitude=2pt}] (Ytg)
                  (Xt)  edge[CentraleRed, bend left = 40,decorate, decoration={snake, pre length = 3pt, post length=5pt, amplitude=2pt}] (Yt);
            
        \end{tikzpicture}
        \caption{FTCG $\mathcal{G}^f \in \mathcal{C} \left( \mathcal{G}^s \right)$}
        \label{subfig:CycleSansAncY:2}
    \end{subfigure}
    \caption{Illustration of Lemma \ref{lemma:gammaGeq1:CycleSansAncY}. Figure \ref{subfig:CycleSansAncY:1} represents a structure contained in the SCG. Figure \ref{subfig:CycleSansAncY:2} represents a FTCG that contains a collider-free backdoor path (in {\color{CentraleRed} red}) that does not pass through the past. The missing arrows are placed at a greater $\gamma$ and are therefore not represented.
}
    \label{fig:CycleSansAncY}
\end{figure}
\end{proof}

The combination of Lemmas \ref{lemma:gammaGeq1:CycleAvecAncY} and \ref{lemma:gammaGeq1:CycleSansAncY} allows us to prove Lemma \ref{lemma:gammaGeq1:PasCycleSurXSansY}.

\begin{restatable}{lemma}{mylemmaPasCycleSansYPourGgeqUn}{}
\label{lemma:gammaGeq1:PasCycleSurXSansY}
    Let $\Gs$ be an SCG such that $X \in \Anc\left( Y, \Gs \right)$ and $P(y_t\mid \Do(x_{t-\gamma})$ is the considered effect for $\gamma \geq 1$. If $\text{Cycles}^>\left(X, \Gs \backslash \{Y \} \right) \neq \emptyset$, then the effect is not identifiable by common backdoor.
\end{restatable}

\begin{proof}
    If  $\text{Cycles}^>(X, \Gs \backslash\left\{ Y\right\}) \neq \emptyset$, then Lemma \ref{lemma:gammaGeq1:CycleAvecAncY} or Lemma \ref{lemma:gammaGeq1:CycleSansAncY} applies. In all cases, the effect is not identifiable by common backdoor.
\end{proof}

\subsection{Proof of Theorem \ref{th:CondNecSuf_IBC_monovarie}: The case \texorpdfstring{$\gamma = 1$}{γ = 1}}

\begin{restatable}{theorem}{mythIBCmonovarieGammaUn}{}
\label{th:CondNecSuf_IBC_monovarie_gamma_un}

    Let $\Gs$ be an SCG such that $X\in \Anc(Y,\Gs)$ and we consider the total effect $P(y_t\mid \Do(x_{t-1}))$. The effect is identifiable by common backdoor in $\Gs$ if and only if:
    \begin{equation}
        \left\{
        \begin{aligned}
            & \text{Pa}(X, \Gs) \cap \Desc\left( X, \Gs \right) \subseteq \{X\},  \\
            \text{or } &
            \Gs_{\mid \Desc\left( X, \Gs \right) \cap \left( \text{Par} \left( X, \Gs \right) \cup \text{Par} \left( Y, \Gs \right) \right)} \in \left\{ X \leftrightarrows Y, \leftselfloop X \leftrightarrows Y \right\};
        \end{aligned}
        \right.
    \label{prop:CondNecSufIBGGammaUn}
    \end{equation}

    \noindent In that case, a common backdoor set is given by $\mathcal{A}_{1}$, where:
   $$
    \mathcal{A}_1 \coloneqq  \mathcal{V}^f \backslash \CD.
   $$

\end{restatable}

\subsubsection{Property \ref{prop:CondNecSufIBGGammaUn} is Necessary}
The goal of this section is to demonstrate that Property \ref{prop:CondNecSufIBGGammaUn} is necessary. The proof of this result begins with Lemma \ref{lemma:gamma1:IBCImpliqueRestrictionPaX}. 

\begin{restatable}{lemma}{mylemmaIBCImpliqueRestrictionPaXGammaUn}
\label{lemma:gamma1:IBCImpliqueRestrictionPaX}
    Let $\gamma = 1$ and $\Gs$ be an SCG such that $X \in \Anc\left( Y, \Gs \right)$. If the effect $P(y_t\mid do(x_{t-1}))$ is identifiable by common backdoor, then $\text{Pa}(X, \Gs) \cap \Desc\left( X, \Gs \right) \subseteq \{ X,Y\}$.
\end{restatable}

\begin{proof}
    We prove the contrapositive. Let $W \in \mathcal{P} \vcentcolon = \text{Pa}(X, \Gs) \cap \Desc\left( X, \Gs \right)$ such that $W \notin \{X, Y \}$. We know that $Y$ is a descendant of $X$. Let $\pi^s$ denote a directed path between $X$ and $Y$.
    \begin{itemize}
        \item If $W \in \pi^s$, then $\Gs$ contains a cycle of size at least two passing through $X$ without passing through $Y$. Lemma \ref{lemma:gammaGeq1:PasCycleSurXSansY} shows that the effect is not identifiable by common backdoor.
        \item Otherwise, there exists an FTCG containing the path $X_{t-1} \leftarrow W_{t-1} \rightarrow X_t \underset{\pi^f}{\rightsquigarrow} Y_t$, where $\pi^f$ reduces to $\pi^s$. Since $W \notin \pi^s$, the arrow $W \rightarrow X$ is not in $\pi^s$, hence $W_t \rightarrow X_t$ is not in $\pi^f$ and Assumption \ref{ass:Consistency_Time} is not violated.
    \end{itemize}
\end{proof}

Depending on the content of $ \mathcal{P} \vcentcolon = \text{Pa}(X, \Gs) \cap \Desc\left( X, \Gs \right)$, common backdoor identification is not always possible. This is especially the case if $\mathcal{P} \in \{\{Y \}, \{X, Y\} \}$. Lemma \ref{lemma:gamma1:IBCImpliqueRestrictionPaYQuandYPaX} helps to distinguish these cases.

\begin{restatable}{lemma}{mylemmaIBCImpliqueRestrictionPaYQuandYPaXGammaUn}
\label{lemma:gamma1:IBCImpliqueRestrictionPaYQuandYPaX}
    Let $\gamma = 1$ and $\Gs$ be an SCG such that $X \in \Anc\left( Y, \Gs \right)$ and $\text{Pa}(X, \Gs) \cap \Desc\left( X, \Gs \right)  \in \{\{Y \}, \{X, Y\} \}$. If the effect $P(y_t\mid do(x_{t-1}))$ is identifiable by common backdoor, then $\text{Pa}(Y, \Gs) \cap \Desc\left( X, \Gs \right) = \{ X\}$.
\end{restatable}

\begin{proof}
We know that $Y$ is a descendant of $X$, so $\text{Pa}(Y, \Gs) \cap \Desc\left( X, \Gs \right) \neq \emptyset$.

We know that $\text{Pa}(X, \Gs) \cap \Desc\left( X, \Gs \right)  \in \{\{Y \}, \{X, Y\} \}$. Hence, $Y$ is a parent of $X$. We prove by contradiction that there is no $W \in \text{Pa}(Y, \Gs) \cap \Desc\left( X, \Gs \right) $ other than $X$. Indeed, if such $W$ exists, we can construct an FTCG containing the path $X_{t-1} \leftarrow Y_{t-1} \leftarrow W_{t-1} \rightarrow Y_t$. This path cannot be blocked by a common backdoor, hence the effect would not be identifiable by common backdoor.

Since $Y$ is a descendant of $X$, $Y$ must have at least one parent in $\Gs_{\mid \Desc\left( X, \Gs \right)}$. Based on the above reasoning, the only possible parent is $X$. Therefore, $\text{Pa}(Y, \Gs) \cap \Desc\left( X, \Gs \right) = \{ X\}$.
\end{proof}

The necessary condition exhibited by Lemma \ref{lemma:gamma1:IBCImpliqueRestrictionPaYQuandYPaX} is rather complex, and a simplification approach has been proposed. This is captured by Lemma \ref{lemma:gamma1:EquivAux}.

\begin{restatable}{lemma}{mylemmaIBCEquivAuxGammaUn}
\label{lemma:gamma1:EquivAux}
    Let $\Gs$ be an SCG such that $X \in \Anc\left( Y, \Gs \right)$. Then the following statements are equivalent.

    \begin{enumerate}
        \item  $\text{Pa}(X, \Gs) \cap \Desc\left( X, \Gs \right)  \in \{\{Y \}, \{X, Y\} \}$ and $\text{Pa}(Y, \Gs) \cap \Desc\left( X, \Gs \right) = \{ X\}$
        \item $\Gs_{ \mid \Desc\left( X, \Gs \right) \cap \left( \text{Pa} \left( X, \Gs \right) \cup \text{Pa} \left( Y, \Gs \right) \right)} \in \left\{ X \leftrightarrows Y, \leftselfloop X \leftrightarrows Y \right\}$.
    \end{enumerate}
\end{restatable}

\begin{proof}
    The second formulation is simply a restatement of the graph structure restricted to the descendants of $X$ and the parents of $X$ and $Y$.
\end{proof}

We can now demonstrate that Property \ref{prop:CondNecSufIBGGammaUn} is necessary.

\begin{restatable}{lemma}{mylemmaIBCCondNecGammaUn}
\label{lemma:gamma1:CondNec}
    Let $\Gs$ be an SCG such that $X\in \Anc(Y,\Gs)$ and consider the total effect $P(y_t\mid do(x_{t-\gamma}))$ for $\gamma = 1$: $P(y_t\mid do(x_{t-1}))$. Property \ref{prop:CondNecSufIBGGammaUn} is necessary for identifiability by common backdoor.
\end{restatable}

\begin{proof}
Suppose the effect is identifiable by common backdoor. Then by Lemma \ref{lemma:gamma1:IBCImpliqueRestrictionPaX}, we know that
$$\text{Pa}(X, \Gs) \cap \Desc\left( X, \Gs \right) \subseteq \{ X,Y\}.$$
\noindent Therefore, we have
$$
    \text{Pa}(X, \Gs) \cap \Desc\left( X, \Gs \right) \subseteq \{ X\} 
    \text{    or    } 
    \text{Pa}(X, \Gs) \cap \Desc\left( X, \Gs \right)  \in \{\{Y \}, \{X, Y\} \}.
$$
\noindent According to Lemma \ref{lemma:gamma1:IBCImpliqueRestrictionPaYQuandYPaX}, we can write:
$$
    \text{Pa}(X, \Gs) \cap \Desc\left( X, \Gs \right) \subseteq \{ X\} 
    \text{    or    } 
    \left\{
        \begin{aligned}
            & \text{Pa}(X, \Gs) \cap \Desc\left( X, \Gs \right)  \in \{\{Y \}, \{X, Y\} \} \text{, and} \\
            & \text{Pa}(Y, \Gs) \cap \Desc\left( X, \Gs \right) = \{ X\}.
        \end{aligned}
    \right.
$$
\noindent Therefore, by Lemma \ref{lemma:gamma1:EquivAux}, we have:
$$
    \left\{
        \begin{aligned}
            & \text{Pa}(X, \Gs) \cap \Desc\left( X, \Gs \right) \subseteq \{X\} \text{, or} \\
            & \Gs_{ \mid \Desc\left( X, \Gs \right) \cap \left( \text{Pa} \left( X, \Gs \right) \cup \text{Pa} \left( Y, \Gs \right) \right)} \in \left\{ X \leftrightarrows Y, \leftselfloop X \leftrightarrows Y \right\}.
        \end{aligned}
    \right.
$$
\end{proof}

\subsubsection{Property \ref{prop:CondNecSufIBGGammaUn} is Sufficient.}

We show that Property \ref{prop:CondNecSufIBGGammaUn} is sufficient by enumerating all backdoor paths that remain in $\CD$. It can be observed that the conditions of Property \ref{prop:CondNecSufIBGGammaUn} and Assumption \ref{ass:Consistency_Time} prevent the construction of a collider-free backdoor path that remains in $\CD$. Then, Theorem \ref{th:equiv_IBC_multivarie} provides a common backdoor set.

\begin{restatable}{lemma}{mylemmansembleBackdoorGammaUn} Let $\Gs$ be an SCG that satisfies Property \ref{prop:CondNecSufIBGGammaUn} such that $X \in \Anc\left( Y, \Gs \right)$, and let $P(y_t\mid \Do(x_{t-1}))$ be the considered effect. Then, the set $\mathcal{A}_{1}$ defined below is a common backdoor set.
\begin{align*}
\mathcal{A}_{1} = \mathcal{V}^f \backslash \CD
\end{align*}
    \label{th:ensemble_backdoor_gamma1}
\end{restatable}

\begin{proof}
Let $\Gs$ be an SCG that satisfies Property \ref{prop:CondNecSufIBGGammaUn} such that $X \in \Anc\left( Y, \Gs \right)$, and let $P(y_t\mid \Do(x_{t-1}))$ be the considered effect. We distinguish two cases:

\begin{itemize}
    \item If $\text{Pa}(X, \Gs) \cap \Desc\left( X, \Gs \right) \subseteq \{X\}$, then there exists no collider-free backdoor path between $X_{t-1}$ and $Y_t$ that remains in $\CD$. We prove this by contradiction: suppose $\pi^f$ is such a collider-free backdoor path, and let $P_{t_p}$ be its second element. Since $P_{t_p}$ is a parent of $X_{t-1}$, we have $t_p \leq t-1$ and $P \in \text{Pa}(X, \Gs)$. Since $P_{t_p}$ belongs to $\CD$, there exists an FTCG that contains $P_{t_p} \rightsquigarrow X_{t-1}$. Thus, $t_p \geq t-1$ and $P \in \Desc(X, \Gs)$. Since $\text{Pa}(X, \Gs) \cap \Desc\left( X, \Gs \right) \subseteq \{X\}$, we conclude that $P = X$ and $P_{t_p} = X_{t-1}$, which contradicts the definition of a path.
    
    \item If $\Gs_D \coloneqq \Gs_{\mid \Desc\left( X, \Gs \right) \cap \left( \text{Par} \left( X, \Gs \right) \cup \text{Par} \left( Y, \Gs \right) \right)} \in \left\{ X \leftrightarrows Y, \leftselfloop X \leftrightarrows Y \right\}$, we prove by contradiction that no FTCG belonging to $\mathcal{C} (\Gs)$ contains a collider-free backdoor path that remains in $\CD$. Let $\Gf \in \mathcal{C} (\Gs)$ be an FTCG that contains $\pi^f$, a collider-free backdoor path that remains in $\CD$. Let $P_{t_p}$ be the second element of $\pi^f$. As seen earlier, we know that $P \in \text{Pa}(X, \Gs) \cap \Desc\left( X, \Gs \right) \subseteq \Desc\left( X, \Gs \right) \cap \left( \text{Par} \left( X, \Gs \right) \cup \text{Par} \left( Y, \Gs \right) \right)$. Thus, $P = Y$, and $\pi^f$ necessarily contains $X_{t-1} \leftarrow Y_{t-1}$. By Assumption \ref{ass:Consistency_Time}, $\Gf$ also contains the arrow $X_{t} \leftarrow Y_{t}$. Furthermore, by Lemma \ref{lemma:gamma1:EquivAux}, we know that $\text{Pa}(Y, \Gs) \cap \Desc\left( X, \Gs \right) = \{ X\}$. Therefore, for any $W \in \mathcal{V}^s$ and $t_1 \in \{ t-1, t \}$, $\pi^f$ cannot end with the arrow $W_{t_1} \rightarrow Y_t$. Hence, there exists $V \in \mathcal{V}^s$ such that $\pi^f$ ends with the arrow $V_{t} \leftarrow Y_{t}$. However, since $\pi^f$ spans from time $t-1$ to time $t$, there exist $F, U \in \mathcal{V}^s$ such that $\pi^f$ contains the arrow $F_{t-1} \rightarrow  U_{t}$. As $\pi^f$ does not contain a collider, all arrows following $F_{t-1} \rightarrow  U_{t}$ are in the direction of $\rightarrow$, which contradicts the fact that $\pi^f$ ends with $V_{t} \leftarrow Y_{t}$.
\end{itemize}

Thus, the effect is identifiable by a common backdoor, and Theorem \ref{th:equiv_IBC_multivarie} ensures that the set $\mathcal{A}_{1}$ defined below is a common backdoor set.

$$
\mathcal{A}_{1} = \mathcal{V}^f \backslash \CD.
$$

\end{proof}

\subsection{Proof of Theorem \ref{th:CondNecSuf_IBC_monovarie}: The case \texorpdfstring{$\gamma \geq 2 $}{γ >= 2}}

\begin{restatable}{theorem}{mythIBCmonovarie_gamma_geq_deux}{}
\label{th:CondNecSuf_IBC_monovarie_gamma_geq_deux}

    Let $\Gs$ be an SCG such that $X\in \Anc(Y,\Gs)$ and we consider the total effect $P(y_t\mid \Do(x_{t-\gamma}))$ for $\gamma \geq 2$. The effect is identifiable by common backdoor in $\Gs$ if and only if:
    \begin{equation}
        \label{prop:CondNecSufIBGGammaGeqDeux}
        \text{Cycles}^>\left(X, \Gs  \right) =  \emptyset.
    \end{equation}

    \noindent In that case, a common backdoor set is given by $\mathcal{A}_{\gamma}$, where:
    
   \begin{equation}
    \mathcal{A}_{\gamma} \coloneqq  \left\{ (Z_{t'})_{t' \leq t- \gamma } \mid Z \in \Anc\left(X, \Gs \right) \right\} \backslash \left\{ X_{t-\gamma} \right\}.
    \end{equation}

\end{restatable}

The proof of Theorem \ref{th:CondNecSuf_IBC_monovarie_gamma_geq_deux} is divided into the following two subsections.

\subsubsection{Property \ref{prop:CondNecSufIBGGammaGeqDeux} is necessary}

Lemma \ref{lemma:gammaGeq1:PasCycleSurXSansY} provides a necessary condition for identifiability by a common backdoor when $\gamma \geq 2$. We can also find an additional necessary condition, which is presented in Lemma \ref{lemma:gammaGeq2:PasCycleXY}.

\begin{restatable}{lemma}{mylemmaPasCycleXYPourGgeqDeux}
\label{lemma:gammaGeq2:PasCycleXY}
    Let $\Gs$ be an SCG such that $X \in \Anc\left( Y, \Gs \right)$ and $P(y_t\mid \Do(x_{t-\gamma})$ is the considered effect for $\gamma \geq 2$. If $\Gs$ contains a cycle that passes through $X$ and $Y$, then the effect is not IBC.
\end{restatable}

\usetikzlibrary {decorations.shapes}
\usetikzlibrary {decorations.text}
\begin{figure}[H]
    \centering
    \begin{subfigure}{0.39\textwidth}
        \centering
        \begin{tikzpicture}[->,>=stealth',shorten >=1pt,auto,node distance=3cm, semithick]
            \tikzstyle{every node}=[fill=white,draw=none,text=black, minimum height = 21pt]
            \node (X)  {$X$};
            \node (Y)  [right of=X] {$Y$};
            \node (Label1) at (1.5,-1) {$\pi_1$};
            \node (Label2) at (1.5,1) {$\pi_2$};
            
            \draw (X) edge[bend right = 35, decorate, decoration={snake, pre length = 3pt, post length=2.5pt, amplitude=2pt}] (Y);

            \draw (Y) edge[bend right = 35] (X);
        \end{tikzpicture}
        \caption{Structure belonging to the SCG.}
        \label{subfig:CycleXY:1}
    \end{subfigure}
    \hfill
    \begin{subfigure}{0.5\textwidth}
        \centering
        \begin{tikzpicture}[scale = 1.5, ->,>=stealth',shorten >=1pt,auto,node distance=2cm, semithick]
            \tikzstyle{every node}=[fill=white,draw=none,text=black, minimum height = 21pt]

            \node (Xtg) {$X_{t-\gamma}$};
            \node (Xtg1) [right of =Xtg] {$X_{t-\gamma + 1}$};
            \node (Dots1) at (2.2, 0) {$\dots$};
            \node (Xt) at (2.7,0) {$X_t$};
            \node (Ytg) [below of=Xtg] {$Y_{t-\gamma}$};
            \node (Ytg1) [right of =Ytg] {$Y_{t-\gamma + 1}$};
            \node (Dots2) [below of =Dots1] {$\dots$};
            \node (Yt) [below of =Xt] {$Y_t$};

            \path  (Ytg) edge[CentraleRed] (Xtg);
            \path  (Ytg) edge[CentraleRed] (Xtg1);
            \path  (Xtg1) edge[CentraleRed, decorate, decoration={snake, pre length = 3pt, post length=2.5pt, amplitude=2pt}] (Yt);

        \end{tikzpicture}
        \caption{FTCG $\mathcal{G}^f \in \mathcal{C} \left( \mathcal{G}^s \right)$}
        \label{subfig:CycleXY:2}
    \end{subfigure}

        \begin{subfigure}{0.39\textwidth}
        \centering
        \begin{tikzpicture}[->,>=stealth',shorten >=1pt,auto,node distance=3cm, semithick]
            \tikzstyle{every node}=[fill=white,draw=none,text=black, minimum height = 21pt]
            \node (X)  {$X$};
            \node (Y)  at (3,0) {$Y$};
            \node (S) at (1,0.65) {$S$};
            \node (Label1) at (1.5,-1) {$\pi_1$};

            \draw (X) edge[bend right = 35, decorate, decoration={snake, pre length = 3pt, post length=2.5pt, amplitude=2pt}] (Y);
            \draw (S) edge[bend right = 20] (X);
            
            \draw (Y) edge[bend right = 20, decorate, decoration={snake, pre length = 3pt, post length=2.5pt, amplitude=2pt}] (S);
        \end{tikzpicture}
        \caption{Structure belonging to the SCG.}
        \label{subfig:CycleXY:3}
    \end{subfigure}
    \hfill
    \begin{subfigure}{0.5\textwidth}
        \centering
        \begin{tikzpicture}[scale = 1.5, ->,>=stealth',shorten >=1pt,auto,node distance=2cm, semithick]
            \tikzstyle{every node}=[fill=white,draw=none,text=black, minimum height = 21pt]

            \node (Xtg) {$X_{t-1}$};
            \node (Xtg1) [right of =Xtg] {$X_{t}$};
            \node (Stg) [below of=Xtg] {$S_{t-\gamma}$};
            \node (Stg1) [right of =Stg] {$S_{t-\gamma + 1}$};
            \node (Dots2) [below of =Dots1] {};
            \node (Dots3) [below of =Dots2] {$\dots$};
            \node (St) [below of =Xt] {$S_t$};
            \node (Yt) [below of = St] {$Y_t$};
            \node (Ytg) [below of = Stg] {$Y_{t- \gamma}$};
            \node (Ytg1) [below of = Stg1] {$Y_{t - \gamma + 1}$};

            \path  (Stg) edge[CentraleRed] (Xtg);
            \path  (Stg1) edge (Xtg1);
            \path  (St) edge (Xt);
            \path  (Stg) edge[CentraleRed] (Xtg1);
            \path  (Xtg1) edge[CentraleRed,bend left = 6, decorate, decoration={snake, pre length = 3pt, post length=2.5pt, amplitude=2pt}] (Yt);

        \end{tikzpicture}
        \caption{FTCG $\mathcal{G}^f \in \mathcal{C} \left( \mathcal{G}^s \right)$}
        \label{subfig:CycleXY:4}
    \end{subfigure}
    \caption{Illustration of Lemma \ref{lemma:gammaGeq2:PasCycleXY}. Figures \ref{subfig:CycleXY:1} and \ref{subfig:CycleXY:3} depict structures contained within the SCG. For clarity, they distinguish the cases based on the length of $\pi_2$. Figures \ref{subfig:CycleXY:2} and \ref{subfig:CycleXY:4} respectively represent an FTCG containing a collider-free backdoor path (in {\color{CentraleRed} red}) that remains within $\CD$. The missing arrows are positioned at a greater $\gamma$ and are thus not depicted. Figures \ref{subfig:CycleXY:2} and \ref{subfig:CycleXY:4} are indeed constructible since $\pi_1$ and $\pi_2$ only share their endpoints.
}
    \label{fig:CycleXY}
\end{figure}
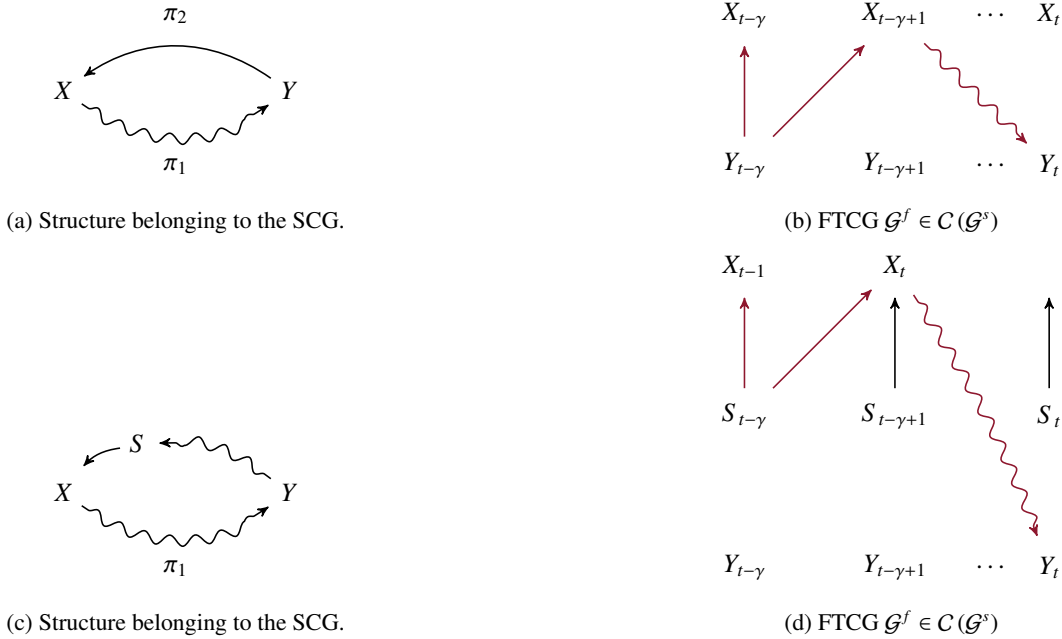

\begin{proof}
Suppose that $\Gs$ contains a cycle passing through $X$ and $Y$. We denote by $\pi_1$ the directed path contained in the cycle from $X$ to $Y$ and by $\pi_2$ the path from $Y$ to $X$. There are two cases:
\begin{itemize}
    \item Either $\pi_1 \cap \pi_2 \neq \{ X, Y \}$. In this case, $\Gs$ contains a cycle of length at least two that passes through $X$. Thus, Lemma \ref{lemma:gammaGeq1:PasCycleSurXSansY} shows that the effect is not identifiable by a common backdoor.
    \item Otherwise, $\pi_1$ and $\pi_2$ only share their endpoints. We can then consider the constructions shown in Figure \ref{fig:CycleXY}, which demonstrates that the effect is not identifiable by a common backdoor.
\end{itemize}
\end{proof}

\begin{restatable}{lemma}{mylemmaPasCyclePourGgeqDeux}
    Let $\Gs$ be an SCG such that $X \in \Anc\left( Y, \Gs \right)$ and $P(y_t\mid \Do(x_{t-\gamma})$ is the considered effect for $\gamma \geq 2$ If $\Gs$ does not satisfy Property \ref{prop:CondNecSufIBGGammaGeqDeux}, then $\Gs$ is not identifiable by a common backdoor.
\end{restatable}

\begin{proof}
    If $\Gs$ does not satisfy Property \ref{prop:CondNecSufIBGGammaGeqDeux}, then Lemma \ref{lemma:gammaGeq1:PasCycleSurXSansY} or Lemma \ref{lemma:gammaGeq2:PasCycleXY} applies. In either case, $\Gs$ is not identifiable by a common backdoor.
\end{proof}

\subsubsection{Property \ref{prop:CondNecSufIBGGammaGeqDeux} is sufficient}

\begin{restatable}{lemma}{mylemmansembleBackdoorGammaGeqdeux}
Let $\Gs$ be an SCG that satisfies Property \ref{prop:CondNecSufIBGGammaGeqDeux}. Then, for $\gamma \geq 2$, the set $\mathcal{A}_{\gamma}$ defined below is a common backdoor set.

\begin{align*}
\mathcal{A}_{\gamma} \vcentcolon= \left\{ (Z_{t'})_{t' \leq t- \gamma } \mid Z \in \Anc\left(X, \Gs \right) \right\} \backslash \left\{ X_{t-\gamma} \right\}.
\end{align*}

    \label{lemma:ensemble_backdoor_gamma_geq_2}
\end{restatable}

\begin{proof}
    $\mathcal{A}_{\gamma}$ contains all the parents of $X_{t-\gamma}$ in any FTCG $\Gf \in \mathcal{C}(\Gs)$. Therefore, $\mathcal{A}_{\gamma}$ blocks all backdoor paths in any FTCG $\Gf \in \mathcal{C}(\Gs).$

    We now show by contradiction that $\mathcal{A}_{\gamma}$ does not contain any descendants of $X_{t-\gamma}$ in any FTCG $\Gf \in \mathcal{C}(\Gs)$. Let $Z_{t'}$ be such a vertex. Then $Z$ is a descendant of $X$ in $\Gs$. However, by the definition of $\mathcal{A}_{\gamma}$, $Z$ is also an ancestor of $X$. Thus, $\Gs$ would contain a cycle of length $\geq 2$, which contradicts Property \ref{prop:CondNecSufIBGGammaGeqDeux}.
\end{proof}

\subsection{Efficient testing for collider-free-backdoor path with a fork. \label{ssct:D2}}
It is possible to compute 
$
\left( \max_{X^i_{t-\gamma_i}} \left\{t^{\NC}_{X^i_{t-\gamma_i}}(S) \right\}\right)_{S \in \mathcal{V}^s}
$
in a single traversal of the graph \(\Gs\). To achieve this, one only needs to modify the initialization of the priority queue \(Q\) in Algorithm \ref{algo:calcul_V_E_acc}. In this optimized strategy, \(Q\) is initialized with 
$
\bigcup_{i \in \{1,\dots, n\}} \left\{P_{t-\gamma_i} \mid P \in \text{Pa}(X^i, \Gs) \right\},
$
as this covers all potential parents of interventions at time \(t - \gamma_i\). This is equivalent to computing 
$
(t^{\NC}_{M_{\star}}(S))_{S \in \mathcal{V}^s},
$
where \(M_{\star}\) is a fictitious vertex added to every FTCG in \(\C(\Gs)\), ensuring that each parent of an intervention also becomes a parent of \(M_{\star}\). The complexity remains in 
\(
\mathcal{O} \left(\left| \mathcal{E}^s \right| + \left| \mathcal{V}^s \right| \log \left| \mathcal{V}^s \right|\right),
\)
as the single traversal avoids redundant checks.

During this computation, it is also possible to efficiently test for the existence of collider-free backdoor paths involving a fork that remains within \(\CD\). Since 
\(
(t^{\NC}_{Y_{t}}(S))_{S \in \mathcal{V}^s}
\)
is already computed, the algorithm can, for each parent \(P\) encountered, directly verify whether 
\(
t_{\NC}(P) \leq t^{\NC}_{M_{\star}}(P) = \max_{X^i_{t-\gamma_i}} \left\{t^{\NC}_{X^i_{t-\gamma_i}}(P)\right\}
\)
and whether 
\(
t_{\NC}(P) \leq t^{\NC}_{Y_t}(F).
\)
This characterizes the existence of a collider-free backdoor path with a fork remaining in \(\CD\) from any intervention to \(Y_t\). Thus, the need for a separate loop to test the existence of backdoor paths with forks is eliminated, further optimizing the computation.

\subsection{A pseudo-linear algorithm for 
 multivaiate IBC.}

By utilizing insights from the previous subsections, it is possible to test for the existence of a directed path in $\mathcal{O} \left(\left| \mathcal{E}^s \right| + \left| \mathcal{V}^s \right| \log \left| \mathcal{V}^s \right| \right)$, as well as for the existence of a backdoor path with a fork in the same time complexity. Therefore, IBC can be computed in $\mathcal{O} \left(n \log n + \left| \mathcal{E}^s \right| + \left| \mathcal{V}^s \right| \log \left| \mathcal{V}^s \right| \right)$.

\section{Theorem \ref{th:CondNecSuf_IBC_monovarie} without Assumption \ref{ass:Consistency_Time}}
\label{sec:app:mono}
When Assumption \ref{ass:Consistency_Time} is not assumed, we lose some power in identification. In the monovariate setting, this is exemplified by Theorem \ref{th:CondNecSuf_IBC_monovarie_whithout_consistency}.
\begin{restatable}{theorem}{mythIBCmonovariewithoutconsistency}{}
\label{th:CondNecSuf_IBC_monovarie_whithout_consistency}

    Let $\Gs$ be an SCG such that $X \in \Anc(Y,\Gs)$ and we consider the total effect $P(y_t\mid \Do(x_{t-\gamma}))$. The effect is identifiable by common backdoor in $\Gs$ if and only if:
\begin{itemize}
     \item For $\gamma = 0$, there is no collider-free backdoor path from $X$ to $Y$ in $\Gs_{\mid \Desc(X, \Gs)}$.

    \item For $\gamma \geq 1$,
   $\text{Cycles}^>\left(X, \Gs  \right) =  \emptyset$.
\end{itemize}
    \noindent In theses cases, a common backdoor set is given by $\mathcal{A}_{\gamma}$, where:
   \begin{align*}
       \mathcal{A}_0 &\coloneqq \bigcup_{\pi^s \text{ backdoor}} \Bigl\{ Z_t  \mid Z \in  \pi^s \backslash \Desc\left( X,\Gs \right) \Bigr\} \cup \Bigl\{ \left( Z_{t'} \right)_{t' < t} \mid Z \in \Gs \Bigr\}\\
    \mathcal{A}_{\gamma} &\coloneqq  \left\{ (Z_{t'})_{t' \leq t - \gamma} \mid Z \in \Anc\left(X, \Gs \right) \right\} \backslash \left\{ X_{t-\gamma} \right\}, \text{for }\gamma \geq 1.
   \end{align*}

\end{restatable}

\begin{proof}

Only the sufficiency part of the case $\gamma = 1$ of Theorem \ref{th:CondNecSuf_IBC_monovarie} necessitates Assumption \ref{ass:Consistency_Time} in its proof. Therefore, we concentrate only on the case $\gamma = 1$ to prove Theorem \ref{th:CondNecSuf_IBC_monovarie_whithout_consistency}. 

Lemma \ref{lemma:gammaGeq1:PasCycleSurXSansY} holds i.e. $\text{Cycles}^>\left(X, \Gs \backslash \{Y \} \right) \neq \emptyset$ is a necessary condition and Lemma \ref{lemma:gammaGeq1:PasCycleXY:sans_consistency} shows that is is necessary that $\Gs$ does not contain any directed cycle that passes through $X$ and $Y$. Therefore $\text{Cycles}^>\left(X, \Gs  \right) =  \emptyset$ is a necessary condition. Since Lemma \ref{lemma:ensemble_backdoor_gamma_geq_2} does not need Assumption \ref{ass:Consistency_Time} to hold true, $\text{Cycles}^>\left(X, \Gs  \right) =  \emptyset$ is also a sufficient condition and $\mathcal{A}_{\gamma}$ is a valid common backdoor set.
\end{proof}

\begin{restatable}{lemma}{mylemmaPasCycleXYPourGUn}
\label{lemma:gammaGeq1:PasCycleXY:sans_consistency}
    Let $\Gs$ be an SCG such that $X \in \Anc\left( Y, \Gs \right)$ and $P(y_t\mid \Do(x_{t-1})$ is the considered effect. Assumption \ref{ass:Consistency_Time} is not satisfied by $\Gs$. If $\Gs$ contains a cycle that passes through $X$ and $Y$, then the effect is not IBC.
\end{restatable}

\usetikzlibrary {decorations.shapes}
\usetikzlibrary {decorations.text}
\begin{figure}[H]
    \centering
    \begin{subfigure}{0.39\textwidth}
        \centering
        \begin{tikzpicture}[->,>=stealth',shorten >=1pt,auto,node distance=3cm, semithick]
            \tikzstyle{every node}=[fill=white,draw=none,text=black, minimum height = 21pt]
            \node (X)  {$X$};
            \node (Y)  [right of=X] {$Y$};
            \node (Label1) at (1.5,-1) {$\pi_1$};
            \node (Label2) at (1.5,1) {$\pi_2$};
            
            \draw (X) edge[bend right = 35, decorate, decoration={snake, pre length = 3pt, post length=2.5pt, amplitude=2pt}] (Y);

            \draw (Y) edge[bend right = 35] (X);
        \end{tikzpicture}
        \caption{Structure belonging to the SCG.}
        \label{subfig:CycleXYsansConsistency:1}
    \end{subfigure}
    \hfill
    \begin{subfigure}{0.5\textwidth}
        \centering
        \begin{tikzpicture}[scale = 1.5, ->,>=stealth',shorten >=1pt,auto,node distance=2cm, semithick]
            \tikzstyle{every node}=[fill=white,draw=none,text=black, minimum height = 21pt]

            \node (Xtg) {$X_{t-1}$};
            \node (Xtg1) [right of =Xtg] {$X_{t}$};
            \node (Ytg) [below of=Xtg] {$Y_{t-1}$};
            \node (Ytg1) [right of =Ytg] {$Y_{t}$};
            
            \path  (Ytg)  edge[CentraleRed] (Xtg);
            \path  (Xtg1) edge[CentraleRed] (Ytg1);
            \path  (Ytg)  edge[CentraleRed] (Xtg1);

        \end{tikzpicture}
        \caption{FTCG $\mathcal{G}^f \in \mathcal{C} \left( \mathcal{G}^s \right)$}
        \label{subfig:CycleXYsansConsistency:2}
    \end{subfigure}

        \begin{subfigure}{0.39\textwidth}
        \centering
        \begin{tikzpicture}[->,>=stealth',shorten >=1pt,auto,node distance=3cm, semithick]
            \tikzstyle{every node}=[fill=white,draw=none,text=black, minimum height = 21pt]
            \node (X)  {$X$};
            \node (Y)  at (3,0) {$Y$};
            \node (S) at (1,0.65) {$S$};
            \node (Label1) at (1.5,-1) {$\pi_1$};

            \draw (X) edge[bend right = 35, decorate, decoration={snake, pre length = 3pt, post length=2.5pt, amplitude=2pt}] (Y);
            \draw (S) edge[bend right = 20] (X);
            
            \draw (Y) edge[bend right = 20, decorate, decoration={snake, pre length = 3pt, post length=2.5pt, amplitude=2pt}] (S);
        \end{tikzpicture}
        \caption{Structure belonging to the SCG.}
        \label{subfig:CycleXYsansConsistency:3}
    \end{subfigure}
    \hfill
    \begin{subfigure}{0.5\textwidth}
        \centering
        \begin{tikzpicture}[scale = 1.5, ->,>=stealth',shorten >=1pt,auto,node distance=2cm, semithick]
            \tikzstyle{every node}=[fill=white,draw=none,text=black, minimum height = 21pt]

            \node (Xtg) {$X_{t-1}$};
            \node (Xtg1) [right of =Xtg] {$X_{t}$};
            \node (Stg) [below of=Xtg] {$S_{t-1}$};
            \node (Stg1) [right of =Stg] {$S_{t}$};
            \node (Ytg) [below of = Stg] {$Y_{t-1}$};
            \node (Ytg1) [below of = Stg1] {$Y_{t}$};

            \path  (Stg) edge[CentraleRed] (Xtg);
            \path  (Stg) edge[CentraleRed] (Xtg1);
            \path  (Xtg1) edge[CentraleRed,bend left = 30, decorate, decoration={snake, pre length = 3pt, post length=2.5pt, amplitude=2pt}] (Ytg1);

        \end{tikzpicture}
        \caption{FTCG $\mathcal{G}^f \in \mathcal{C} \left( \mathcal{G}^s \right)$}
        \label{subfig:CycleXYsansConsistency:4}
    \end{subfigure}
    \caption{Illustration of Lemma \ref{lemma:gammaGeq1:PasCycleXY:sans_consistency}. Figures \ref{subfig:CycleXYsansConsistency:1} and \ref{subfig:CycleXYsansConsistency:3} depict structures contained within the SCG. For clarity, they distinguish the cases based on the length of $\pi_2$. Figures \ref{subfig:CycleXYsansConsistency:2} and \ref{subfig:CycleXYsansConsistency:4} respectively represent an FTCG containing a collider-free backdoor path (in {\color{CentraleRed} red}) that remains within $\CD$. The missing arrows are positioned at a greater $\gamma$ and are thus not depicted. Figures \ref{subfig:CycleXYsansConsistency:2} and \ref{subfig:CycleXYsansConsistency:4} are indeed constructible since $\pi_1$ and $\pi_2$ only share their endpoints and because Assumption \ref{ass:Consistency_Time} does not need to be satisfied.
}
    \label{fig:CycleXYsansConsistency}
\end{figure}
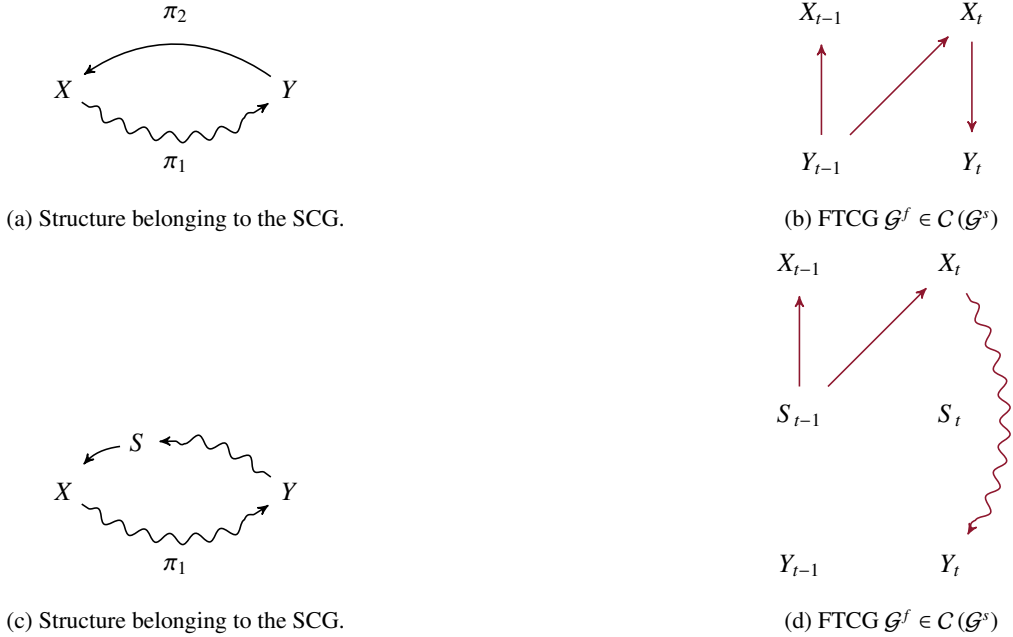

\begin{proof}
Suppose that $\Gs$ contains a cycle passing through $X$ and $Y$. We denote by $\pi_1$ the directed path contained in the cycle from $X$ to $Y$ and by $\pi_2$ the path from $Y$ to $X$. There are two cases:
\begin{itemize}
    \item Either $\pi_1 \cap \pi_2 \neq \{ X, Y \}$. In this case, $\Gs$ contains a cycle of length at least two that passes through $X$. Thus, Lemma \ref{lemma:gammaGeq1:PasCycleSurXSansY} shows that the effect is not identifiable by a common backdoor.
    \item Otherwise, $\pi_1$ and $\pi_2$ only share their endpoints. We can then consider the construction shown in Figure \ref{fig:CycleXYsansConsistency}, which demonstrates that the effect is not identifiable by a common backdoor.
\end{itemize}
\end{proof}

\section{An efficient algorithm for IBC under Assumption \ref{ass:Consistency_Time}. \label{sct:IBC_pseudo_lineaire}}

The algorithm presented in Algorithm \ref{algo:calcul_IBC_with_consistency} aims to determine whether the causal effect $P(y_t \mid \text{do}(x^1_{t-\gamma_1}), \dots, \text{do}(x^n_{t-\gamma_n}))$ is identifiable by common backdoor under Assumption \ref{ass:Consistency_Time}. Its correctness arises from Theorem \ref{th:equiv_IBC_multivarie}. Algorithm \ref{algo:calcul_IBC_with_consistency} tests the existence of an FTCG which contains a collider-free backdoor path from an intervention to $Y_t$ which remains in $\CD$. It starts by checking the existence of a directed path and then test the existence of backdoor paths with a fork using the characterisation of Lemma \ref{lemma:equiv_existence_chemin_fork_NC}. 

The algorithm focuses on checking the existence of backdoor path with a fork. To do so, it uses the characterisation of Lemma \ref{lemma:equiv_existence_chemin_fork_NC}. In line 1, it checks the backdoor path whose fork reduces to $F \neq Y$. Thanks to the strategy presented in Subsection \ref{ssct:D2}, this can be done in a single traversing. Then, on line 2, the algorithm checks that last case of case \ref{lemma:equiv_existence_chemin_fork_NC:2a}. It works on path whose fork reduces to $Y$ and whose interventions are at time $t- \gamma_i \neq t_{\NC}(Y)$. Thanks to the strategy presented in Subsection \ref{ssct:D2}, this can be done in a single traversing. Then, the algorithm checks the cases \ref{lemma:equiv_existence_chemin_fork_NC:2b}. It starts by the case \ref{lemma:equiv_existence_chemin_fork_NC:2b:i} on line 3. Since $\forall X^i_{t-\gamma_i} \in \mathcal{X}, t-\gamma_i = t_{\NC}(Y)$, finding a directed path from $Y_{t_{\NC}(Y)}$ to $X^i_{t_{\NC}(Y)}$ without using $X^i_{t_{\NC}(Y)} \leftarrow Y_{t_{\NC}(Y)} $ is equivalent of finding a directed path from $Y$ to $X^i$ in $\Gs$ without using $Y \rightarrow X^i$. This can be done by a BFS algorithm in $\Gs$. The algorithm finishes by checking the case \ref{lemma:equiv_existence_chemin_fork_NC:2b:ii}. $\mathcal{X}'$ represents the set of interventions at time $t_{\NC}(Y)$ for which $Y_{t_{\NC}(Y)}$ is $X^i_{t-\gamma_i}$-$\NC$-accessible. If $\left| \mathcal{X}' \right| \geq 2$, since $Y_{t_{\NC}(Y)}$ is $Y_t$-$\NC$-accessible, there is at least one $X^i_{t-\gamma_i} \in \mathcal{X}'$ that $Y_{t_{\NC}(Y)}$ is $Y_t$ - $\NC$-accessible without using $X^i_t \rightarrow Y_t$. Indeed,  $Y_{t_{\NC}(Y)}$ is $Y_t$-$\NC$-accessible, therefore, there exists a FTCG which contains a directed path from $Y_{t_{\NC}(Y)}$ to $Y_t$. This path cannot use simultaneously an arrow $X^i_t \rightarrow Y_t$ and an arrow $X^j_t \rightarrow Y_t$.  This reasoning explains the third line of the algorithm. The line 5 is reached when $\left| \mathcal{X}' \right| = 1$. In this case, a traversing algorithm starting from $Y_{t_{\NC}(Y)}$ can test if $Y_{t_{\NC}(Y)}$ is $Y_t$ - $\NC$-  accessible without using $X^i_t \rightarrow Y_t$. All the subcases of case  \ref{lemma:equiv_existence_chemin_fork_NC:2b:ii} have been checked by the algorithm

Algorithm \ref{algo:calcul_IBC_with_consistency} calls 5 traversing of $\Gs$. Therefore, its complexity is indeed pseudo-linear.

\begin{algorithm}[t]

\caption{Computation of identifiability by common backdoor under Assumption \ref{ass:Consistency_Time}.}
\label{algo:calcul_IBC_with_consistency}
\SetNlSty{}{\relsize{2.5}}{.}
\SetKwInput{KwData}{Input}
\SetKwInput{KwResult}{Output}

\KwData{ $\Gs$ an SCG and $P(y_t \mid \text{do}(x^1_{t-\gamma_1}), \dots, \text{do}(x^n_{t-\gamma_n}))$ the considered effect.}
\KwResult{A boolean indicating whether the effect is identifiable by common backdoor or not.}

\tcp{Test all all forks $F \neq Y$ using strategy from \ref{ssct:D2}:}

\nl \If(\tcc*[f]{\ref{ssct:D2}}){there exist $X^i_{t- \gamma_i}$ and $F$ such that $F_{t_{\NC}(F)}$ is well defined, $X^i_{t- \gamma_i}$-$\NC$-accessible and $Y_t$-$\NC$-accessible}{
    \Return{False}
}

\tcp{Test for $F = Y$:}
\If{$Y_{t_{\NC}(Y)}$ is $Y_t$-$\NC$-accessible}{
    \tcp{Last case of \ref{lemma:equiv_existence_chemin_fork_NC:2a}}
    
    \nl \If(\tcc*[f]{\ref{ssct:D2}}){there exist $X^i_{t- \gamma_i}$ such that $t- \gamma_i \neq t_{\NC}(Y)$ and $Y_{t_{\NC}(Y)}$ is $X^i_{t- \gamma_i}$-$\NC$-accessible}{ 
        \Return{False}
    }
    
    \tcp{Test for \ref{lemma:equiv_existence_chemin_fork_NC:2b}:}
    
    $\mathcal{X} \gets \{X^i_{t-\gamma_i} \mid t-\gamma_i = t_{\NC}(Y) \}$ \;
    
    \tcp{Case \ref{lemma:equiv_existence_chemin_fork_NC:2b:i}:}
    
    \nl \If(\tcc*[f]{BFS}){ $\exists~ X^i_{t-\gamma_i} \in \mathcal{X}$ such that $Y_{t_{\NC}(Y)}$ is $X^i_{t- \gamma_i}$ - $\NC$-accessible without using $X^i_{t- \gamma_i} \leftarrow Y_{t- \gamma_i}$ } {
        \Return{False}
    }
    \tcp{Case \ref{lemma:equiv_existence_chemin_fork_NC:2b:ii}:}
    $\mathcal{X}' \gets \{X^i_{t-\gamma_i} \mid t-\gamma_i = t_{\NC}(Y) \text{ and } X^i \in \Desc(Y, \Gs) \}$ \tcc*{BFS}

    \nl \If{$\left| \mathcal{X}' \right| \geq 2$} {
        \Return{False}
    } 
    \nl \If(\tcc*[f]{BFS}){$\mathcal{X}' = \{ X^i_{t-\gamma_i}\}$ \textbf{and} $Y_{t_{\NC}(Y)}$ is $Y_t$ - $\NC$-  accessible without using $X^i_t \rightarrow Y_t$.} {
        \Return{False}
        
    }
}
\Return{True}
\end{algorithm}

\section{Proofs of Section \ref{sec:notions}}

\mypropbackdooradjustmentmultivarie*

\begin{proof}
    We prove Proposition \ref{prop:backdoor_adjustment_multivarie} using do-calculus. We begin by marginalizing over $\mathcal{Z}$. Let us denote symbolically $\textbf{\text{do}(x)} \vcentcolon= \text{do}(x^1), \dots, \text{do}(x^n)$ and $\textbf{Y} = \left( Y^1, \dots , Y^m\right)$.

    $$P\left(y^1, \dots, y^m \mid \textbf{\text{do}(x)} \right)
= \sum_{z\in \mathcal{Z}} P\left(y^1, \dots, y^m  \mid \textbf{\text{do}(x)}, z\right) P(z \mid \textbf{\text{do}(x)}).$$

\noindent $\mathcal{Z}$ contains no descendant of any $X^i$, hence $\forall i,   ~ \mathcal{Z} \ind_{\mathcal{G}_{\overline{X^i}}} X^i$. We apply \footnote{We may apply the rules in any order. In any case, $\mathcal{Z} \ind_{\mathcal{G}_{\overline{X^i}}} X^i$ implies the desired graphical independence.} \textbf{rule 3} of do-calculus n times and obtain:

    $$P\left(y^1, \dots, y^m \mid \textbf{\text{do}(x)} \right)
= \sum_z P\left(y^1, \dots, y^m  \mid \textbf{\text{do}(x)}, z\right) P(z)$$

\noindent Now, for all $ i \in \left[ \left| 1,n \right| \right]$, we have $\textbf{Y} \ind_{\mathcal{G}_{\overline{X^{i+1}}, \dots, \overline{X^n} \underline{X^i}}} X^i \mid \mathcal{Z}, \{ X^j\}_{j\neq i}$. Indeed, let us proceed by contradiction. If this is false, then there exists a backdoor path $\pi$ from $X^i$ to $Y$ that is not blocked by $\mathcal{Z} \cup \{ X^j\}_{j\neq i}$ in $\mathcal{G}_{\overline{X^{i+1}}, \dots, \overline{X^n}}$. Two cases arise:

\begin{itemize}
    \item If $\pi$ does not pass through $\{X^1, \dots, X^{i-1} \}$, then $\pi$ would not be blocked by $\mathcal{Z}$ in $\mathcal{G}_i$ \footnote{Note that $\pi$ indeed exists in $\mathcal{G}_i$ since it does not pass through $\{X^1, \dots, X^{i-1} \}$.}. This contradicts the fact that $\mathcal{Z}$ satisfies the multivariate backdoor criterion.
    \item Otherwise, $\pi$ passes through $\{X^1, \dots, X^{i-1} \}$. Necessarily, all encounters between $\pi$ and $\{X^1, \dots, X^{i-1} \}$ form colliders, because otherwise, $\pi$ would be blocked by conditioning on $\{ X^j\}_{j\neq i}$. We then consider $X^{i_0}$, the last element of $\{X^1, \dots, X^{i-1} \}$ visited by $\pi$, and $\pi'$, the sub-path of $\pi$ from $X^{i_0}$ to $\textbf{Y}$. $\pi'$ is a backdoor path from $X^{i_0}$ to $\textbf{Y}$, unblocked by $\mathcal{Z} \cup \{ X^j\}_{j\neq i} \backslash \{ X^{i_0} \}$ in $\mathcal{G}_{\overline{X^{i+1}}, \dots, \overline{X^n}}$. Moreover, $\pi'$ does not pass through $\{X^1, \dots, X^{i} \} \backslash \{ X^{i_0} \}$ \footnote{Since $X^{i_0}$ is the last node visited by $\pi$ in $\{X^1, \dots, X^{i-1} \}$, and $\pi$ is a path, it does not pass through a same vertex twice. }. Therefore, $\pi'$ is not blocked by $\mathcal{Z} \cup \{X^{j}\}_{j \neq i_0}$ in $\mathcal{G}_{ \overline{(X^{j})_{j \neq i}}}$, which contradicts the fact that $\mathcal{Z}$ satisfies the multivariate backdoor criterion.
\end{itemize}
We can thus apply \textbf{rule 2} of do-calculus n times, obtaining:

$$P\left(y^1, \dots, y^m \mid \text{do}(x^1), \dots, \text{do}(x^n) \right)
= \sum_z P\left(y^1, \dots, y^m  \mid x^1, \dots, x^n, z\right) P(z).$$
\end{proof}

\mylemmaobservationsmultiples*

\begin{proof}
    We proceed directly by equivalence for the two conditions of the multivariate backdoor criterion.
    \begin{itemize}
        \item Condition \ref{def:critere_backdoor:1}: The condition does not depend on $j$ so the equivalence holds.

        \item Condition \ref{def:critere_backdoor:2}: $\forall (X^i, Y^j), ~~ X^i \ind_{\Gf_i} Y^j \mid \mathcal{Z} \Leftrightarrow  \forall Y^j, ~~ \left( \forall X^i, ~~ X^i \ind_{\Gf_i} Y^j \mid \mathcal{Z} \right)$.
    \end{itemize}
\end{proof}

\end{document}